\documentclass[11 pt,onecolumn]{IEEEtran}

\IEEEoverridecommandlockouts     

\usepackage{graphicx,algorithm,algorithmic,subfigure,latexsym,amsmath,amsfonts,amssymb,cite,mathrsfs}
\usepackage{epsfig,color}
\usepackage{psfig,psfrag}
\usepackage{pstricks}
\usepackage{epstopdf}

\newtheorem{lemma}{Lemma}[section]
\newtheorem{theorem}{Theorem}[section]
\newtheorem{corollary}{Corollary}[section]
\newtheorem{remark}{Remark}[section]
\newtheorem{definition}{Definition}[section]

\def\rrr#1\\{\par
\medskip\hbox{\vbox{\parindent=2em\hsize=6.12in
\hangindent=4em\hangafter=1#1}}}
\allowdisplaybreaks

\begin{document}



%

\title{\LARGE Convergence Analysis and Parallel Computing Implementation for the Multiagent Coordination Optimization Algorithm}

\author{Qing Hui and Haopeng Zhang\\
\small{Control Science and Engineering Laboratory\\
Department of Mechanical Engineering \\
Texas Tech University \\
Lubbock, TX 79409-1021\\
{\tt\small (qing.hui@ttu.edu; haopeng.zhang@ttu.edu)}\\
Technical Report CSEL-06-13, June 2013}
\thanks{This work was supported by the Defense Threat Reduction Agency, Basic Research Award \#HDTRA1-10-1-0090 and Fundamental Research Award \#HDTRA1-13-1-0048, to Texas Tech University.}
}

\maketitle
\IEEEpeerreviewmaketitle
\pagestyle{empty}
\thispagestyle{empty}
\baselineskip 22pt

\begin{abstract}
In this report, a novel variation of Particle Swarm Optimization (PSO) algorithm, called Multiagent Coordination Optimization (MCO), is implemented in a parallel computing way for practical use by introducing MATLAB built-in function $\mathtt{parfor}$ into MCO. Then we rigorously analyze the global convergence of MCO by means of semistability theory. Besides sharing global optimal solutions with the PSO algorithm, the MCO algorithm integrates cooperative swarm behavior of multiple agents into the update formula by sharing velocity and position information between neighbors to improve its performance. Numerical evaluation of the parallel MCO algorithm is provided in the report by running the proposed algorithm on supercomputers in the High Performance Computing Center at Texas Tech University. In particular, the optimal value and consuming time are compared with PSO and serial MCO by solving several benchmark functions in the literature, respectively. Based on the simulation results, the performance of the parallel MCO is not only superb compared with PSO for solving many nonlinear, noncovex optimization problems, but also is of high efficiency by saving the computational time. 
\end{abstract}

\section{Introduction}
Particle Swarm Optimization (PSO) is a well developed
swarm intelligence method that optimizes a nonlinear or linear objective function iteratively by trying to improve a candidate solution with regards to
a given measure of quality. Motivated by a simplified social model,
the algorithm is first introduced by Kennedy and Eberhart in \cite{488968}, where some very primitive analysis of the convergence of PSO
is also provided. Since the PSO algorithm requires
only elementary mathematical operations and is computationally
efficient in terms of both memory requirements and speed, it solves many optimization problems quite efficiently, particularly some nonlinear, nonconvex optimization problems. Consequently,
the application of PSO has been widely seen from
interdisciplinary subjects ranging from computer science, engineering,
biology, to mathematics, economy \cite{1388490,1282114}, etc.
Several applications are reviewed in \cite{934374}, which includes
evolving neural networks, and reactive power and voltage control.

The mechanism of the PSO algorithm can be briefly explained as follows. The algorithm searches the
solution space of an objective function by updating the individual
solution vectors called particles. In the beginning, each particle
is assigned to a position in the solution space and a velocity
randomly. Each particle has a memory of its previous best value and
the corresponding previous best position. In addition, every
particle in the swarm can know the global best value among all
particles and the corresponding global best position. During every
iteration, the velocity of each particle is updated so that the
particle is guided by the previous best position of the particle and
the global best position stochastically.

As the PSO algorithm is used more extensively, more research efforts are devoted to its refinement. To improve the efficiency of the PSO
algorithm, the selection of the parameters becomes crucial.
References \cite{IoanCristian2003317,springerlink:10} study the relationship
between convergence rate and parameter selection while \cite{wight} focuses
on the impact of inertia weight and maximum velocity of PSO in which an adaptive inertia weight is equipped to guarantee the
convergence rate near a minimum. On the other hand, some variations
of PSO are proposed to improve the various aspects of the
algorithm, not limited to efficiency. In particular, \cite{785513} presents a simple variation with the addition of a
standard selection mechanism from evolutionary computation to improve performance. The authors in \cite{MOPSO:CC:2002,Hu02multiobjectiveoptimization,1304847} expand PSO to multiobjective optimization by augmenting the objective functions. More recently, a
new simple-structure variation of PSO is proposed \cite{ZSQ:DSCC:2011} to improve convergence. Unlike the
standard PSO, in this algorithm the particles can not
only communicate with each other via the objective function but also via
a new variable named ``quantizer" which displays a better
convergence than the standard PSO by evaluating some standard test functions in the literature.

All the above PSO variants focus either on some highly mathematical skills or on nature-inspired structures to improve their performance, lacking the fundamental understanding of how these algorithms work for \textit{general} problems. Thus, to address this issue, we
need to look at the swarm intelligence algorithm design from a new perspective since the traditional way of looking to \textit{natural} network systems appearing in nature for
inspiration does not provide a satisfactory answer. In particular, the new algorithms need to have \textit{robustness} properties on the practical uncertainty of distributed network implementation with communication constraints. Furthermore, due to the real-time implementation requirement for many network optimization systems in harsh or even adversarial environments, these new algorithms need to have faster (or even finite-time) convergence properties compared with the existing algorithms. Last but not least, these new algorithms need to have a capability of dealing with dynamical systems and control problems instead of just static optimization problems. In particular, it is favorable to use these new algorithms to \textit{modify} (control) the dynamic behavior of engineered network systems due to the inherent similarity between swarm optimization in computational intelligence \cite{BDT:1999} and cooperative networks in control theory \cite{LF:CDC:2001,Passino:CSM:2002,GP:TAC:2003,Gazi:TR:2005,LPP:TAC:2003,XTB:AUT:2006,HH:NARWA:2009}.

Multiagent Coordination Optimization
(MCO) algorithms are inspired by swarm intelligence and consensus protocols for multiagent coordination in \cite{HHB:TAC:2008,HHB:TAC:2009,Hui:IJC:2010,HH:AUT:2008}. Unlike the standard PSO, this new algorithm is a new optimization technique based not only on swarm
intelligence \cite{BDT:1999} which simulates the bio-inspired
behavior, but also on cooperative control of autonomous agents. Similar to PSO, the MCO algorithm
starts with a set of random solutions for agents which can
communicate with each other. The agents then move through the
solution space based on the evaluation of their cost functional and
neighbor-to-neighbor rules like multiagent consensus protocols
\cite{HHB:TAC:2008,HHB:TAC:2009,HH:AUT:2008,Hui:IJC:2010,Hui:TAC:2011,Hui:AUT:2011}.
By adding a distributed control term and gradient-based adaptation, we hope that the convergence speed of MCO can be accelerated and the convergence time of MCO can be improved compared with the existing techniques. Moreover, this new algorithm will be more suitable to distributed and parallel computation for solving large-scale physical network optimization problems by means of high performance computing facilities. 

In this report, we first implement MCO in a parallel computing way by introducing MATLAB\textsuperscript{\textregistered} built-in function $\mathtt{parfor}$ into MCO. Then we rigorously analyze the global convergence of MCO by means of \textit{semistability theory} \cite{HHB:TAC:2008,Hui:TAC:2013}. Besides sharing global optimal solutions with the PSO algorithm, the MCO algorithm incorporates cooperative swarm behavior of multiple agents into the update formula by sharing velocity and position information between neighbors to improve its performance. Numerical evaluation of the parallel MCO algorithm is provided by running the proposed algorithm on supercomputers in the High Performance Computing Center at Texas Tech University. In particular, the optimal solution and consuming time are compared with PSO and serial MCO by solving several benchmark functions in the literature, respectively. Based on the simulation results, the performance of the parallel MCO is not only superb compared with PSO by solving many nonlinear, nonconvex optimization problems, but also is of high efficiency by saving the computational time. 

This report is organized as follows. In Section \ref{mp}, some notions and notation in graph theory are introduced. In Section \ref{pmco} the realization of the parallel MCO algorithm in the MATLAB environment is described in details. The convergence results are developed in Section \ref{scr}. The numerical evaluation of the parallel MCO algorithm is then presented in Section \ref{ne}. Finally, Section \ref{co} concludes the report.

\section{Mathematical Preliminaries}\label{mp}
Graph theory is a powerful tool to investigate the topological change of large-scale network systems.
In this report, we use graph-related notation to describe our network topology based MCO algorithm. More specifically, let $\mathcal{G}(t)= (\mathcal{V},
\mathcal{E}(t), \mathcal{A}(t))$ denote a \textit{node-fixed dynamic directed graph} (or \textit{node-fixed dynamic digraph}) with the set
of vertices $\mathcal{V}= \{v_1,v_2,\ldots,v_{n}\}$ and
$\mathcal{E}(t)\subseteq \mathcal{V}\times \mathcal{V}$
represent the set of edges, where $t\in\overline{\mathbb{Z}}_{+}=\{0,1,2,\ldots\}$. The time-varying matrix $\mathcal{A}(t)$ with
nonnegative adjacency elements $a_{i,j}(t)$ serves as the weighted
adjacency matrix. The node index of $\mathcal{G}(t)$ is denoted as a
finite index set $\mathcal{N}=\{1,2,\ldots,n\}$. An edge of $\mathcal{G}(t)$ is
denoted by $e_{i,j}(t)=(v_i,v_j)$ and the adjacency elements associated
with the edges are positive. We assume $e_{i,j}(t)\in
\mathcal{E}(t)\Leftrightarrow a_{i,j}(t)=1$ and $a_{i,i}(t)=0$ for all $i\in
\mathcal{N}$. The set of neighbors of the node $v_i$ is denoted by
$\mathcal{N}^{i}(t)=\{v_j \in \mathcal {V}:(v_i,v_j)\in \mathcal {E}(t),
j=1,2, \ldots, |\mathcal{N}|, j\not = i\}$, where $|\mathcal{N}|$ denotes the cardinality of $\mathcal{N}$. The degree matrix of a node-fixed dynamic digraph
$\mathcal{G}(t)$ is defined as
 \begin{equation}
 \Delta(t) =[\delta_{i,j}(t)]_{i,j=1,2,\ldots,|\mathcal{N}|}
\end{equation} where\
\begin{eqnarray*}
\delta_{i,j}(t)=\left\{
      \begin{array}{ll}
        \sum_{j=1}^{|\mathcal{N}|} a_{i,j}(t), & \hbox{if $i=j$,} \\
        0, & \hbox{if $i\neq j$.}
      \end{array}
    \right.
\end{eqnarray*}
The \textit{Laplacian matrix} of the node-fixed dynamic digraph $\mathcal{G}(t)$ is defined by
\begin{equation}
L(t)=\Delta(t) - \mathcal {A}(t).
\end{equation} If $L(t)=L^{\mathrm{T}}(t)$, then $\mathcal{G}(t)$ is called a \textit{node-fixed dynamic undirected graph} (or simply \textit{node-fixed dynamic graph}). If there is a path from any node to any other node in a node-fixed dynamic digraph,
then we call the dynamic digraph \textit{strongly connected}. Analogously, if there is a path from any node to any other node in a node-fixed dynamic graph,
then we call the dynamic graph \textit{connected}. From now on we use short notations $L_{t},\mathcal{G}_{t},\mathcal{N}^{i}_{t}$ to denote $L(t),\mathcal{G}(t),\mathcal{N}^{i}(t)$, respectively. The following result due to Proposition 1 of \cite{AC:LAA:2005} is a property about the eigenvalue distribution of a Laplacian matrix.

\begin{lemma}[\hspace{-0.01em}\cite{AC:LAA:2005}]\label{lemma_argument}
Consider the Laplacian matrix $L_{t}$ for a node-fixed dynamic digraph or graph $\mathcal{G}_{t}$ with the index set $\mathcal{N}$, $|\mathcal{N}|\geq 2$. Let $\lambda_{t}\in{\mathrm{spec}}(L_{t})$, where ${\mathrm{spec}}(A)$ denotes the spectrum of $A$. Then for every $t\in\overline{\mathbb{Z}}_{+}$,
\begin{eqnarray}
-\Big(\frac{\pi}{2}-\frac{\pi}{|\mathcal{N}|}\Big)\leq\arg\lambda_{t}\leq\Big(\frac{\pi}{2}-\frac{\pi}{|\mathcal{N}|}\Big),
\end{eqnarray} where $\arg\lambda$ denotes the argument of $\lambda\in\mathbb{C}$, where $\mathbb{C}$ denotes the set of complex numbers.
\end{lemma}

A direct consequence from Lemma~\ref{lemma_argument} is that ${\mathrm{Re}}\,\lambda_{t}\geq0$, where ${\mathrm{Re}}\,\lambda$ denotes the real part of $\lambda\in\mathbb{C}$. Moreover, if $\mathcal{G}_{t}$ is an undirected graph, then $\lambda_{t}$ is real and $\lambda_{t}\geq0$.

\section{Parallel Multiagent Coordination Optimization Algorithm}\label{pmco}

\subsection{Multiagent Coordination Optimization with Node-Fixed Dynamic Graph Topology}

The MCO algorithm with static graph topology, proposed in \cite{ZH:CEC:2013} to solve a given optimization problem $\min_{\textbf{x}\in\mathbb{R}^{n}}f(\textbf{x})$, can be described in a vector form as follows:
\begin{eqnarray}
\textbf{v}_{k}(t+1)&=&\textbf{v}_{k}(t)+\eta\sum_{j\in\mathcal{N}^{k}}(\textbf{v}_{j}(t)-\textbf{v}_{k}(t))+\mu\sum_{j\in\mathcal{N}^{k}}(\textbf{x}_{j}(t)-\textbf{x}_{k}(t))+\kappa(\textbf{p}(t)-\textbf{x}_{i}(t)),\label{PSO_1}\\
\textbf{x}_{k}(t+1)&=&\textbf{x}_{k}(t)+\textbf{v}_{k}(t+1),\label{PSO_2}\\
\textbf{p}(t+1)&=&\left\{\begin{array}{ll}
\textbf{p}(t)+\kappa(\textbf{x}_{\min}(t)-\textbf{p}(t)), & {\mathrm{if}}\,\,\textbf{p}(t)\not\in\mathcal{Z},\\
\textbf{x}_{\min}(t), & {\mathrm{if}}\,\,\textbf{p}(t)\in\mathcal{Z},\\
\end{array}\right.\label{PSO_3}
\end{eqnarray}
where $k=1,\ldots,q$, $t\in\overline{\mathbb{Z}}_{+}$, $\textbf{v}_{k}(t)\in\mathbb{R}^{n}$ and $\textbf{x}_{k}(t)\in\mathbb{R}^{n}$ are the velocity
and position of particle $k$ at iteration $t$, respectively,
$\textbf{p}(t)\in\mathbb{R}^{n}$ is the position
of the global best value that the swarm of the particles can achieve
so far,
$\eta$, $\mu$, and $\kappa$ are three scalar random coefficients which are usually
selected in uniform distribution in the range $[0,1]$, $\mathcal{Z}=\{\textbf{y}\in\mathbb{R}^{n}:f(\textbf{x}_{\min})<f(\textbf{y})\}$, and $\textbf{x}_{\min}=\arg\min_{1\leq k\leq q}f(\textbf{x}_{k})$. In this report, we allow node-fixed dynamic graph topology in MCO so that $\mathcal{N}^{k}$ in (\ref{PSO_1}) becomes $\mathcal{N}^{k}(t)=\mathcal{N}_{t}^{k}$. A natural question arising from (\ref{PSO_1})--(\ref{PSO_3}) is the following: Can we always guarantee the convergence of (\ref{PSO_1})--(\ref{PSO_3}) for a given optimization problem $\min_{\textbf{x}\in\mathbb{R}^{n}}f(\textbf{x})$? Here convergence means that all the limits $\lim_{t\to\infty}\textbf{x}_{k}(t)$, $\lim_{t\to\infty}\textbf{v}_{k}(t)$, and $\lim_{t\to\infty}\textbf{p}(t)$ exist for every $k=1,\ldots,q$. This report tries to answer this question by giving some sufficient conditions to guarantee the convergence of (\ref{PSO_1})--(\ref{PSO_3}).  

\subsection{Parallel Implementation of MCO}

In this section, a parallel implementation of the MCO algorithm is introduced, which is described as Algorithm \ref{MCO} in the MATLAB language format. The command $\mathtt{matlabpool}$ opens or closes a pool of MATLAB sessions for parallel computation, and enables the parallel language features within the MATLAB  language (e.g., $\mathtt{parfor}$) by starting a parallel job which connects this MATLAB client with a number of labs.

The command $\mathtt{parfor}$ executes code loop in parallel. Part of the $\mathtt{parfor}$ body is executed on the MATLAB client (where the $\mathtt{parfor}$ is issued) and part is executed in parallel on MATLAB workers. The necessary data on which $\mathtt{parfor}$ operates is sent from the client to workers, where most of the computation happens, and the results are sent back to the client and pieced together. In Algorithm \ref{MCO}, the command $\mathtt{parfor}$ is used for loop of the update formula of all particles. Since the update formula needs the neighbors' information, so two temporary variables $C$ and $D$ are introduced for storing the global information of  position and velocity, respectively, and $L$ is the Laplacian matrix for the communication topology $\mathcal{G}$ for MCO.

\begin{algorithm}[ht]
\caption{Parallel MCO Algorithm}\label{MCO}
\begin{center}
\small\tt\begin{algorithmic}
\FOR{each agent $i=1,\ldots,q$}
  \STATE Initialize the agent's position with a uniformly distributed random vector: $x_{i}\sim U(\underline{x},\overline{x})\in
\mathbf{R}^{n\times 1}$, where $\underline{x}$ and $\overline{x}$ are the lower and upper boundaries of the search space;
  \STATE Initialize the agent's velocity: $v_{i}\sim U(\underline{v},\overline{v})$, where $\underline{v}$ and $\overline{v} \in
\mathbf{R}^{n\times 1}$ are the lower and upper boundaries of the search speed;
  \STATE Update the agent's best known position to its initial position: $p_{i}\leftarrow x_{i}$;
  \STATE If $f(p_{i})<f(p)$ update the multiagent network's best known position: $p\leftarrow p_i$.
\ENDFOR
  \REPEAT
   \STATE $k \leftarrow k+1$;\\
       \FOR {each agent $i=1,\ldots,q$}
       \STATE $C=[x_1,x_2,\cdots, x_q]^{\rm{T}}$, $D=[v_1,v_2,\cdots, v_q]^{\rm{T}}$;\\
       \STATE $\mathtt{parfor}$ {each agent $i=1,\ldots,q$}
      \STATE Choose random parameters: $\eta\sim U(0,1)$, $\mu\sim U(0,1)$, $\kappa\sim U(0,1)$;
      \STATE Update the agent's velocity: $v_{i}\leftarrow v_{i}+\eta(L_{k}(i,:)D)^{\rm{T}}+\mu(L_{k}(i,:)C)^{\rm{T}}+\kappa(p-x_{i})$;
      \STATE  Update the agent's position: $x_{i}\leftarrow x_{i}+v_{i}$;\\
          $\mathtt{endparfor}$\\
    \FOR{$f(x_{i})<f(p_{i})$}
    \STATE Update the agent's best known position: $p_{i}\leftarrow x_{i}$;
    \STATE Update the multiagent network's best known position: $p\leftarrow p+\kappa(p_{i}-p)$;
    \STATE If $f(p_{i})<f(p)$ update the multiagent network's best known position: $p\leftarrow p_{i}$;\\
     \ENDFOR
     \ENDFOR
   \UNTIL{$k$ is large enough or the value of $f$ has small change}
  \RETURN{$p$}
\end{algorithmic}
\end{center}
\end{algorithm}

\section{Convergence analysis}\label{scr}
In this section, we present some theoretic results on global convergence of the iterative process in Algorithm~\ref{MCO}. In particular, we view the randomized MCO algorithm as a discrete-time switched linear system and then use semistability theory to rigorously show its global convergence. 
To proceed with presentation, let $\mathbb{R}$ denote the set of real numbers.

\begin{lemma}\label{lemma_EW}
Let $n,q$ be positive integers and $q\geq 2$. For every $j=1,\ldots,q$, let $E_{n\times nq}^{[j]}\in\mathbb{R}^{n\times nq}$ denote a block-matrix whose $j$th block-column is $I_{n}$ and the rest block-elements are all zero matrices, i.e., $E_{n\times nq}^{[j]}=[\textbf{0}_{n\times n},\ldots,\textbf{0}_{n\times n},I_{n},\textbf{0}_{n\times n},\ldots,\textbf{0}_{n\times n}]$, $j=1,\ldots,q$, where $I_{m}\in\mathbb{R}^{m\times m}$ denotes the $m\times m$ identity matrix and $\textbf{0}_{m\times n}$ denotes the $m\times n$ zero matrix. Define $W^{[j]}=(\textbf{1}_{q\times 1}\otimes I_{n})E_{n\times nq}^{[j]}$ for every $j=1,\ldots,q$, where $\otimes$ denotes the Kronecker product and $\textbf{1}_{m\times n}$ denotes the $m\times n$ matrix whose entries are all ones. Then the following statements hold:
\begin{itemize}
\item[$i$)] For every $j=1,\ldots,q$, $W^{[j]}$ is an idempotent matrix, i.e., $(W^{[j]})^{2}=W^{[j]}$, and ${\mathrm{rank}}(W^{[j]}-I_{nq})=nq-n$, where ${\mathrm{rank}}(A)$ denotes the rank of $A$. 
\item[$ii$)] For any $\textbf{w}=[w_{1},\ldots,w_{q}]^{\mathrm{T}}\in\mathbb{R}^{q}$, $W^{[j]}(\textbf{w}\otimes\textbf{e}_{i})=w_{j}\textbf{1}_{q\times 1}\otimes\textbf{e}_{i}$ for every $j=1,\ldots,q$ and every $i=1,\ldots,n$. In particular, $W^{[j]}(\textbf{1}_{q\times 1}\otimes\textbf{e}_{i})=\textbf{1}_{q\times 1}\otimes\textbf{e}_{i}$ and $\ker(W^{[j]}-I_{nq})={\mathrm{span}}\{\textbf{1}_{q\times 1}\otimes\textbf{e}_{1},\ldots,\textbf{1}_{q\times 1}\otimes\textbf{e}_{n}\}$ for every $j=1,\ldots,q$ and every $i=1,\ldots,n$, where $[\textbf{e}_{1},\ldots,\textbf{e}_{n}]=I_{n}$, $\ker(A)$ denotes the kernel of $A$, and ${\mathrm{span}}\,S$ denotes the span of $S$.  
\item[$iii$)]
$E_{n\times nq}^{[j]}(\textbf{1}_{q\times 1}\otimes\textbf{e}_{i})=\textbf{e}_{i}$, $E_{n\times nq}^{[j]}(\textbf{1}_{q\times 1}\otimes I_{n})=I_{n}$, and $(\textbf{1}_{q\times 1}\otimes I_{n})\textbf{e}_{i}=\textbf{1}_{q\times 1}\otimes\textbf{e}_{i}$ for every $j=1,\ldots,q$ and every $i=1,\ldots,n$.
\end{itemize}
\end{lemma}

\begin{IEEEproof}
$i$) First note that by Fact 7.4.3 of \cite[p.~445]{Bernstein:2009}, $W^{[j]}=(\textbf{1}_{q\times 1}\otimes I_{n})E_{n\times nq}^{[j]}=\textbf{1}_{q\times 1}\otimes E_{n\times nq}^{[j]}$ for every $j=1,\ldots,q$. Now it follows from Fact 7.4.20 of \cite[p.~446]{Bernstein:2009} that 
\begin{eqnarray}\label{Wj2}
&&\hspace{-2.3em}(W^{[j]})^{2}=(\textbf{1}_{q\times 1}\otimes E_{n\times nq}^{[j]})(\textbf{1}_{q\times 1}\otimes E_{n\times nq}^{[j]})=(\textbf{1}_{q\times 1}\otimes[\textbf{0}_{n\times n},\ldots,\textbf{0}_{n\times n},I_{n},\textbf{0}_{n\times n},\ldots,\textbf{0}_{n\times n}])^{2}\nonumber\\
&&\hspace{0.9em}=[\textbf{1}_{q\times 1}\otimes\textbf{0}_{n\times n},\ldots,\textbf{1}_{q\times 1}\otimes\textbf{0}_{n\times n},\textbf{1}_{q\times 1}\otimes I_{n},\textbf{1}_{q\times 1}\otimes\textbf{0}_{n\times n},\ldots,\textbf{1}_{q\times 1}\otimes\textbf{0}_{n\times n}]^{2}\nonumber\\
&&\hspace{0.8em}=\small\left[\begin{array}{ccccccc}
\textbf{0}_{n\times n} & \ldots & \textbf{0}_{n\times n} & I_{n} & \textbf{0}_{n\times n} & \ldots & \textbf{0}_{n\times n}\\
\vdots & \ddots & \vdots & \vdots & \vdots & \ddots & \vdots\\
\textbf{0}_{n\times n} & \ldots & \textbf{0}_{n\times n} & I_{n} & \textbf{0}_{n\times n} & \ldots & \textbf{0}_{n\times n}\\
\end{array}\right]\small\left[\begin{array}{ccccccc}
\textbf{0}_{n\times n} & \ldots & \textbf{0}_{n\times n} & I_{n} & \textbf{0}_{n\times n} & \ldots & \textbf{0}_{n\times n}\\
\vdots & \ddots & \vdots & \vdots & \vdots & \ddots & \vdots\\
\textbf{0}_{n\times n} & \ldots & \textbf{0}_{n\times n} & I_{n} & \textbf{0}_{n\times n} & \ldots & \textbf{0}_{n\times n}\\
\end{array}\right]\nonumber\\
&&\hspace{0.8em}=\small\left[\begin{array}{ccccccc}
\textbf{0}_{n\times n} & \ldots & \textbf{0}_{n\times n} & I_{n} & \textbf{0}_{n\times n} & \ldots & \textbf{0}_{n\times n}\\
\vdots & \ddots & \vdots & \vdots & \vdots & \ddots & \vdots\\
\textbf{0}_{n\times n} & \ldots & \textbf{0}_{n\times n} & I_{n} & \textbf{0}_{n\times n} & \ldots & \textbf{0}_{n\times n}\\
\end{array}\right]=W^{[j]},
\end{eqnarray} which shows that $W^{[j]}$ is idempotent. 

Next, it follows from (\ref{Wj2}) that ${\mathrm{rank}}(W^{[j]})=n$ for every $j=1,\ldots,q$. By Sylvester's inequality, we have ${\mathrm{rank}}(W^{[j]}-I_{nq})+{\mathrm{rank}}(W^{[j]})\leq{\mathrm{rank}}((W^{[j]})^{2}-W^{[j]})+nq=nq$, and hence, ${\mathrm{rank}}(W^{[j]}-I_{nq})\leq nq-n$ for every $j=1,\ldots,q$. On the other hand, since $I_{nq}-W^{[j]}+W^{[j]}=I_{nq}$, it follows from $iv$) of Fact 2.10.17 of \cite[p.~127]{Bernstein:2009} that ${\mathrm{rank}}(I_{nq}-W^{[j]})+{\mathrm{rank}}(W^{[j]})\geq{\mathrm{rank}}(I_{nq}-W^{[j]}+W^{[j]})={\mathrm{rank}}(I_{nq})=nq$, which implies that ${\mathrm{rank}}(I_{nq}-W^{[j]})\geq nq-n$ for every $j=1,\ldots,q$. Thus, ${\mathrm{rank}}(W^{[j]}-I_{nq})=nq-n$ for every $j=1,\ldots,q$. 

$ii$) It follows from (\ref{Wj2}) that for every $j=1,\ldots,q$ and every $i=1,\ldots,n$,
\begin{eqnarray*}
W^{[j]}(\textbf{1}_{q\times 1}\otimes\textbf{e}_{i})=\small\left[\begin{array}{ccccccc}
\textbf{0}_{n\times n} & \ldots & \textbf{0}_{n\times n} & I_{n} & \textbf{0}_{n\times n} & \ldots & \textbf{0}_{n\times n}\\
\vdots & \ddots & \vdots & \vdots & \vdots & \ddots & \vdots\\
\textbf{0}_{n\times n} & \ldots & \textbf{0}_{n\times n} & I_{n} & \textbf{0}_{n\times n} & \ldots & \textbf{0}_{n\times n}\\
\end{array}\right]\small\left[\begin{array}{c}
\textbf{e}_{i}\\
\vdots \\
\textbf{e}_{i}\\
\end{array}\right]=\small\left[\begin{array}{c}
\textbf{e}_{i}\\
\vdots \\
\textbf{e}_{i}\\
\end{array}\right]=\textbf{1}_{q\times 1}\otimes\textbf{e}_{i},
\end{eqnarray*} namely, $(W^{[j]}-I_{nq})(\textbf{1}_{q\times 1}\otimes\textbf{e}_{i})=\textbf{0}_{nq\times 1}$ for every $j=1,\ldots,q$. Since by $i$), ${\mathrm{rank}}(W^{[j]}-I_{nq})=nq-n$ for every $j=1,\ldots,q$, it follows from Corollary 2.5.5 of \cite[p.~105]{Bernstein:2009} that ${\mathrm{def}}(W^{[j]}-I_{nq})=nq-{\mathrm{rank}}(W^{[j]}-I_{nq})=n$ for every $j=1,\ldots,q$, where ${\mathrm{def}}(A)=\dim\ker(A)$ denotes the defect of $A$ and $\dim S$ denotes the dimension of a subspace $S$. Note that $\textbf{1}_{q\times 1}\otimes\textbf{e}_{i}$, $i=1,\ldots,n$, are linearly independent, it follows that $\ker(W^{[j]}-I_{nq})={\mathrm{span}}\{\textbf{1}_{q\times 1}\otimes\textbf{e}_{1},\ldots,\textbf{1}_{q\times 1}\otimes\textbf{e}_{n}\}$ for every $j=1,\ldots,q$. 

Finally, for any $\textbf{w}=[w_{1},\ldots,w_{q}]^{\mathrm{T}}\in\mathbb{R}^{q}$, it follows from (\ref{Wj2}) that
\begin{eqnarray*}
W^{[j]}(\textbf{w}\otimes\textbf{e}_{i})=\small\left[\begin{array}{ccccccc}
\textbf{0}_{n\times n} & \ldots & \textbf{0}_{n\times n} & I_{n} & \textbf{0}_{n\times n} & \ldots & \textbf{0}_{n\times n}\\
\vdots & \ddots & \vdots & \vdots & \vdots & \ddots & \vdots\\
\textbf{0}_{n\times n} & \ldots & \textbf{0}_{n\times n} & I_{n} & \textbf{0}_{n\times n} & \ldots & \textbf{0}_{n\times n}\\
\end{array}\right]\small\left[\begin{array}{c}
w_{1}\textbf{e}_{i}\\
\vdots \\
w_{q}\textbf{e}_{i}\\
\end{array}\right]=\small\left[\begin{array}{c}
w_{j}\textbf{e}_{i}\\
\vdots \\
w_{j}\textbf{e}_{i}\\
\end{array}\right]=w_{j}\textbf{1}_{q\times 1}\otimes\textbf{e}_{i}
\end{eqnarray*} for every $j=1,\ldots,q$ and every $i=1,\ldots,n$.

$iii$) For every $j=1,\ldots,q$ and every $i=1,\ldots,n$,
$E_{n\times nq}^{[j]}(\textbf{1}_{q\times 1}\otimes\textbf{e}_{i})=[\textbf{0}_{n\times n},\ldots,\textbf{0}_{n\times n},I_{n},\textbf{0}_{n\times n},\ldots,\textbf{0}_{n\times n}]$\\$[\textbf{e}_{i}^{\mathrm{T}},\ldots,\textbf{e}_{i}^{\mathrm{T}}]^{\mathrm{T}}=\textbf{e}_{i}$ and 
$E_{n\times nq}^{[j]}(\textbf{1}_{q\times 1}\otimes I_{n})=[\textbf{0}_{n\times n},\ldots,\textbf{0}_{n\times n},I_{n},\textbf{0}_{n\times n},\ldots,\textbf{0}_{n\times n}][I_{n},\ldots,I_{n}]^{\mathrm{T}}=I_{n}$. Finally, by Fact 7.4.3 of \cite[p.~445]{Bernstein:2009}, $(\textbf{1}_{q\times 1}\otimes I_{n})\textbf{e}_{i}=\textbf{1}_{q\times 1}\otimes\textbf{e}_{i}$ for every $i=1,\ldots,n$.
\end{IEEEproof}

Next, we use some graph notions to state a result on the rank of certain matrices related to the matrix form of the iterative process in Algorithm~\ref{MCO}. 

\begin{lemma}\label{lemma_Arank}
Define a (possibly infinite) series of matrices $A^{[j]}_{k}$, $j=1,\ldots,q$, $k=0,1,2,\ldots$, as follows:
\begin{eqnarray}\label{Amatrix}
A_{k}^{[j]}=\small\left[\begin{array}{ccc}
\textbf{0}_{nq\times nq} & I_{nq} & \textbf{0}_{nq\times n} \\
-\mu_{k} L_{k}\otimes I_{n}-\kappa_{k} I_{nq} & -\eta_{k} L_{k}\otimes I_{n} & \kappa_{k} \textbf{1}_{q\times 1}\otimes I_{n} \\
\kappa_{k} E_{n\times nq}^{[j]} & \textbf{0}_{n\times nq} & -\kappa_{k} I_{n} \\
\end{array}\right],
\end{eqnarray} where $\mu_{k},\eta_{k},\kappa_{k}\geq0$, $k\in\overline{\mathbb{Z}}_{+}$, $L_{k}\in\mathbb{R}^{q\times q}$ denotes the Laplacian matrix of a node-fixed dynamic digraph $\mathcal{G}_{k}$, and $E_{n\times nq}^{[j]}\in\mathbb{R}^{n\times nq}$ is defined in Lemma~\ref{lemma_EW}. 
\begin{itemize}
\item[$i$)] If $\mu_{k}=0$ and $\kappa_{k}=0$, then ${\mathrm{rank}}(A_{k}^{[j]})=nq$ and $\ker(A_{k}^{[j]})=\{[\sum_{i=1}^{n}\sum_{l=1}^{q}\alpha_{il}(\textbf{e}_{i}\otimes\textbf{g}_{l})^{\mathrm{T}},\textbf{0}_{1\times nq},\sum_{i=1}^{n}\beta_{i}\\\textbf{e}_{i}^{\mathrm{T}}]^{\mathrm{T}}:\forall\alpha_{il}\in\mathbb{R},\forall\beta_{i}\in\mathbb{R},i=1,\ldots,n,l=1,\ldots,q\}$ for every $j=1,\ldots,q$, $k\in\overline{\mathbb{Z}}_{+}$, where $[\textbf{g}_{1},\ldots,\textbf{g}_{q}]=I_{q}$.
\item[$ii$)] If $\kappa_{k}\neq0$, then ${\mathrm{rank}}(A_{k}^{[j]})=2nq$ and $\ker(A_{k}^{[j]})=\{[\sum_{i=1}^{n}\alpha_{i}(\textbf{1}_{q\times 1}\otimes\textbf{e}_{i})^{\mathrm{T}},\textbf{0}_{1\times nq},\sum_{i=1}^{n}\alpha_{i}\textbf{e}_{i}^{\mathrm{T}}]^{\mathrm{T}}:\forall\alpha_{i}\in\mathbb{R},i=1,\ldots,n\}$ for every $j=1,\ldots,q$, $k\in\overline{\mathbb{Z}}_{+}$.
\item[$iii$)] If $\mu_{k}\neq0$ and $\kappa_{k}=0$, then ${\mathrm{rank}}(A_{k}^{[j]})=n(q+{\mathrm{rank}}(L_{k}))$ and $\ker(A_{k}^{[j]})=\{[\sum_{l=0}^{q-1-{\mathrm{rank}}(L_{k})}\sum_{i=1}^{n}\alpha_{li}(\textbf{w}_{l}\otimes\textbf{e}_{i})^{\mathrm{T}},\textbf{0}_{1\times nq},\sum_{i=1}^{n}\beta_{i}\textbf{e}_{i}^{\mathrm{T}}]^{\mathrm{T}}:\forall\alpha_{li}\in\mathbb{R},\forall\beta_{i}\in\mathbb{R},l=0,1,\ldots,q-1-{\mathrm{rank}}(L_{k}),i=1,\ldots,n\}$ for every $j=1,\ldots,q$, $k\in\overline{\mathbb{Z}}_{+}$.
\end{itemize}
\end{lemma}

\begin{IEEEproof}
First, it follows from (\ref{Amatrix}) that $\ker(A_{k}^{[j]})=\{[\textbf{z}_{1}^{\rm{T}},\textbf{z}_{2}^{\rm{T}},\textbf{z}_{3}^{\rm{T}}]^{\rm{T}}\in\mathbb{R}^{2nq+n}:\textbf{z}_{2}=\textbf{0}_{nq\times 1},-\mu_{k} (L_{k}\otimes I_{n})\textbf{z}_{1}-\kappa_{k}\textbf{z}_{1}-\eta_{k}(L_{k}\otimes I_{n})\textbf{z}_{2}+\kappa_{k}(\textbf{1}_{q\times 1}\otimes I_{n})\textbf{z}_{3}=\textbf{0}_{nq\times 1}, \kappa_{k}E_{n\times nq}^{[j]}\textbf{z}_{1}-\kappa_{k}\textbf{z}_{3}=\textbf{0}_{n\times 1}\}$, $k\in\overline{\mathbb{Z}}_{+}$, where $\textbf{z}_{1},\textbf{z}_{2}\in\mathbb{R}^{nq}$ and $\textbf{z}_{3}\in\mathbb{R}^{n}$. 

$i$) If $\mu_{k}=0$ and $\kappa_{k}=0$, then $\textbf{z}_{1}\in\mathbb{R}^{nq}$ and $\textbf{z}_{3}\in\mathbb{R}^{n}$ in $\ker(A_{k}^{[j]})$ can be chosen arbitrarily in $\mathbb{R}^{nq}$ and $\mathbb{R}^{n}$, respectively. Thus, $\textbf{z}_{1}$ and $\textbf{z}_{3}$ can be represented as  $\textbf{z}_{1}=\sum_{i=1}^{n}\sum_{l=1}^{q}\alpha_{il}(\textbf{e}_{i}\otimes\textbf{g}_{l})$ and $\textbf{z}_{3}=\sum_{i=1}^{n}\beta_{i}\textbf{e}_{i}$, where $\alpha_{li},\beta_{i}\in\mathbb{R}$. In this case, $\ker(A_{k}^{[j]})=\{[\textbf{z}_{1}^{\rm{T}},\textbf{z}_{2}^{\rm{T}},\textbf{z}_{3}^{\rm{T}}]^{\rm{T}}\in\mathbb{R}^{2nq+n}:\textbf{z}_{1}=\sum_{i=1}^{n}\sum_{l=1}^{q}\alpha_{il}(\textbf{e}_{i}\otimes\textbf{g}_{l}),\textbf{z}_{2}=\textbf{0}_{nq\times 1},\textbf{z}_{3}=\sum_{i=1}^{n}\beta_{i}\textbf{e}_{i},\forall\alpha_{li}\in\mathbb{R},\forall\beta_{i}\in\mathbb{R},i=1,\ldots,n,l=1,\ldots,q\}$ and ${\mathrm{def}}(A_{k}^{[j]})=nq+n$ for every $j=1,\ldots,q$, $k\in\overline{\mathbb{Z}}_{+}$. By Corollary 2.5.5 of \cite[p.~105]{Bernstein:2009}, ${\mathrm{rank}}(A_{k}^{[j]})=2nq+n-{\mathrm{def}}(A_{k}^{[j]})=nq$ for every $j=1,\ldots,q$, $k\in\overline{\mathbb{Z}}_{+}$.

$ii$) We consider two cases on $\mu_{k}$. 

\textit{Case 1.} If $\mu_{k}=0$ and $\kappa_{k}\neq0$, then 
substituting $\textbf{z}_{2}=\textbf{0}_{nq\times 1}$ and $\textbf{z}_{3}=E_{n\times nq}^{[j]}\textbf{z}_{1}$ into $-\kappa_{k}\textbf{z}_{1}-\eta_{k}(L_{k}\otimes I_{n})\textbf{z}_{2}+\kappa_{k}(\textbf{1}_{q\times 1}\otimes I_{n})\textbf{z}_{3}=\textbf{0}_{nq\times 1}$ yields 
\begin{eqnarray}\label{z1p}
\kappa_{k}(W^{[j]}-I_{nq})\textbf{z}_{1}=\textbf{0}_{nq\times 1},
\end{eqnarray}
where $W^{[j]}$ is defined in Lemma~\ref{lemma_EW}. Since, by $ii$) of Lemma~\ref{lemma_EW}, $\ker(W^{[j]}-I_{nq})={\mathrm{span}}\{\textbf{1}_{q\times 1}\otimes\textbf{e}_{1},\ldots,\textbf{1}_{q\times 1}\otimes\textbf{e}_{n}\}$ for every $j=1,\ldots,q$ and every $i=1,\ldots,n$, it follows from (\ref{z1p}) that $\textbf{z}_{1}$ can be represented as $\textbf{z}_{1}=\sum_{i=1}^{n}\alpha_{i}\textbf{1}_{q\times 1}\otimes\textbf{e}_{i}$, where $\alpha_{i}\in\mathbb{R}$. Furthermore, it follows from $iii$) of Lemma~\ref{lemma_EW} that $\textbf{z}_{3}=E_{n\times nq}^{[j]}\textbf{z}_{1}=\sum_{i=1}^{n}\alpha_{i}E_{n\times nq}^{[j]}(\textbf{1}_{q\times 1}\otimes\textbf{e}_{i})=\sum_{i=1}^{n}\alpha_{i}\textbf{e}_{i}$ for every $j=1,\ldots,q$. Thus, $\ker(A_{k}^{[j]})=\{[\sum_{i=1}^{n}\alpha_{i}(\textbf{1}_{q\times 1}\otimes\textbf{e}_{i})^{\mathrm{T}},\textbf{0}_{1\times nq},\\\sum_{i=1}^{n}\alpha_{i}\textbf{e}_{i}^{\mathrm{T}}]^{\mathrm{T}}:\forall\alpha_{i}\in\mathbb{R},i=1,\ldots,n\}$ for every $j=1,\ldots,q$, $k\in\overline{\mathbb{Z}}_{+}$, which implies that ${\mathrm{def}}(A_{k}^{[j]})=n$ for every $j=1,\ldots,q$, $k\in\overline{\mathbb{Z}}_{+}$. Therefore, in this case ${\mathrm{rank}}(A_{k}^{[j]})=2nq+n-{\mathrm{def}}(A_{k}^{[j]})=2nq$.

\textit{Case 2.} If $\mu_{k}\neq0$ and $\kappa_{k}\neq0$, then we claim that $\kappa_{k}/\mu_{k}\not\in{\mathrm{spec}}(-L_{k})$. To see this, it follows from Lemma~\ref{lemma_argument} that for any $\lambda_{k}\in{\mathrm{spec}}(-L_{k})$, ${\mathrm{Re}}\,\lambda_{k}\leq0$. Furthermore, note that $L_{k}\textbf{1}_{q\times 1}=\textbf{0}_{q\times 1}$. Thus, if $\kappa_{k}\neq0$, then $0<\kappa_{k}/\mu_{k}\not\in{\mathrm{spec}}(-L_{k})$.
Now, substituting $\textbf{z}_{2}=\textbf{0}_{nq\times 1}$ and $\textbf{z}_{3}=E_{n\times nq}^{[j]}\textbf{z}_{1}$ into $-\mu_{k} (L_{k}\otimes I_{n})\textbf{z}_{1}-\kappa_{k}\textbf{z}_{1}-\eta_{k}(L_{k}\otimes I_{n})\textbf{z}_{2}+\kappa_{k}(\textbf{1}_{q\times 1}\otimes I_{n})\textbf{z}_{3}=\textbf{0}_{nq\times 1}$ yields 
 \begin{eqnarray}\label{z1}
 (-\mu_{k} L_{k}\otimes I_{n}-\kappa_{k} I_{nq}+\kappa_{k} W^{[j]})\textbf{z}_{1}=\textbf{0}_{nq\times 1}, \quad k\in\overline{\mathbb{Z}}_{+}.
 \end{eqnarray}
Note that $(L_{k}\otimes I_{n})W^{[j]}=(L_{k}\otimes I_{n})(\textbf{1}_{q\times 1}\otimes E_{n\times nq}^{[j]})=L_{k}\textbf{1}_{q\times 1}\otimes E_{n\times nq}^{[j]}=\textbf{0}_{q\times 1}\otimes E_{n\times nq}^{[j]}=\textbf{0}_{nq\times nq}$, $k\in\overline{\mathbb{Z}}_{+}$. Pre-multiplying $-L_{k}\otimes I_{n}$ on both sides of (\ref{z1}) yields 
 $(\mu_{k}(L_{k}\otimes I_{n})^{2}+\kappa_{k} L_{k}\otimes I_{n})\textbf{z}_{1}=(\mu_{k} L_{k}\otimes I_{n}+\kappa_{k} I_{nq})(L_{k}\otimes I_{n})\textbf{z}_{1}=\textbf{0}_{nq\times 1}$, $k\in\overline{\mathbb{Z}}_{+}$. Since $\kappa_{k}/\mu_{k}\not\in{\mathrm{spec}}(-L_{k})$ for every $k\in\overline{\mathbb{Z}}_{+}$, it follows that $\det(\mu_{k} L_{k}\otimes I_{n}+\kappa_{k} I_{nq})\neq0$, $k\in\overline{\mathbb{Z}}_{+}$, where $\det$ denotes the determinant. Hence, $(L_{k}\otimes I_{n})\textbf{z}_{1}=\textbf{0}_{nq\times 1}$, $k\in\overline{\mathbb{Z}}_{+}$. 

Let $\textbf{w}_{0}=\textbf{1}_{q\times 1}$. Note that $L_{k}\textbf{w}_{0}=\textbf{0}_{q\times 1}$, it follows from Fact 7.4.22 of \cite[p.~446]{Bernstein:2009} that $(L_{k}\otimes I_{n})(\textbf{w}_{0}\otimes\textbf{e}_{i})=\textbf{0}_{q\times 1}$ for every $i=1,\ldots,n$, and hence, ${\mathrm{span}}\{\textbf{w}_{0}\otimes\textbf{e}_{1},\ldots,\textbf{w}_{0}\otimes\textbf{e}_{n}\}\subseteq\ker(L_{k}\otimes I_{n})$. Next, let ${\mathrm{span}}\{\textbf{w}_{1},\ldots,\textbf{w}_{q-1-{\mathrm{rank}}(L_{k})}\}=\ker(L_{k})\backslash{\mathrm{span}}\{\textbf{w}_{0}\}$, it follows that $\bigcup_{i=0}^{q-1-{\mathrm{rank}}(L_{k})}{\mathrm{span}}\{\textbf{w}_{i}\otimes\textbf{e}_{1},\ldots,\textbf{w}_{i}\otimes\textbf{e}_{n}\}=\ker(L_{k}\otimes I_{n})$, $k\in\overline{\mathbb{Z}}_{+}$. Hence, $\textbf{z}_{1}=\sum_{l=0}^{q-1-{\mathrm{rank}}(L_{k})}\sum_{i=1}^{n}\alpha_{li}\textbf{w}_{l}\otimes\textbf{e}_{i}$, where $\alpha_{li}\in\mathbb{R}$ and $\alpha_{li}=0$ for every $i=1,\ldots,n$ if $\textbf{w}_{l}=\textbf{0}_{q\times 1}$ for some $l\in\{1,\ldots,q-1-{\mathrm{rank}}(L_{k})\}$. Substituting this $\textbf{z}_{1}$ into the left-hand side of (\ref{z1}) yields $(-\mu_{k} L_{k}\otimes I_{n}-\kappa_{k} I_{nq}+\kappa_{k} W^{[j]})\textbf{z}_{1}=\kappa_{k}(W^{[j]}-I_{nq})\textbf{z}_{1}=\kappa_{k}(W^{[j]}-I_{nq})(\sum_{l=0}^{q-1-{\mathrm{rank}}(L_{k})}\sum_{i=1}^{n}\alpha_{li}\textbf{w}_{l}\otimes\textbf{e}_{i})=\kappa_{k}\sum_{l=0}^{q-1-{\mathrm{rank}}(L_{k})}\sum_{i=1}^{n}\alpha_{li}W^{[j]}\textbf{w}_{l}\otimes\textbf{e}_{i}-\kappa_{k}\sum_{l=0}^{q-1-{\mathrm{rank}}(L_{k})}\sum_{i=1}^{n}\alpha_{li}\textbf{w}_{l}\otimes\textbf{e}_{i}$. Note that it follows from $ii$) of Lemma~\ref{lemma_EW} that $W^{[j]}\textbf{w}_{0}\otimes\textbf{e}_{i}=\textbf{w}_{0}\otimes\textbf{e}_{i}$ for every $j=1,\ldots,q$ and every $i=1,\ldots,n$. Let $\textbf{w}_{l}=[w_{l1},\ldots,w_{lq}]^{\mathrm{T}}\in\mathbb{R}^{q}$ for every $l=1,\ldots,q-1-{\mathrm{rank}}(L_{k})$, then it follows from $ii$) of Lemma~\ref{lemma_EW} that
\begin{eqnarray*}
&&\kappa_{k}\sum_{l=0}^{q-1-{\mathrm{rank}}(L_{k})}\sum_{i=1}^{n}\alpha_{li}W^{[j]}\textbf{w}_{l}\otimes\textbf{e}_{i}-\kappa_{k}\sum_{l=0}^{q-1-{\mathrm{rank}}(L_{k})}\sum_{i=1}^{n}\alpha_{li}\textbf{w}_{l}\otimes\textbf{e}_{i}\nonumber\\
&&=\kappa_{k}\sum_{l=1}^{q-1-{\mathrm{rank}}(L_{k})}\sum_{i=1}^{n}\alpha_{li}(W^{[j]}\textbf{w}_{l}\otimes\textbf{e}_{i}-\textbf{w}_{l}\otimes\textbf{e}_{i})\nonumber\\
&&=\kappa_{k}\sum_{l=1}^{q-1-{\mathrm{rank}}(L_{k})}\sum_{i=1}^{n}\alpha_{li}(w_{lj}\textbf{w}_{0}\otimes\textbf{e}_{i}-\textbf{w}_{l}\otimes\textbf{e}_{i})\nonumber\\
&&=\kappa_{k}\sum_{l=1}^{q-1-{\mathrm{rank}}(L_{k})}\sum_{i=1}^{n}\alpha_{li}(w_{lj}\textbf{w}_{0}-\textbf{w}_{l})\otimes\textbf{e}_{i}\nonumber\\
&&=\kappa_{k}\sum_{l=1}^{q-1-{\mathrm{rank}}(L_{k})}\sum_{i=1}^{n}\alpha_{li}w_{lj}\textbf{w}_{0}\otimes\textbf{e}_{i}+\kappa_{k}\sum_{l=1}^{q-1-{\mathrm{rank}}(L_{k})}\sum_{i=1}^{n}(-\alpha_{li})\textbf{w}_{l}\otimes\textbf{e}_{i}.
\end{eqnarray*} Note that $\textbf{w}_{l}\otimes\textbf{e}_{i}$, $l=0,1,\ldots,q-1-{\mathrm{rank}}(L_{k})$, $i=1,\ldots,n$, are linearly independent.  Hence, $\textbf{z}_{1}$ satisfies (\ref{z1}) if and only if $\alpha_{li}=0$ for every $i=1,\ldots,n$ and every $l=1,\ldots,q-1-{\mathrm{rank}}(L_{k})$. In this case, we have $\textbf{z}_{1}=\sum_{i=1}^{n}\alpha_{0i}\textbf{w}_{0}\otimes\textbf{e}_{i}$.

Note that by $iii$) of Lemma~\ref{lemma_EW}, $\textbf{z}_{3}=E_{n\times nq}^{[j]}\textbf{z}_{1}=\sum_{i=1}^{n}\alpha_{0i}E_{n\times nq}^{[j]}(\textbf{1}_{q\times 1}\otimes\textbf{e}_{i})=\sum_{i=1}^{n}\alpha_{0i}\textbf{e}_{i}$ for every $j=1,\ldots,q$. Thus, $\ker(A_{k}^{[j]})=\{[\sum_{i=1}^{n}\alpha_{i}(\textbf{1}_{q\times 1}\otimes\textbf{e}_{i})^{\mathrm{T}},\textbf{0}_{1\times nq},\sum_{i=1}^{n}\alpha_{i}\textbf{e}_{i}^{\mathrm{T}}]^{\mathrm{T}}:\forall\alpha_{i}\in\mathbb{R},i=1,\ldots,n\}$ for every $j=1,\ldots,q$, $k\in\overline{\mathbb{Z}}_{+}$. Clearly $\dim\ker(A_{k}^{[j]})=n$ for every $j=1,\ldots,q$, $k\in\overline{\mathbb{Z}}_{+}$. Therefore, it follows from Corollary 2.5.5 of \cite[p.~105]{Bernstein:2009} that  ${\mathrm{rank}}(A_{k}^{[j]})=2nq+n-{\mathrm{def}}(A_{k}^{[j]})=2nq$ for every $j=1,\ldots,q$, $k\in\overline{\mathbb{Z}}_{+}$.

$iii)$ If $\mu_{k}\neq0$ and $\kappa_{k}=0$, then $\textbf{z}_{2}=\textbf{0}_{nq\times 1}$, $-\mu_{k}(L_{k}\otimes I_{n})\textbf{z}_{1}=\textbf{0}_{nq\times 1}$, and $\textbf{z}_{3}$ in $\ker(A_{k}^{[j]})$ can be chosen arbitrarily in $\mathbb{R}^{n}$. Thus, $\textbf{z}_{3}$ can be represented as $\textbf{z}_{3}=\sum_{i=1}^{n}\beta_{i}\textbf{e}_{i}$, where $\beta_{i}\in\mathbb{R}$. In this case, since $(L_{k}\otimes I_{n})\textbf{z}_{1}=\textbf{0}_{nq\times 1}$, $k\in\overline{\mathbb{Z}}_{+}$, it follows from the similar arguments as in Case 2 of $ii$) that $\textbf{z}_{1}=\sum_{l=0}^{q-1-{\mathrm{rank}}(L_{k})}\sum_{i=1}^{n}\alpha_{li}\textbf{w}_{l}\otimes\textbf{e}_{i}$. Therefore, $\ker(A_{k}^{[j]})=\{[\sum_{l=0}^{q-1-{\mathrm{rank}}(L_{k})}\sum_{i=1}^{n}\alpha_{li}(\textbf{w}_{l}\otimes\textbf{e}_{i})^{\mathrm{T}},\textbf{0}_{1\times nq},\sum_{i=1}^{n}\beta_{i}\textbf{e}_{i}^{\mathrm{T}}]^{\mathrm{T}}:\forall\alpha_{li}\in\mathbb{R},\forall\beta_{i}\in\mathbb{R},l=0,1,\ldots,q-1-{\mathrm{rank}}(L_{k}),i=1,\ldots,n\}$ for every $j=1,\ldots,q$, $k\in\overline{\mathbb{Z}}_{+}$. Clearly $\dim\ker(A_{k}^{[j]})=n(q-{\mathrm{rank}}(L_{k}))+n$ for every $j=1,\ldots,q$, $k\in\overline{\mathbb{Z}}_{+}$. Therefore, it follows from Corollary 2.5.5 of \cite[p.~105]{Bernstein:2009} that  ${\mathrm{rank}}(A_{k}^{[j]})=2nq+n-{\mathrm{def}}(A_{k}^{[j]})=n(q+{\mathrm{rank}}(L_{k}))$ for every $j=1,\ldots,q$, $k\in\overline{\mathbb{Z}}_{+}$.
\end{IEEEproof}

It follows from Lemma~\ref{lemma_Arank} that 0 is an eigenvalue of $A_{k}^{[j]}$ for every $j=1,\ldots,q$ and every $k\in\overline{\mathbb{Z}}_{+}$. Next, we further investigate some relationships of the null spaces between a row-addition transformed matrix of $A_{k}^{[j]}$ and $A_{k}^{[j]}$ itself in order to unveil an important property of this eigenvalue 0 later.

\begin{lemma}\label{lemma_Ah}
Consider the (possibly infinitely many) matrices $A_{k}^{[j]}+h_{k}A_{{\mathrm{c}}k}$, $j=1,\ldots,q$, $k=0,1,2,\ldots$, where $A_{k}^{[j]}$ is defined by (\ref{Amatrix}) in Lemma~\ref{lemma_Arank},
\begin{eqnarray}\label{Ac}
A_{{\mathrm{c}}k}=\small\left[\begin{array}{ccc}
-\mu_{k} L_{k}\otimes I_{n}-\kappa_{k} I_{nq} & -\eta_{k} L_{k}\otimes I_{n} & \kappa_{k} \textbf{1}_{q\times 1}\otimes I_{n} \\
\textbf{0}_{nq\times nq} & \textbf{0}_{nq\times nq} & \textbf{0}_{nq\times n} \\
\textbf{0}_{n\times nq} & \textbf{0}_{n\times nq} & \textbf{0}_{n\times n} \\
\end{array}\right],
\end{eqnarray} and $\mu_{k},\kappa_{k},\eta_{k},h_{k}\geq0$, $k\in\overline{\mathbb{Z}}_{+}$. Then $\ker(A_{k}^{[j]})=\ker(A_{k}^{[j]}+h_{k}A_{{\mathrm{c}}k})$ and $\ker(A_{k}^{[j]}(A_{k}^{[j]}+h_{k}A_{{\mathrm{c}}k}))=\ker((A_{k}^{[j]}+h_{k}A_{{\mathrm{c}}k})^{2})$ for every $j=1,\ldots,q$ and every $k\in\overline{\mathbb{Z}}_{+}$.
\end{lemma}

\begin{IEEEproof}
To show that $\ker(A_{k}^{[j]})=\ker(A_{k}^{[j]}+h_{k}A_{{\mathrm{c}}k})$, note that for every $j=1,\ldots,q$, $\ker(A_{k}^{[j]})=\{[\textbf{z}_{1}^{\rm{T}},\textbf{z}_{2}^{\rm{T}},\textbf{z}_{3}^{\rm{T}}]^{\rm{T}}\in\mathbb{R}^{2nq+n}:\textbf{z}_{2}=\textbf{0}_{nq\times 1},-\mu_{k} (L_{k}\otimes I_{n})\textbf{z}_{1}-\kappa_{k}\textbf{z}_{1}-\eta_{k}(L_{k}\otimes I_{n})\textbf{z}_{2}+\kappa_{k}(\textbf{1}_{q\times 1}\otimes I_{n})\textbf{z}_{3}=\textbf{0}_{nq\times 1}, \kappa_{k}E_{n\times nq}^{[j]}\textbf{z}_{1}-\kappa_{k}\textbf{z}_{3}=\textbf{0}_{n\times 1}\}$, $k\in\overline{\mathbb{Z}}_{+}$. Alternatively, for every $j=1,\ldots,q$ and every $k\in\overline{\mathbb{Z}}_{+}$, let  $\textbf{y}=[\textbf{y}_{1}^{\rm{T}},\textbf{y}_{2}^{\rm{T}},\textbf{y}_{3}^{\rm{T}}]^{\rm{T}}\in\ker(A_{k}^{[j]}+h_{k}A_{{\mathrm{c}}k})$, where $\textbf{y}_{1},\textbf{y}_{2}\in\mathbb{R}^{nq}$ and $\textbf{y}_{3}\in\mathbb{R}^{n}$, we have
\begin{eqnarray}
h_{k}(-\mu_{k} L_{k}\otimes I_{n}-\kappa_{k} I_{nq})\textbf{y}_{1}+h_{k}(-\eta_{k} L_{k}\otimes I_{n})\textbf{y}_{2}+\textbf{y}_{2}+h_{k}(\kappa_{k} \textbf{1}_{q\times 1}\otimes I_{n})\textbf{y}_{3}=\textbf{0}_{nq\times 1},\label{y_1p}\\
(-\mu_{k} L_{k}\otimes I_{n}-\kappa_{k} I_{nq})\textbf{y}_{1}+(-\eta_{k} L_{k}\otimes I_{n})\textbf{y}_{2}+(\kappa_{k} \textbf{1}_{q\times 1}\otimes I_{n})\textbf{y}_{3}=\textbf{0}_{nq\times 1},\label{y_2p}\\
\kappa_{k} E_{n\times nq}^{[j]}\textbf{y}_{1}-\kappa_{k}\textbf{y}_{3}=\textbf{0}_{n\times 1}.\label{y_3p}
\end{eqnarray} Substituting (\ref{y_2p}) into (\ref{y_1p}) yields $\textbf{y}_{2}=\textbf{0}_{nq\times 1}$. Together with (\ref{y_2p}) and (\ref{y_3p}), we have $\textbf{y}\in\ker(A_{k}^{[j]})$, which implies that $\ker(A_{k}^{[j]}+h_{k}A_{{\mathrm{c}}k})\subseteq\ker(A_{k}^{[j]})$ for every $j=1,\ldots,q$ and every $k\in\overline{\mathbb{Z}}_{+}$. On the other hand, if $\textbf{y}\in\ker(A_{k}^{[j]})$, then $\textbf{y}_{2}=\textbf{0}_{nq\times 1}$, $-\mu_{k} (L_{k}\otimes I_{n})\textbf{y}_{1}-\kappa_{k}\textbf{y}_{1}-\eta_{k}(L_{k}\otimes I_{n})\textbf{y}_{2}+\kappa_{k}(\textbf{1}_{q\times 1}\otimes I_{n})\textbf{y}_{3}=\textbf{0}_{nq\times 1}$, and $\kappa_{k}E_{n\times nq}^{[j]}\textbf{y}_{1}-\kappa_{k}\textbf{y}_{3}=\textbf{0}_{n\times 1}$. Clearly in this case, (\ref{y_1p})--(\ref{y_3p}) hold, i.e., $\textbf{y}\in\ker(A_{k}^{[j]}+h_{k}A_{{\mathrm{c}}k})$, which implies that $\ker(A_{k}^{[j]})\subseteq\ker(A_{k}^{[j]}+h_{k}A_{{\mathrm{c}}k})$ for every $j=1,\ldots,q$ and every $k\in\overline{\mathbb{Z}}_{+}$. Thus, $\ker(A_{k}^{[j]})=\ker(A_{k}^{[j]}+h_{k}A_{{\mathrm{c}}k})$ for every $j=1,\ldots,q$ and every $k\in\overline{\mathbb{Z}}_{+}$.

Finally, to show that $\ker(A_{k}^{[j]}(A_{k}^{[j]}+h_{k}A_{{\mathrm{c}}k}))=\ker((A_{k}^{[j]}+h_{k}A_{{\mathrm{c}}k})^{2})$, note that $\ker((A_{k}^{[j]}+h_{k}A_{{\mathrm{c}}k})^{2})=\ker((A_{k}^{[j]}+h_{k}A_{{\mathrm{c}}k})(A_{k}^{[j]}+h_{k}A_{{\mathrm{c}}k}))$ for every $j=1,\ldots,q$ and every $k\in\overline{\mathbb{Z}}_{+}$. Let $\textbf{y}\in\ker((A_{k}^{[j]}+h_{k}A_{{\mathrm{c}}k})(A_{k}^{[j]}+h_{k}A_{{\mathrm{c}}k}))$, then $(A_{k}^{[j]}+h_{k}A_{{\mathrm{c}}k})\textbf{y}\in\ker(A_{k}^{[j]}+h_{k}A_{{\mathrm{c}}k})=\ker(A_{k}^{[j]})$, and hence, $\textbf{y}\in\ker((A_{k}^{[j]}+h_{k}A_{{\mathrm{c}}k})^{2})$, which implies that $\ker(A_{k}^{[j]}(A_{k}^{[j]}+h_{k}A_{{\mathrm{c}}k}))\subseteq\ker(A_{k}^{[j]}(A_{k}^{[j]}+h_{k}A_{{\mathrm{c}}k}))$ for every $j=1,\ldots,q$ and every $k\in\overline{\mathbb{Z}}_{+}$. Alternatively, let $\textbf{z}\in\ker(A_{k}^{[j]}(A_{k}^{[j]}+h_{k}A_{{\mathrm{c}}k}))$, then $(A_{k}^{[j]}+h_{k}A_{{\mathrm{c}}k})\textbf{z}\in\ker(A_{k}^{[j]})=\ker(A_{k}^{[j]}+h_{k}A_{{\mathrm{c}}k})$, and hence, $\textbf{z}\in\ker((A_{k}^{[j]}+h_{k}A_{{\mathrm{c}}k})^{2})$, which implies that $\ker(A_{k}^{[j]}(A_{k}^{[j]}+h_{k}A_{{\mathrm{c}}k}))\subseteq\ker((A_{k}^{[j]}+h_{k}A_{{\mathrm{c}}k})^{2})$ for every $j=1,\ldots,q$ and every $k\in\overline{\mathbb{Z}}_{+}$. Thus, $\ker(A_{k}^{[j]}(A_{k}^{[j]}+h_{k}A_{{\mathrm{c}}k}))=\ker((A_{k}^{[j]}+h_{k}A_{{\mathrm{c}}k})^{2})$ for every $j=1,\ldots,q$ and every $k\in\overline{\mathbb{Z}}_{+}$.
\end{IEEEproof}

Next, we assert that 0 is semisimple for $A_{k}^{[j]}+h_{k}A_{{\mathrm{c}}k}$. Recall from Definition 5.5.4 of \cite[p.~322]{Bernstein:2009} that $0$ is semisimple if its geometric multiplicity and algebraic multiplicity are equal.

\begin{lemma}\label{lemma_semisimple}
Consider the (possibly infinitely many) matrices $A_{k}^{[j]}+h_{k}A_{{\mathrm{c}}k}$, $j=1,\ldots,q$, $k=0,1,2,\ldots$, defined in Lemma~\ref{lemma_Ah}, where $\mu_{k},\kappa_{k},\eta_{k},h_{k}\geq0$, $k\in\overline{\mathbb{Z}}_{+}$.
\begin{itemize}
\item[$i$)] If $\kappa_{k}=0$ and $\mu_{k}=0$, then ${\mathrm{rank}}(A_{k}^{[j]}+h_{k}A_{{\mathrm{c}}k})=nq$ and 0 is not a semisimple eigenvalue of $A_{k}^{[j]}+h_{k}A_{{\mathrm{c}}k}$ for every $j=1,\ldots,q$, $k\in\overline{\mathbb{Z}}_{+}$.
\item[$ii$)] If $\kappa_{k}=0$ and $\mu_{k}\neq0$, then ${\mathrm{rank}}(A_{k}^{[j]}+h_{k}A_{{\mathrm{c}}k})=n(q+{\mathrm{rank}}(L_{k}))$ and 0 is not a semisimple eigenvalue of $A_{k}^{[j]}+h_{k}A_{{\mathrm{c}}k}$ for every $j=1,\ldots,q$, $k\in\overline{\mathbb{Z}}_{+}$. 
\item[$iii$)] If $\kappa_{k}\neq0$, then ${\mathrm{rank}}(A_{k}^{[j]}+h_{k}A_{{\mathrm{c}}k})=2nq$ and 0 is a semisimple eigenvalue of $A_{k}^{[j]}+h_{k}A_{{\mathrm{c}}k}$ for every $j=1,\ldots,q$, $k\in\overline{\mathbb{Z}}_{+}$.
\end{itemize}
\end{lemma}

\begin{IEEEproof}
First, it follows from Lemma~\ref{lemma_Ah} that $\ker(A_{k}^{[j]}+h_{k}A_{{\mathrm{c}}k})=\ker(A_{k}^{[j]})$, and hence ${\mathrm{def}}(A_{k}^{[j]}+h_{k}A_{{\mathrm{c}}k})={\mathrm{def}}(A_{k}^{[j]})$ for every $j=1,\ldots,q$ and every $k\in\overline{\mathbb{Z}}_{+}$. Thus,  ${\mathrm{rank}}(A_{k}^{[j]}+h_{k}A_{{\mathrm{c}}k})=2nq+n-{\mathrm{def}}(A_{k}^{[j]}+h_{k}A_{{\mathrm{c}}k})=2nq+n-{\mathrm{def}}(A_{k}^{[j]})={\mathrm{rank}}(A_{k}^{[j]})$ for every $j=1,\ldots,q$ and every $k\in\overline{\mathbb{Z}}_{+}$.
Therefore, all the rank conclusions on $A_{k}^{[j]}+h_{k}A_{{\mathrm{c}}k}$ in $i$)--$iii$) directly follow from Lemma~\ref{lemma_Arank}.

Next, it follows from these rank conclusions on $A_{k}^{[j]}+h_{k}A_{{\mathrm{c}}k}$ that $A_{k}^{[j]}+h_{k}A_{{\mathrm{c}}k}$ has an eigenvalue 0 for every $j=1,\ldots,q$ and every $k\in\overline{\mathbb{Z}}_{+}$. Now we want to further investigate whether 0 is a semisimple eigenvalue of $A_{k}^{[j]}+h_{k}A_{{\mathrm{c}}k}$ or not for every $j=1,\ldots,q$, $k\in\overline{\mathbb{Z}}_{+}$. To this end, we need to study the relationship between $\ker(A_{k}^{[j]})$ and $\ker(A_{k}^{[j]}(A_{k}^{[j]}+h_{k}A_{{\mathrm{c}}k}))$ for every $j=1,\ldots,q$, $k\in\overline{\mathbb{Z}}_{+}$.
  
Noting that $(L_{k}\otimes I_{n})(\textbf{1}_{q\times 1}\otimes I_{n})=(L_{k}\textbf{1}_{q\times 1})\otimes I_{n}=\textbf{0}_{nq\times n}$ and by $iii$) of Lemma~\ref{lemma_EW}, $E_{n\times nq}^{[j]}(\textbf{1}_{q\times 1}\otimes I_{n})=I_{n}$, we have 
\begin{eqnarray*}
&&(A_{k}^{[j]})^{2}=\small\left[\begin{array}{ccc}
-\mu_{k} L_{k}\otimes I_{n}-\kappa_{k} I_{nq} & -\eta_{k} L_{k}\otimes I_{n} & \kappa_{k} \textbf{1}_{q\times 1}\otimes I_{n} \\
\eta_{k}\mu_{k}(L_{k}\otimes I_{n})^{2}+\eta_{k}\kappa_{k} L_{k}\otimes I_{n}+\kappa_{k}^{2}W^{[j]} & \eta_{k}^{2}(L_{k}\otimes I_{n})^{2}-\mu_{k} L_{k}\otimes I_{n}-\kappa_{k} I_{nq} & -\kappa_{k}^{2} \textbf{1}_{q\times 1}\otimes I_{n} \\
-\kappa_{k}^{2} E_{n\times nq}^{[j]} & \kappa_{k} E_{n\times nq}^{[j]} & \kappa_{k}^{2} I_{n} \\
\end{array}\right],\\
&&A_{k}^{[j]}A_{{\mathrm{c}}k}=\small\left[\begin{array}{ccc}
\textbf{0}_{nq\times nq} & \textbf{0}_{nq\times nq} & \textbf{0}_{nq\times n} \\
\mu_{k}^{2} (L_{k}\otimes I_{n})^{2}+2\mu_{k}\kappa_{k}(L_{k}\otimes I_{n})+\kappa_{k}^{2} I_{nq} & \mu_{k}\eta_{k}(L_{k}\otimes I_{n})^{2}+\kappa_{k}\eta_{k}L_{k}\otimes I_{n} & -\kappa_{k}^{2} \textbf{1}_{q\times 1}\otimes I_{n} \\
-\kappa_{k}\mu_{k}E_{n\times nq}^{[j]}(L_{k}\otimes I_{n})-\kappa_{k}^{2}E_{n\times nq} ^{[j]} & -\kappa_{k}\eta_{k}E_{n\times nq}^{[j]}(L_{k}\otimes I_{n}) & \kappa_{k}^{2}I_{n} \\
\end{array}\right].
\end{eqnarray*} Thus, for every $j=1,\ldots,q$ and every $k\in\overline{\mathbb{Z}}_{+}$, let $\textbf{y}=[\textbf{y}_{1}^{\rm{T}},\textbf{y}_{2}^{\rm{T}},\textbf{y}_{3}^{\rm{T}}]^{\rm{T}}\in\ker(A_{k}^{[j]}(A_{k}^{[j]}+h_{k}A_{{\mathrm{c}}k}))$, where $\textbf{y}_{1},\textbf{y}_{2}\in\mathbb{R}^{nq}$ and $\textbf{y}_{3}\in\mathbb{R}^{n}$, we have 
\begin{eqnarray}
(-\mu_{k} L_{k}\otimes I_{n}-\kappa_{k} I_{nq})\textbf{y}_{1}-(\eta_{k} L_{k}\otimes I_{n})\textbf{y}_{2}+(\kappa_{k} \textbf{1}_{q\times 1}\otimes I_{n})\textbf{y}_{3}=\textbf{0}_{nq\times 1},\label{y_1}\\
(\eta_{k}\mu_{k}(L_{k}\otimes I_{n})^{2}+\eta_{k}\kappa_{k} L_{k}\otimes I_{n}+\kappa_{k}^{2}W^{[j]})\textbf{y}_{1}+(\eta_{k}^{2}(L_{k}\otimes I_{n})^{2}-\mu_{k} L_{k}\otimes I_{n}-\kappa_{k} I_{nq})\textbf{y}_{2}\nonumber\\
+(-\kappa_{k}^{2} \textbf{1}_{q\times 1}\otimes I_{n})\textbf{y}_{3}\nonumber\\
+h_{k}(\mu_{k}^{2} (L_{k}\otimes I_{n})^{2}+2\mu_{k}\kappa_{k}(L_{k}\otimes I_{n})+\kappa_{k}^{2} I_{nq})\textbf{y}_{1}+h_{k}(\mu_{k}\eta_{k}(L_{k}\otimes I_{n})^{2}+\kappa_{k}\eta_{k}L_{k}\otimes I_{n})\textbf{y}_{2}\nonumber\\
+h_{k}(-\kappa_{k}^{2} \textbf{1}_{q\times 1}\otimes I_{n})\textbf{y}_{3}=\textbf{0}_{nq\times 1},\label{y_2}\\
-\kappa_{k}^{2} E_{n\times nq}^{[j]}\textbf{y}_{1}+\kappa_{k}E_{n\times nq}^{[j]}\textbf{y}_{2}+\kappa_{k}^{2}\textbf{y}_{3}\nonumber\\
+h_{k}(-\kappa_{k}\mu_{k}E_{n\times nq}^{[j]}(L_{k}\otimes I_{n})-\kappa_{k}^{2}E_{n\times nq} ^{[j]})\textbf{y}_{1}+h_{k}(-\kappa_{k}\eta_{k}E_{n\times nq}^{[j]}(L_{k}\otimes I_{n}))\textbf{y}_{2}+h_{k}\kappa_{k}^{2}\textbf{y}_{3}=\textbf{0}_{n\times 1}.\label{y_3}
\end{eqnarray} Now we consider two cases on $\kappa_{k}$.

\textit{Case 1.} $\kappa_{k}=0$. In this case, (\ref{y_3}) becomes trivial and (\ref{y_1}) and (\ref{y_2}) become
\begin{eqnarray}
(-\mu_{k} L_{k}\otimes I_{n})\textbf{y}_{1}-(\eta_{k} L_{k}\otimes I_{n})\textbf{y}_{2}=\textbf{0}_{nq\times 1},\label{y_1k}\\
\eta_{k}\mu_{k}(L_{k}\otimes I_{n})^{2}\textbf{y}_{1}+(\eta_{k}^{2}(L_{k}\otimes I_{n})^{2}-\mu_{k} L_{k}\otimes I_{n})\textbf{y}_{2}\nonumber\\
+h_{k}\mu_{k}^{2} (L_{k}\otimes I_{n})^{2}\textbf{y}_{1}+h_{k}\mu_{k}\eta_{k}(L_{k}\otimes I_{n})^{2}\textbf{y}_{2}=\textbf{0}_{nq\times 1}.\label{y_2k}
\end{eqnarray} 

If $\mu_{k}=0$, then it follows from (\ref{y_1k}) and (\ref{y_2k}) that $-(\eta_{k}L_{k}\otimes I_{n})\textbf{y}_{2}=\textbf{0}_{nq\times 1}$ and $\eta_{k}^{2}(L_{k}\otimes I_{n})^{2}\textbf{y}_{2}=\textbf{0}_{nq\times 1}$. Hence, either $\eta_{k}=0$ or $(L_{k}\otimes I_{n})\textbf{y}_{2}=\textbf{0}_{nq\times 1}$. If $\eta_{k}=0$, then $\textbf{y}_{1},\textbf{y}_{2}\in\mathbb{R}^{nq}$ and $\textbf{y}_{3}\in\mathbb{R}^{n}$ can be chosen arbitrarily. Thus, $\ker(A_{k}^{[j]}(A_{k}^{[j]}+h_{k}A_{{\mathrm{c}}k}))=\mathbb{R}^{2nq+n}$, and it follows from $i$) of Lemma~\ref{lemma_Arank} that $\ker(A_{k}^{[j]}(A_{k}^{[j]}+h_{k}A_{{\mathrm{c}}k}))\neq\ker(A_{k}^{[j]})$. By Lemma~\ref{lemma_Ah}, we have $\ker((A_{k}^{[j]}+h_{k}A_{{\mathrm{c}}k})^{2})\neq\ker(A_{k}^{[j]}+h_{k}A_{{\mathrm{c}}k})$. Now, by Proposition 5.5.8 of \cite[p.~323]{Bernstein:2009}, 0 is not semisimple. Alternatively, if $\eta_{k}\neq0$, then $(L_{k}\otimes I_{n})\textbf{y}_{2}=\textbf{0}_{nq\times 1}$ and $\textbf{y}_{1}\in\mathbb{R}^{nq}$ and $\textbf{y}_{3}\in\mathbb{R}^{n}$ can be chosen arbitrarily. Using the similar arguments as in the proof of Case 2 in $ii$) of Lemma~\ref{lemma_Arank}, it follows that $\textbf{y}_{2}=\sum_{l=0}^{q-1-{\mathrm{rank}}(L_{k})}\sum_{i=1}^{n}\alpha_{li}\textbf{w}_{l}\otimes\textbf{e}_{i}$, where $\alpha_{li}\in\mathbb{R}$.
Hence, $\ker(A_{k}^{[j]}(A_{k}^{[j]}+h_{k}A_{{\mathrm{c}}k}))=\{[\sum_{i=1}^{n}\sum_{r=1}^{q}\beta_{ir}(\textbf{e}_{i}\otimes\textbf{e}_{r})^{\mathrm{T}},\sum_{l=0}^{q-1-{\mathrm{rank}}(L_{k})}\sum_{i=1}^{n}\alpha_{li}(\textbf{w}_{l}\otimes\textbf{e}_{i})^{\mathrm{T}},\sum_{i=1}^{n}\gamma_{i}\textbf{e}_{i}^{\mathrm{T}}]^{\mathrm{T}}:\forall\alpha_{li}\in\mathbb{R},\forall\beta_{ir}\in\mathbb{R},\forall\gamma_{i}\in\mathbb{R},i=1,\ldots,n,r=1,\ldots,q,l=0,\ldots,q-1-{\mathrm{rank}}(L_{k})\}$ for every $j=1,\ldots,q$, $k\in\overline{\mathbb{Z}}_{+}$. Clearly it follows from $i$) of Lemma~\ref{lemma_Arank} that $\ker(A_{k}^{[j]}(A_{k}^{[j]}+h_{k}A_{{\mathrm{c}}k}))\neq\ker(A_{k}^{[j]})$. By Lemma~\ref{lemma_Ah}, we have $\ker((A_{k}^{[j]}+h_{k}A_{{\mathrm{c}}k})^{2})\neq\ker(A_{k}^{[j]}+h_{k}A_{{\mathrm{c}}k})$. Now, by Proposition 5.5.8 of \cite[p.~323]{Bernstein:2009}, 0 is not semisimple.

If $\mu_{k}\neq0$, then substituting (\ref{y_1k}) into (\ref{y_2k}) yields $-\mu_{k}(L_{k}\otimes I_{n})\textbf{y}_{2}=\textbf{0}_{nq\times 1}$. Substituting this equation into (\ref{y_1k}) yields $-\mu_{k}(L_{k}\otimes I_{n})\textbf{y}_{1}=\textbf{0}_{nq\times 1}$. Using the similar arguments as in the proof of Case 2 in $ii$) of Lemma~\ref{lemma_Arank}, it follows that $\textbf{y}_{1}=\sum_{l=0}^{q-1-{\mathrm{rank}}(L_{k})}\sum_{i=1}^{n}\alpha_{li}\textbf{w}_{l}\otimes\textbf{e}_{i}$ and $\textbf{y}_{2}=\sum_{l=0}^{q-1-{\mathrm{rank}}(L_{k})}\sum_{i=1}^{n}\beta_{li}\textbf{w}_{l}\otimes\textbf{e}_{i}$, where $\alpha_{li},\beta_{li}\in\mathbb{R}$. Note that $\textbf{y}_{3}\in\mathbb{R}^{n}$ can be chosen arbitrarily, and hence, $\ker(A_{k}^{[j]}(A_{k}^{[j]}+h_{k}A_{{\mathrm{c}}k}))=\{[\sum_{l=0}^{q-1-{\mathrm{rank}}(L_{k})}\sum_{i=1}^{n}\alpha_{li}(\textbf{w}_{l}\otimes\textbf{e}_{i})^{\mathrm{T}},\sum_{l=0}^{q-1-{\mathrm{rank}}(L_{k})}\sum_{i=1}^{n}\beta_{li}(\textbf{w}_{l}\otimes\textbf{e}_{i})^{\mathrm{T}},\sum_{i=1}^{n}\gamma_{i}\textbf{e}_{i}^{\mathrm{T}}]^{\mathrm{T}}:\forall\alpha_{li}\in\mathbb{R},\forall\beta_{li}\in\mathbb{R},\forall\gamma_{i}\in\mathbb{R},i=1,\ldots,n,l=0,\ldots,q-1-{\mathrm{rank}}(L_{k})\}$ for every $j=1,\ldots,q$, $k\in\overline{\mathbb{Z}}_{+}$. Clearly it follows from $iii$) of Lemma~\ref{lemma_Arank} that $\ker(A_{k}^{[j]}(A_{k}^{[j]}+h_{k}A_{{\mathrm{c}}k}))\neq\ker(A_{k}^{[j]})$. By Lemma~\ref{lemma_Ah}, we have $\ker((A_{k}^{[j]}+h_{k}A_{{\mathrm{c}}k})^{2})\neq\ker(A_{k}^{[j]}+h_{k}A_{{\mathrm{c}}k})$. Now, by Proposition 5.5.8 of \cite[p.~323]{Bernstein:2009}, 0 is not semisimple.

\textit{Case 2.} $\kappa_{k}\neq0$. In this case, substituting (\ref{y_1}) into (\ref{y_2}) and (\ref{y_3}) yields
\begin{eqnarray}
(\eta_{k}\mu_{k}(L_{k}\otimes I_{n})^{2}+\eta_{k}\kappa_{k} L_{k}\otimes I_{n}+\kappa_{k}^{2}W^{[j]})\textbf{y}_{1}+(\eta_{k}^{2}(L_{k}\otimes I_{n})^{2}-\mu_{k} L_{k}\otimes I_{n}-\kappa_{k} I_{nq})\textbf{y}_{2}\nonumber\\
+(-\kappa_{k}^{2} \textbf{1}_{q\times 1}\otimes I_{n})\textbf{y}_{3}\nonumber\\
+h_{k}(\mu_{k}^{2} (L_{k}\otimes I_{n})^{2}+\mu_{k}\kappa_{k}(L_{k}\otimes I_{n}))\textbf{y}_{1}+h_{k}\mu_{k}\eta_{k}(L_{k}\otimes I_{n})^{2}\textbf{y}_{2}=\textbf{0}_{nq\times 1},\label{y12p}\\
-\kappa_{k}^{2} E_{n\times nq}^{[j]}\textbf{y}_{1}+\kappa_{k}E_{n\times nq}^{[j]}\textbf{y}_{2}+\kappa_{k}^{2}\textbf{y}_{3}=\textbf{0}_{n\times 1}.\label{y13p}
\end{eqnarray}
Note that $(L_{k}\otimes I_{n})W^{[j]}=(L_{k}\otimes I_{n})(\textbf{1}_{q\times 1}\otimes E_{n\times nq}^{[j]})=L_{k}\textbf{1}_{q\times 1}\otimes E_{n\times nq}^{[j]}=\textbf{0}_{q\times 1}\otimes E_{n\times nq}^{[j]}=\textbf{0}_{nq\times nq}$. Pre-multiplying $-L_{k}\otimes I_{n}$ on both sides of (\ref{y_1}) yields $(\mu_{k}(L_{k}\otimes I_{n})^{2}+\kappa_{k} L_{k}\otimes I_{n})\textbf{y}_{1}+\eta_{k}(L_{k}\otimes I_{n})^{2}\textbf{y}_{2}=\textbf{0}_{nq\times 1}$. Substituting this equation into (\ref{y12p}) yields
\begin{eqnarray}
\kappa_{k}^{2}W^{[j]}\textbf{y}_{1}+(-\mu_{k} L_{k}\otimes I_{n}-\kappa_{k} I_{nq})\textbf{y}_{2}+(-\kappa_{k}^{2} \textbf{1}_{q\times 1}\otimes I_{n})\textbf{y}_{3}=\textbf{0}_{nq\times 1}.\label{y12pp}
\end{eqnarray} Finally, substituting (\ref{y13p}) into (\ref{y_1}) and (\ref{y12pp}) by eliminating $\textbf{y}_{3}$ yields
\begin{eqnarray}
(-\mu_{k} L_{k}\otimes I_{n}-\kappa_{k} I_{nq}+\kappa_{k} W^{[j]})\textbf{y}_{1}-(\eta_{k} L_{k}\otimes I_{n}+W^{[j]})\textbf{y}_{2}=\textbf{0}_{nq\times 1},\label{y12}\\
(-\mu_{k} L_{k}\otimes I_{n}-\kappa_{k}I_{nq}+\kappa_{k}W^{[j]})\textbf{y}_{2}=\textbf{0}_{nq\times 1}.\label{y21}
\end{eqnarray} Note that (\ref{y21}) is identical to (\ref{z1}). Then it follows from the similar arguments as in the proof of Case 2 of $ii$) of Lemma~\ref{lemma_Arank} that $\textbf{y}_{2}=\sum_{i=1}^{n}\beta_{i}\textbf{1}_{q\times 1}\otimes\textbf{e}_{i}$, where $\beta_{i}\in\mathbb{R}$. Clearly $\textbf{y}_{2}\in\ker(L_{k}\otimes I_{n})$. Next, substituting this $\textbf{y}_{2}$ into $(\mu_{k}(L_{k}\otimes I_{n})^{2}+\kappa_{k} L_{k}\otimes I_{n})\textbf{y}_{1}+\eta_{k}(L_{k}\otimes I_{n})^{2}\textbf{y}_{2}=\textbf{0}_{nq\times 1}$ yields $(\mu_{k}(L_{k}\otimes I_{n})^{2}+\kappa_{k} L_{k}\otimes I_{n})\textbf{y}_{1}=\textbf{0}_{nq\times 1}$. If $\mu_{k}=0$, then $\kappa_{k}(L_{k}\otimes I_{n})\textbf{y}_{1}=\textbf{0}_{nq\times 1}$. Otherwise, if $\mu_{k}\neq0$, then $\det(\mu_{k}(L_{k}\otimes I_{n})+\kappa_{k} I_{nq})\neq0$, which implies that $(L_{k}\otimes I_{n})\textbf{y}_{1}=\textbf{0}_{nq\times 1}$. Again, it follows from the similar arguments as in the proof of $ii$) of Lemma~\ref{lemma_Arank} that $\textbf{y}_{1}=\sum_{i=1}^{n}\gamma_{i}\textbf{1}_{q\times 1}\otimes\textbf{e}_{i}$, where $\gamma_{i}\in\mathbb{R}$. Clearly it follows from $ii$) of Lemma~\ref{lemma_EW} that $W^{[j]}\textbf{y}_{1}=\textbf{y}_{1}$ and $W^{[j]}\textbf{y}_{2}=\textbf{y}_{2}$ for every $j=1,\ldots,q$. Now substituting $\textbf{y}_{1}$ and $\textbf{y}_{2}$ into the left-hand side of (\ref{y12}) yields $(-\mu_{k} L_{k}\otimes I_{n}-\kappa_{k} I_{nq}+\kappa_{k} W^{[j]})\textbf{y}_{1}-(\eta_{k} L_{k}\otimes I_{n}+W^{[j]})\textbf{y}_{2}=-\textbf{y}_{2}=-\sum_{i=1}^{n}\beta_{i}(\textbf{1}_{q\times 1}\otimes\textbf{e}_{i})$. Thus, (\ref{y12}) holds if and only if $\sum_{i=1}^{n}\beta_{i}(\textbf{1}_{q\times 1}\otimes\textbf{e}_{i})=\textbf{0}_{nq\times 1}$, which implies that $\beta_{i}=0$ for every $i=1,\ldots,n$, that is, $\textbf{y}_{2}=\textbf{0}_{nq\times 1}$. Then it follows from (\ref{y13p}) and $iii$) of Lemma~\ref{lemma_EW} that $\textbf{y}_{3}=E_{n\times nq}^{[j]}\textbf{y}_{1}=\sum_{i=1}^{n}\gamma_{i}E_{n\times nq}^{[j]}(\textbf{1}_{q\times 1}\otimes\textbf{e}_{i})=\sum_{i=1}^{n}\gamma_{i}\textbf{e}_{i}$. Clearly such $\textbf{y}_{1}=\sum_{i=1}^{n}\gamma_{i}\textbf{1}_{q\times 1}\otimes\textbf{e}_{i}$, $\textbf{y}_{2}=\textbf{0}_{nq\times 1}$, and $\textbf{y}_{3}=\sum_{i=1}^{n}\gamma_{i}\textbf{e}_{i}$ satisfy (\ref{y_1})--(\ref{y_3}). Thus, $\ker(A_{k}^{[j]}(A_{k}^{[j]}+h_{k}A_{{\mathrm{c}}k}))=\{[\sum_{i=1}^{n}\gamma_{i}(\textbf{1}_{q\times 1}\otimes\textbf{e}_{i})^{\mathrm{T}},\textbf{0}_{1\times nq},\sum_{i=1}^{n}\gamma_{i}\textbf{e}_{i}^{\mathrm{T}}]^{\mathrm{T}}:\forall\gamma_{i}\in\mathbb{R},i=1,\ldots,n\}=\ker(A_{k}^{[j]})$, where the last step follows from $ii$) of Lemma~\ref{lemma_Arank}. By Lemma~\ref{lemma_Ah}, we have $\ker((A_{k}^{[j]}+h_{k}A_{{\mathrm{c}}k})^{2})=\ker(A_{k}^{[j]}+h_{k}A_{{\mathrm{c}}k})$. Now, by Proposition 5.5.8 of \cite[p.~323]{Bernstein:2009}, 0 is semisimple.
\end{IEEEproof}

It follows from Lemma~\ref{lemma_semisimple} that for every $j=1,\ldots,q$, 0 is a semisimple eigenvalue of $A_{k}^{[j]}+h_{k}A_{{\mathrm{c}}k}$ defined in Lemma~\ref{lemma_Ah}, where $\mu_{k},\kappa_{k},\eta_{k},h_{k}\geq0$, if and only if $\kappa_{k}\neq0$, $k\in\overline{\mathbb{Z}}_{+}$.
To proceed, let $\mathbb{C}^{n}$ (respectively $\mathbb{C}^{m\times n}$) denote the set of complex vectors (respectively matrices). Using Lemmas \ref{lemma_EW}--\ref{lemma_semisimple}, one can show the following complete result about the nonzero eigenvalue and eigenspace structures of $A_{k}^{[j]}+h_{k}A_{{\mathrm{c}}k}$. 

\begin{lemma}\label{lemma_A}
Consider the (possibly infinitely many) matrices $A_{k}^{[j]}+h_{k}A_{{\mathrm{c}}k}$, $j=1,\ldots,q$, $k=0,1,2,\ldots$, defined by (\ref{Amatrix}) in Lemma~\ref{lemma_Arank} and (\ref{Ac}) in Lemma~\ref{lemma_semisimple}, where $\mu_{k},\kappa_{k},\eta_{k}\geq0$ and  $h_{k}>0$, $k\in\overline{\mathbb{Z}}_{+}$. 
\begin{itemize}
\item[$i$)] Then for every $j=1,\ldots,q$, ${\mathrm{spec}}(A_{k}^{[j]}+h_{k}A_{{\mathrm{c}}k})\subseteq\{0,-\kappa_{k},-\frac{\kappa_{k}(1+h_{k})}{2}\pm\frac{1}{2}\sqrt{\kappa_{k}^{2}(1+h_{k})^{2}-4\kappa_{k}},\lambda\in\mathbb{C}:\forall \frac{\lambda^{2}+\kappa_{k} h_{k}\lambda+\kappa_{k}}{\eta_{k}\lambda+\mu_{k} h_{k}\lambda+\mu_{k}}\in{\mathrm{spec}}(-L_{k})\}=\{0,-\kappa_{k},-\frac{\kappa_{k}(1+h_{k})}{2}\pm\frac{1}{2}\sqrt{\kappa_{k}^{2}(1+h_{k})^{2}-4\kappa_{k}},-\frac{\kappa_{k}h_{k}}{2}\pm\frac{1}{2}\sqrt{\kappa_{k}^{2}h_{k}^{2}-4\kappa_{k}},\lambda\\\in\mathbb{C}:\forall \frac{\lambda^{2}+\kappa_{k} h_{k}\lambda+\kappa_{k}}{\eta_{k}\lambda+\mu_{k} h_{k}\lambda+\mu_{k}}\in{\mathrm{spec}}(-L_{k})\backslash\{0\}\}$. 
\item[$ii$)] If $1\not\in{\mathrm{spec}}((\frac{\mu_{k}}{\lambda_{1,2}\kappa_{k}}+\frac{\mu_{k} h_{k}}{\kappa_{k}}+\frac{\eta_{k}}{\kappa_{k}})L_{k})$, then $\lambda_{1,2}=-\frac{\kappa_{k}(1+h_{k})}{2}\pm\frac{1}{2}\sqrt{\kappa_{k}^{2}(1+h_{k})^{2}-4\kappa_{k}}$ are the eigenvalues of $A_{k}^{[j]}+h_{k}A_{{\mathrm{c}}k}$. The corresponding eigenspace is given by 
\begin{eqnarray}\label{egns1}
&&\hspace{-2em}\ker\Big(A_{k}^{[j]}+h_{k}A_{{\mathrm{c}}k}-\lambda_{1,2} I_{2nq+n}\Big)\nonumber\\
&&\hspace{-2em}=\Big\{\Big[\frac{1+h_{k}\lambda_{1,2}^{*}}{\lambda_{1,2}^{*}}\sum_{l=0}^{q-1-{\mathrm{rank}}(L_{k})}\sum_{i=1}^{n}\omega_{li}(\textbf{w}_{l}\otimes\textbf{e}_{i})^{\mathrm{T}},\sum_{l=0}^{q-1-{\mathrm{rank}}(L_{k})}\sum_{i=1}^{n}\omega_{li}(\textbf{w}_{l}\otimes\textbf{e}_{i})^{\mathrm{T}},\nonumber\\
&&-\sum_{l=0}^{q-1-{\mathrm{rank}}(L_{k})}\sum_{i=1}^{n}\omega_{li}w_{lj}\textbf{e}_{i}^{\mathrm{T}}\Big]^{*}:\forall\omega_{li}\in\mathbb{C},i=1,\ldots,n,l=0,1,\ldots,q-1-{\mathrm{rank}}(L_{k})\Big\},
\end{eqnarray} where $\textbf{x}^{*}$ denotes the complex conjugate transpose of $\textbf{x}\in\mathbb{C}^{n}$.
\item[$iii$)] If $1\in{\mathrm{spec}}((\frac{\mu_{k}}{\lambda_{1,2}\kappa_{k}}+\frac{\mu_{k} h_{k}}{\kappa_{k}}+\frac{\eta_{k}}{\kappa_{k}})L_{k})$, and $h_{k}\kappa_{k}\neq 1$, then $\lambda_{1,2}=-\frac{\kappa_{k}(1+h_{k})}{2}\pm\frac{1}{2}\sqrt{\kappa_{k}^{2}(1+h_{k})^{2}-4\kappa_{k}}$ are the eigenvalues of $A_{k}^{[j]}+h_{k}A_{{\mathrm{c}}k}$. The corresponding eigenspace is given by 
\begin{eqnarray}\label{egns2}
&&\ker\Big(A_{k}^{[j]}+h_{k}A_{{\mathrm{c}}k}-\lambda_{1,2} I_{2nq+n}\Big)\nonumber\\
&&=\Big\{\Big[\frac{1+h_{k}\lambda_{1,2}^{*}}{\lambda_{1,2}^{*}}\sum_{i=1}^{n}\sum_{l=1}^{q}\varpi_{li}((\textbf{g}_{l}-G_{k}^{+}G_{k}\textbf{g}_{l})\otimes \textbf{e}_{i})^{\mathrm{T}}-\frac{1+h_{k}\lambda_{1,2}^{*}}{\kappa_{k}\lambda_{1,2}^{*}}\sum_{i=1}^{n}\omega_{0i}(\textbf{1}_{q\times 1}\otimes\textbf{e}_{i})^{\mathrm{T}},\nonumber\\
&&\sum_{i=1}^{n}\sum_{l=1}^{q}\varpi_{li}((\textbf{g}_{l}-G_{k}^{+}G_{k}\textbf{g}_{l})\otimes \textbf{e}_{i})^{\mathrm{T}}-\frac{1}{\kappa_{k}}\sum_{i=1}^{n}\omega_{0i}(\textbf{1}_{q\times 1}\otimes\textbf{e}_{i})^{\mathrm{T}},\nonumber\\
&&\frac{\kappa_{k}+\kappa_{k}h_{k}\lambda_{1,2}^{*}}{\lambda_{1,2}^{*}(\lambda_{1,2}^{*}+\kappa_{k})}\sum_{i=1}^{n}\sum_{l=1}^{q}\varpi_{li}(\textbf{g}_{j}^{\mathrm{T}}\textbf{g}_{l}-\textbf{g}_{j}^{\mathrm{T}}G_{k}^{+}G_{k}\textbf{g}_{l})\textbf{e}_{i}^{\mathrm{T}}-\frac{1+h_{k}\lambda_{1,2}^{*}}{\lambda_{1,2}^{*}(\lambda_{1,2}^{*}+\kappa_{k})}\sum_{i=1}^{n}\omega_{0i}\textbf{e}_{i}^{\mathrm{T}}\Big]^{*}:\nonumber\\
&&\forall\omega_{0i}\in\mathbb{C},\forall\varpi_{li}\in\mathbb{C},i=1,\ldots,n,l=1,\ldots,q\Big\},
\end{eqnarray} where $G_{k}=(\frac{\mu_{k}}{\lambda_{1,2}}+\mu_{k} h_{k}+\eta_{k})L_{k}-\kappa_{k}I_{q}$ and $A^{+}$ denotes the Moore-Penrose generalized inverse of $A$.
\item[$iv$)] If $\frac{\kappa_{k}}{\lambda_{4}}+\lambda_{4}+\kappa_{k} h_{k}\neq0$, $\lambda_{4}\neq-\kappa_{k}$, $\frac{\mu_{k}}{\lambda_{4}}+\mu_{k} h_{k}+\eta_{k}\neq0$, and $\frac{\lambda_{4}^{2}+\kappa_{k} h_{k}\lambda_{4}+\kappa_{k}}{\eta_{k}\lambda+\mu_{k} h_{k}\lambda_{4}+\mu_{k}}\in{\mathrm{spec}}(-L_{k})$, then $\lambda_{4}$ are the eigenvalues of $A_{k}^{[j]}+h_{k}A_{{\mathrm{c}}k}$. The corresponding eigenspace is given by 
\begin{eqnarray}\label{egns3}
&&\hspace{-3em}\ker\Big(A_{k}^{[j]}+h_{k}A_{{\mathrm{c}}k}-\lambda_{4} I_{2nq+n}\Big)\nonumber\\
&&\hspace{-3em}=\Big\{\Big[\frac{1+h_{k}\lambda_{4}^{*}}{\lambda_{4}^{*}}\sum_{i=1}^{n}\sum_{l=1}^{q}\varpi_{li}\Big(\textbf{g}_{l}-F_{k}^{+}F_{k}\textbf{g}_{l}+\frac{\kappa_{k}^{2}(1+h_{k}\lambda_{4})}{\lambda_{4}(\lambda_{4}+\kappa_{k})}(\textbf{g}_{j}^{\mathrm{T}}F_{k}\textbf{g}_{l})F_{k}^{+}\psi_{k}-\frac{\kappa_{k}^{2}(1+h_{k}\lambda_{4})}{\lambda_{4}(\lambda_{4}+\kappa_{k})}(\textbf{g}_{j}^{\mathrm{T}}\textbf{g}_{l})\psi_{k}\Big)^{*}\otimes \textbf{e}_{i}^{\mathrm{T}},\nonumber\\
&&\hspace{-3em}\sum_{i=1}^{n}\sum_{l=1}^{q}\varpi_{li}\Big(\textbf{g}_{l}-F_{k}^{+}F_{k}\textbf{g}_{l}+\frac{\kappa_{k}^{2}(1+h_{k}\lambda_{4})}{\lambda_{4}(\lambda_{4}+\kappa_{k})}(\textbf{g}_{j}^{\mathrm{T}}F_{k}\textbf{g}_{l})F_{k}^{+}\psi_{k}-\frac{\kappa_{k}^{2}(1+h_{k}\lambda_{4})}{\lambda_{4}(\lambda_{4}+\kappa_{k})}(\textbf{g}_{j}^{\mathrm{T}}\textbf{g}_{l})\psi_{k}\Big)^{*}\otimes \textbf{e}_{i}^{\mathrm{T}},\nonumber\\
&&\hspace{-3em}\frac{\kappa_{k}+\kappa_{k}h_{k}\lambda_{4}^{*}}{\lambda_{4}^{*}(\lambda_{4}^{*}+\kappa_{k})}\sum_{i=1}^{n}\sum_{l=1}^{q}\varpi_{li}\Big(\textbf{g}_{j}^{\mathrm{T}}\textbf{g}_{l}-\textbf{g}_{j}^{\mathrm{T}}F_{k}^{+}F_{k}\textbf{g}_{l}+\frac{\kappa_{k}^{2}(1+h_{k}\lambda_{4})}{\lambda_{4}(\lambda_{4}+\kappa_{k})}(\textbf{g}_{j}^{\mathrm{T}}F_{k}\textbf{g}_{l})\textbf{g}_{j}^{\mathrm{T}}F_{k}^{+}\psi_{k}\nonumber\\
&&\hspace{-3em}-\frac{\kappa_{k}^{2}(1+h_{k}\lambda_{4})}{\lambda_{4}(\lambda_{4}+\kappa_{k})}(\textbf{g}_{j}^{\mathrm{T}}\textbf{g}_{l})\textbf{g}_{j}^{\mathrm{T}}\psi_{k}\Big)^{*}\otimes \textbf{e}_{i}^{\mathrm{T}}\Big]^{*}:\varpi_{li}\in\mathbb{C},i=1,\ldots,n,l=1,\ldots,q\Big\},
\end{eqnarray} where $F_{k}=(\frac{\mu_{k}}{\lambda_{4}}+\mu_{k} h_{k}+\eta_{k})L_{k}+(\frac{\kappa_{k}}{\lambda_{4}}+\lambda_{4}+\kappa_{k} h_{k})I_{q}$ and
\begin{eqnarray}\label{psik}
\psi_{k}=\small\left\{\begin{array}{ll}
(\frac{\kappa_{k}^{2}(1+h_{k}\lambda_{4})}{\lambda_{4}(\lambda_{4}+\kappa_{k})}\textbf{g}_{j}^{\mathrm{T}}-\frac{\kappa_{k}^{2}(1+h_{k}\lambda_{4})}{\lambda_{4}(\lambda_{4}+\kappa_{k})}\textbf{g}_{j}^{\mathrm{T}}F_{k}^{+}F_{k})^{+}, & \frac{\kappa_{k}^{2}(1+h_{k}\lambda_{4}^{*})}{\lambda_{4}^{*}(\lambda_{4}^{*}+\kappa_{k})}\textbf{g}_{j}\neq \frac{\kappa_{k}^{2}(1+h_{k}\lambda_{4}^{*})}{\lambda_{4}^{*}(\lambda_{4}^{*}+\kappa_{k})}F_{k}^{+}F_{k}\textbf{g}_{j},\\
\frac{\kappa_{k}^{2}(1+h_{k}\lambda_{4})}{\lambda_{4}(\lambda_{4}+\kappa_{k})}(1+|\frac{\kappa_{k}^{2}(1+h_{k}\lambda_{4})}{\lambda_{4}(\lambda_{4}+\kappa_{k})}|^{2}\textbf{g}_{j}^{\mathrm{T}}(F_{k}^{\mathrm{T}}F_{k})^{+}\textbf{g}_{j})^{-1}(F_{k}^{\mathrm{T}}F_{k})^{+}\textbf{g}_{j}, & \frac{\kappa_{k}^{2}(1+h_{k}\lambda_{4}^{*})}{\lambda_{4}^{*}(\lambda_{4}^{*}+\kappa_{k})}\textbf{g}_{j}= \frac{\kappa_{k}^{2}(1+h_{k}\lambda_{4}^{*})}{\lambda_{4}^{*}(\lambda_{4}^{*}+\kappa_{k})}F_{k}^{+}F_{k}\textbf{g}_{j}. \\
\end{array}\right.
\end{eqnarray}
\item[$v$)] If $\frac{\mu_{k}}{\lambda_{5,6}}+\mu_{k} h_{k}+\eta_{k}\neq0$, $\lambda_{5,6}\neq-\kappa_{k}$, and $\frac{\kappa_{k}}{\lambda_{5,6}}+\lambda_{5,6}+\kappa_{k} h_{k}=0$, then $\lambda_{5,6}=-\frac{\kappa_{k}h_{k}}{2}\pm\frac{1}{2}\sqrt{\kappa_{k}^{2}h_{k}^{2}-4\kappa_{k}}$ are the eigenvalues of $A_{k}^{[j]}+h_{k}A_{{\mathrm{c}}k}$. The corresponding eigenspace is given by the form (\ref{egns3}) with $\lambda_{4}$ being replaced by $\lambda_{5,6}$.
\item[$vi$)] If $\frac{\mu_{k}}{\lambda_{5,6}}+\mu_{k} h_{k}+\eta_{k}=0$, $\lambda_{5,6}\neq-\kappa_{k}$, $\mu_{k}=0$, and $\frac{\kappa_{k}}{\lambda_{5,6}}+\lambda_{5,6}+\kappa_{k} h_{k}=0$, then $\lambda_{5,6}$ are the eigenvalues of $A_{k}^{[j]}+h_{k}A_{{\mathrm{c}}k}$. The corresponding eigenspace is given by
\begin{eqnarray}\label{egns4}
&&\hspace{-1em}\ker\Big(A_{k}^{[j]}+h_{k}A_{{\mathrm{c}}k}-\lambda_{5,6} I_{2nq+n}\Big)\nonumber\\
&&\hspace{-1em}=\Big\{\Big[\frac{1+h_{k}\lambda_{5,6}^{*}}{\lambda_{5,6}^{*}}\sum_{i=1}^{n}\sum_{l=1}^{q}\varpi_{li}(\textbf{g}_{l}-(\textbf{g}_{j}^{\mathrm{T}}\textbf{g}_{l})\textbf{g}_{j})^{\mathrm{T}}\otimes \textbf{e}_{i}^{\mathrm{T}},\sum_{i=1}^{n}\sum_{l=1}^{q}\varpi_{li}(\textbf{g}_{l}-(\textbf{g}_{j}^{\mathrm{T}}\textbf{g}_{l})\textbf{g}_{j})^{\mathrm{T}}\otimes \textbf{e}_{i}^{\mathrm{T}},\textbf{0}_{1\times n}\Big]^{*}:\nonumber\\
&&\hspace{-1em}\varpi_{li}\in\mathbb{C},i=1,\ldots,n,l=1,\ldots,q\Big\}.
\end{eqnarray}
\item[$vii$)] If $1\in{\mathrm{spec}}(\frac{\eta_{k}}{\kappa_{k}}L_{k})$ and $\kappa_{k}h_{k}=1$, then $\lambda_{3}=-\kappa_{k}$ is an eigenvalue of $A_{k}^{[j]}+h_{k}A_{{\mathrm{c}}k}$. The corresponding eigenspace is given by
\begin{eqnarray}\label{egns5}
&&\ker\Big(A_{k}^{[j]}+h_{k}A_{{\mathrm{c}}k}-\lambda_{3} I_{2nq+n}\Big)\nonumber\\
&&=\Big\{\Big[\textbf{0}_{1\times nq},\sum_{i=1}^{n}\sum_{l=1}^{q}\alpha_{li}(\textbf{g}_{l}\otimes\textbf{e}_{i})^{\mathrm{T}},\sum_{i=1}^{n}\sum_{l=1}^{q}\frac{\eta_{k}}{\kappa_{k}}\alpha_{li}(L_{k}\textbf{g}_{l}\otimes\textbf{e}_{i})^{\mathrm{T}}-\sum_{i=1}^{n}\sum_{l=1}^{q}\alpha_{li}(\textbf{g}_{l}\otimes\textbf{e}_{i})^{\mathrm{T}}\Big]^{*}:\nonumber\\
&&\forall\alpha_{li}\in\mathbb{C},i=1,\ldots,n,l=1,\ldots,q\Big\}.
\end{eqnarray}
\item[$viii$)] If $\frac{\mu_{k}}{\kappa_{k}}(\kappa_{k} h_{k}-1)+\eta_{k}=0$ and $h_{k}=1+\frac{1}{\kappa_{k}}$, then $\lambda_{3}=-\kappa_{k}$ is an eigenvalue of $A_{k}^{[j]}+h_{k}A_{{\mathrm{c}}k}$. The corresponding eigenspace is given by
\begin{eqnarray}\label{egns6}
&&\ker\Big(A_{k}^{[j]}+h_{k}A_{{\mathrm{c}}k}-\lambda_{3} I_{2nq+n}\Big)=\Big\{\Big[\textbf{0}_{1\times nq},\sum_{i=1}^{n}\sum_{l=1}^{q}\alpha_{li}(\textbf{g}_{l}-(\textbf{g}_{j}^{\mathrm{T}}\textbf{g}_{l})\textbf{g}_{j})^{\mathrm{T}}\otimes \textbf{e}_{i}^{\mathrm{T}},\textbf{0}_{1\times n}\Big]^{*}:\nonumber\\
&&\forall\alpha_{li}\in\mathbb{C},i=1,\ldots,n,l=1,\ldots,q\Big\}.
\end{eqnarray}
\item[$ix$)] If $1\in{\mathrm{spec}}(\frac{\mu_{k}+\eta_{k}}{\kappa_{k}}L_{k})$ and $h_{k}=1+\frac{1}{\kappa_{k}}$, then $\lambda_{3}=-\kappa_{k}$ is an eigenvalue of $A_{k}^{[j]}+h_{k}A_{{\mathrm{c}}k}$. The corresponding eigenspace is given by
\begin{eqnarray}\label{egns7}
&&\ker\Big(A_{k}^{[j]}+h_{k}A_{{\mathrm{c}}k}-\lambda_{3} I_{2nq+n}\Big)\nonumber\\
&&=\Big\{\Big[\textbf{0}_{1\times nq},\frac{\kappa_{k}}{\mu_{k}+\eta_{k}}\sum_{i=1}^{n}\beta_{i}(L_{k}^{+}\textbf{1}_{q\times 1}\otimes\textbf{e}_{i})^{\mathrm{T}}-\frac{\kappa_{k}}{\mu_{k}+\eta_{k}}\sum_{i=1}^{n}\beta_{i}(L_{k}^{+}\varphi_{k}\otimes\textbf{e}_{i})^{\mathrm{T}}\nonumber\\
&&+\sum_{l=1}^{q}\sum_{i=1}^{n}\gamma_{li}(\textbf{g}_{l}-L_{k}^{+}L_{k}\textbf{g}_{l}+(\textbf{g}_{j}^{\mathrm{T}}L_{k}\textbf{g}_{l})L_{k}^{+}\varphi_{k}-(\textbf{g}_{j}^{\mathrm{T}}\textbf{g}_{l})\varphi_{k})^{\mathrm{T}}\otimes \textbf{e}_{i}^{\mathrm{T}},\sum_{i=1}^{n}\beta_{i}\textbf{e}_{i}^{\mathrm{T}}\Big]^{*}:\nonumber\\
&&\beta_{i}\in\mathbb{C},\gamma_{li}\in\mathbb{C},i=1,\ldots,n,l=1,\ldots,q\Big\},
\end{eqnarray} where 
\begin{eqnarray}\label{varphik}
\varphi_{k}=\small\left\{\begin{array}{ll}
(\textbf{g}_{j}^{\mathrm{T}}-\textbf{g}_{j}^{\mathrm{T}}L_{k}^{+}L_{k})^{+}, & \textbf{g}_{j}\neq L_{k}^{+}L_{k}\textbf{g}_{j},\\
(1+\textbf{g}_{j}^{\mathrm{T}}(L_{k}^{\mathrm{T}}L_{k})^{+}\textbf{g}_{j})^{-1}(L_{k}^{\mathrm{T}}L_{k})^{+}\textbf{g}_{j}, & \textbf{g}_{j}= L_{k}^{+}L_{k}\textbf{g}_{j}. \\
\end{array}\right.
\end{eqnarray}
\item[$x$)] If $1\in{\mathrm{spec}}(\frac{\mu_{k}(\kappa_{k}h_{k}-1)+\eta_{k}\kappa_{k}}{\kappa_{k}(-\kappa_{k} h_{k}+1+\kappa_{k})}L_{k})$ and $\kappa_{k}h_{k}\neq 1$, then $\lambda_{3}=-\kappa_{k}$ is an eigenvalue of $A_{k}^{[j]}+h_{k}A_{{\mathrm{c}}k}$. The corresponding eigenspace is given by
\begin{eqnarray}\label{egns8}
&&\ker\Big(A_{k}^{[j]}+h_{k}A_{{\mathrm{c}}k}-\lambda_{3} I_{2nq+n}\Big)\nonumber\\
&&=\Big\{\Big[\textbf{0}_{1\times nq},\sum_{i=1}^{n}\sum_{l=1}^{q}\varpi_{li}\Big(\textbf{g}_{l}-M_{k}^{+}M_{k}\textbf{g}_{l}+(\textbf{g}_{j}^{\mathrm{T}}M_{k}\textbf{g}_{l})M_{k}^{+}\phi_{k}-(\textbf{g}_{j}^{\mathrm{T}}\textbf{g}_{l})\phi_{k}\Big)^{\mathrm{T}}\otimes\textbf{e}_{i}^{\mathrm{T}},\textbf{0}_{1\times n}\Big]^{*}:\nonumber\\
&&\varpi_{li}\in\mathbb{C},i=1,\ldots,n,l=1,\ldots,q\Big\},
\end{eqnarray} where $M_{k}=(\frac{\mu_{k}}{\kappa_{k}}(\kappa_{k} h_{k}-1)+\eta_{k})L_{k}+(\kappa_{k} h_{k}-1-\kappa_{k})I_{q}$ and 
\begin{eqnarray}\label{phik}
\phi_{k}=\small\left\{\begin{array}{ll}
(\textbf{g}_{j}^{\mathrm{T}}-\textbf{g}_{j}^{\mathrm{T}}M_{k}^{+}M_{k})^{+}, & \textbf{g}_{j}\neq M_{k}^{+}M_{k}\textbf{g}_{j},\\
(1+\textbf{g}_{j}^{\mathrm{T}}(M_{k}^{\mathrm{T}}M_{k})^{+}\textbf{g}_{j})^{-1}(M_{k}^{\mathrm{T}}M_{k})^{+}\textbf{g}_{j}, & \textbf{g}_{j}= M_{k}^{+}M_{k}\textbf{g}_{j}. \\
\end{array}\right.
\end{eqnarray}
\end{itemize}
\end{lemma}

\begin{IEEEproof}
For a fixed $j\in\{1,\ldots,q\}$ and a fixed $k\in\overline{\mathbb{Z}}_{+}$, let $\textbf{x}\in\mathbb{C}^{2nq+n}$ be an eigenvector of the
corresponding eigenvalue $\lambda\in\mathbb{C}$ for $A_{k}^{[j]}+h_{k}A_{{\mathrm{c}}k}$. We partition
$\textbf{x}$ into
$\textbf{x}=[\textbf{x}_{1}^{*},\textbf{x}_{2}^{*},\textbf{x}_{3}^{*}]^{*}\neq\textbf{0}_{(2nq+n)\times 1}$,
where $\textbf{x}_{1},\textbf{x}_{2}\in\mathbb{C}^{nq}$, and $\textbf{x}_{3}\in\mathbb{C}^{n}$. It follows from
$(A_{k}^{[j]}+h_{k}A_{{\mathrm{c}}k})\textbf{x}=\lambda\textbf{x}$ that
\begin{eqnarray}
h_{k}(-\mu_{k} L_{k}\otimes I_{n}-\kappa_{k} I_{nq})\textbf{x}_{1}+h_{k}(-\eta_{k} L_{k}\otimes I_{n})\textbf{x}_{2}+\textbf{x}_{2}+h_{k}(\kappa_{k} \textbf{1}_{q\times 1}\otimes I_{n})\textbf{x}_{3}=\lambda\textbf{x}_{1},\label{Aeig_1}\\
(-\mu_{k} L_{k}\otimes I_{n}-\kappa_{k} I_{nq})\textbf{x}_{1}+(-\eta_{k} L_{k}\otimes I_{n})\textbf{x}_{2}+(\kappa_{k} \textbf{1}_{q\times 1}\otimes I_{n})\textbf{x}_{3}=\lambda \textbf{x}_{2},\label{x2}\\
\kappa_{k} E_{n\times nq}^{[j]}\textbf{x}_{1}-\kappa_{k}\textbf{x}_{3}=\lambda \textbf{x}_{3}.\label{Aeig_3}
\end{eqnarray} Note that it follows from Lemma~\ref{lemma_semisimple} that $A_{k}^{[j]}+h_{k}A_{{\mathrm{c}}k}$ has an eigenvalue 0. Now we assume that $\lambda\neq0$. 

Substituting (\ref{x2}) into (\ref{Aeig_1}) yields
$\textbf{x}_{1}=\frac{1+h_{k}\lambda}{\lambda}\textbf{x}_{2}$. Replacing $\textbf{x}_{1}$ in (\ref{x2}) and (\ref{Aeig_3}) with $\textbf{x}_{1}=\frac{1+h_{k}\lambda}{\lambda}\textbf{x}_{2}$ yields
\begin{eqnarray}
-\Big[\Big(\frac{\mu_{k}}{\lambda}+\mu_{k} h_{k}+\eta_{k}\Big)(L_{k}\otimes I_{n})+\Big(\frac{\kappa_{k}}{\lambda}+\lambda+\kappa_{k} h_{k}\Big)I_{nq}\Big]\textbf{x}_{2}+\kappa_{k}(\textbf{1}_{q\times 1}\otimes I_{n})\textbf{x}_{3}=\textbf{0}_{nq\times 1},\label{hx2}\\
\Big(\frac{\kappa_{k}}{\lambda}+\kappa_{k} h_{k}\Big)E_{n\times nq}^{[j]}\textbf{x}_{2}-(\lambda+\kappa_{k})\textbf{x}_{3}=\textbf{0}_{n\times 1}.\label{hx3}
\end{eqnarray} Clearly $\textbf{x}_{2}\neq\textbf{0}_{nq\times 1}$. Thus, (\ref{hx2}) and (\ref{hx3}) have nontrivial solutions if and only if
\begin{eqnarray}\label{detcon}
\det\small\left[\begin{array}{cc}
\Big(\frac{\mu_{k}}{\lambda}+\mu_{k} h_{k}+\eta_{k}\Big)(L_{k}\otimes I_{n})+\Big(\frac{\kappa_{k}}{\lambda}+\lambda+\kappa_{k} h_{k}\Big)I_{nq} & -\kappa_{k}(\textbf{1}_{q\times 1}\otimes I_{n})\\
\Big(\frac{\kappa_{k}}{\lambda}+\kappa_{k} h_{k}\Big)E_{n\times nq}^{[j]} & -(\lambda+\kappa_{k})I_{n}
\end{array}\right]=0.
\end{eqnarray} 

If $\det\Big[\Big(\frac{\mu_{k}}{\lambda}+\mu_{k} h_{k}+\eta_{k}\Big)(L_{k}\otimes I_{n})+\Big(\frac{\kappa_{k}}{\lambda}+\lambda+\kappa_{k} h_{k}\Big)I_{nq}\Big]\neq0$, then pre-multiplying $-L_{k}\otimes I_{n}$ on both sides of (\ref{hx2}) yields
\begin{eqnarray*}
\Big[\Big(\frac{\mu_{k}}{\lambda}+\mu_{k} h_{k}+\eta_{k}\Big)(L_{k}\otimes I_{n})+\Big(\frac{\kappa_{k}}{\lambda}+\lambda+\kappa_{k} h_{k}\Big)I_{nq}\Big](L_{k}\otimes I_{n})\textbf{x}_{2}=\textbf{0}_{nq\times 1},
\end{eqnarray*} which implies that $(L_{k}\otimes I_{n})\textbf{x}_{2}=\textbf{0}_{nq\times 1}$. Now following the similar arguments as in the proof of Case 2 of $ii$) in Lemma~\ref{lemma_Arank}, we have $\textbf{x}_{2}=\sum_{l=0}^{q-1-{\mathrm{rank}}(L_{k})}\sum_{i=1}^{n}\varpi_{li}(\textbf{w}_{l}\otimes\textbf{e}_{i})$, where $\varpi_{li}\in\mathbb{C}$ and not all $\varpi_{li}$ are zero. Substituting this expression of $\textbf{x}_{2}$ into (\ref{hx2}) and (\ref{hx3}) by using $iii$) of Lemma~\ref{lemma_EW} yields
\begin{eqnarray} 
\kappa_{k}\textbf{x}_{3}=\left(\frac{\kappa_{k}}{\lambda}+\lambda+\kappa_{k}h_{k}\right)\sum_{l=0}^{q-1-{\mathrm{rank}}(L_{k})}\sum_{i=1}^{n}\varpi_{li}w_{lj}\textbf{e}_{i}.\label{x3e1}\\
(\lambda+\kappa_{k})\textbf{x}_{3}=\left(\frac{\kappa_{k}}{\lambda}+\kappa_{k} h_{k}\right)\sum_{l=0}^{q-1-{\mathrm{rank}}(L_{k})}\sum_{i=1}^{n}\varpi_{li}w_{lj}\textbf{e}_{i}.\label{x3e2}
\end{eqnarray}
Furthermore, substituting (\ref{x3e1}) into (\ref{x3e2}) yields
\begin{eqnarray*}
\lambda\textbf{x}_{3}=-\lambda\sum_{l=0}^{q-1-{\mathrm{rank}}(L_{k})}\sum_{i=1}^{n}\varpi_{li}w_{lj}\textbf{e}_{i},
\end{eqnarray*} which implies that $\textbf{x}_{3}=-\sum_{l=0}^{q-1-{\mathrm{rank}}(L_{k})}\sum_{i=1}^{n}\varpi_{li}w_{lj}\textbf{e}_{i}$ since $\lambda\neq0$. Finally, substituting the obtained expressions for $\textbf{x}_{2}$ and $\textbf{x}_{3}$ into (\ref{hx3}), or substituting the obtained expression for $\textbf{x}_{3}$ into either (\ref{x3e1}) or (\ref{x3e2}) yields
\begin{eqnarray}\label{eqn_ei} 
\left(\frac{\kappa_{k}}{\lambda}+\kappa_{k} h_{k}+\lambda+\kappa_{k}\right)\sum_{l=0}^{q-1-{\mathrm{rank}}(L_{k})}\sum_{i=1}^{n}\varpi_{li}w_{lj}\textbf{e}_{i}=-\left(\frac{\kappa_{k}}{\lambda}+\kappa_{k} h_{k}+\lambda+\kappa_{k}\right)\textbf{x}_{3}=\textbf{0}_{n\times 1}.
\end{eqnarray}
In this case, (\ref{hx2}) and (\ref{hx3}) have nontrivial solutions if and only if (\ref{eqn_ei}) holds, which implies that  $\frac{\kappa_{k}}{\lambda}+\kappa_{k} h_{k}+\lambda+\kappa_{k}=0$ since $\textbf{x}_{3}\neq\textbf{0}_{n\times 1}$, and hence, $\kappa_{k}\neq0$. Let $\lambda_{1,2}$ denote the two solutions to $\frac{\kappa_{k}}{\lambda}+\kappa_{k} h_{k}+\lambda+\kappa_{k}=0$. Then 
\begin{eqnarray}\label{lambda1}
\lambda_{1,2}=-\frac{\kappa_{k}(1+h_{k})}{2}\pm\frac{1}{2}\sqrt{\kappa_{k}^{2}(1+h_{k})^{2}-4\kappa_{k}}.
\end{eqnarray} In this case, note that
\begin{eqnarray}\label{detlambda}
&&\det\Big[\Big(\frac{\mu_{k}}{\lambda_{1,2}}+\mu_{k} h_{k}+\eta_{k}\Big)(L_{k}\otimes I_{n})+\Big(\frac{\kappa_{k}}{\lambda_{1,2}}+\lambda_{1,2}+\kappa_{k} h_{k}\Big)I_{nq}\Big]\nonumber\\
&&=\det\Big[\Big(\frac{\mu_{k}}{\lambda_{1,2}}+\mu_{k} h_{k}+\eta_{k}\Big)(L_{k}\otimes I_{n})-\kappa_{k}I_{nq}\Big]\nonumber\\
&&=\kappa_{k}^{nq}\det\Big[\Big(\frac{\mu_{k}}{\lambda_{1,2}\kappa_{k}}+\frac{\mu_{k} h_{k}}{\kappa_{k}}+\frac{\eta_{k}}{\kappa_{k}}\Big)(L_{k}\otimes I_{n})-I_{nq}\Big].
\end{eqnarray} Hence, $\det\Big[\Big(\frac{\mu_{k}}{\lambda_{1,2}}+\mu_{k} h_{k}+\eta_{k}\Big)(L_{k}\otimes I_{n})+\Big(\frac{\kappa_{k}}{\lambda_{1,2}}+\lambda_{1,2}+\kappa_{k} h_{k}\Big)I_{nq}\Big]\neq0$ if and only if $1\not\in{\mathrm{spec}}((\frac{\mu_{k}}{\lambda_{1,2}\kappa_{k}}+\frac{\mu_{k} h_{k}}{\kappa_{k}}+\frac{\eta_{k}}{\kappa_{k}})L_{k})$. Thus, if 
$1\not\in{\mathrm{spec}}((\frac{\mu_{k}}{\lambda_{1,2}\kappa_{k}}+\frac{\mu_{k} h_{k}}{\kappa_{k}}+\frac{\eta_{k}}{\kappa_{k}})L_{k})$, then $\lambda_{1,2}$ given by (\ref{lambda1}) are indeed the eigenvalues of $A_{k}^{[j]}+h_{k}A_{{\mathrm{c}}k}$ and the corresponding eigenvectors for $\lambda_{1,2}$ are given by 
\begin{eqnarray}
\textbf{x}&=&\Big[\frac{1+h_{k}\lambda_{1,2}^{*}}{\lambda_{1,2}^{*}}\sum_{l=0}^{q-1-{\mathrm{rank}}(L_{k})}\sum_{i=1}^{n}\varpi_{li}(\textbf{w}_{l}\otimes\textbf{e}_{i})^{\mathrm{T}},\sum_{l=0}^{q-1-{\mathrm{rank}}(L_{k})}\sum_{i=1}^{n}\varpi_{li}(\textbf{w}_{l}\otimes\textbf{e}_{i})^{\mathrm{T}},\nonumber\\
&&-\sum_{l=0}^{q-1-{\mathrm{rank}}(L_{k})}\sum_{i=1}^{n}\varpi_{li}w_{lj}\textbf{e}_{i}^{\mathrm{T}}\Big]^{*},
\end{eqnarray} where $\varpi_{li}\in\mathbb{C}$ and not all of $\varpi_{li}$ are zero. Therefore,
$\ker\Big(A_{k}^{[j]}+h_{k}A_{{\mathrm{c}}k}-\lambda_{1,2} I_{2nq+n}\Big)$ is given by (\ref{egns1}). 

Alternatively, if $\det\Big[\Big(\frac{\mu_{k}}{\lambda}+\mu_{k} h_{k}+\eta_{k}\Big)(L_{k}\otimes I_{n})+\Big(\frac{\kappa_{k}}{\lambda}+\lambda+\kappa_{k} h_{k}\Big)I_{nq}\Big]=0$, then in this case, we consider two additional cases for (\ref{detcon}): 

\textit{Case 1.} If $\lambda\neq-\kappa_{k}$, then it follows from Proposition 2.8.4 of \cite[p.~116]{Bernstein:2009} that (\ref{detcon}) is equivalent to $\det\Big(\Big(\frac{\mu_{k}}{\lambda}+\mu_{k} h_{k}+\eta_{k}\Big)(L_{k}\otimes I_{n})+\Big(\frac{\kappa_{k}}{\lambda}+\lambda+\kappa_{k} h_{k}\Big)I_{nq}-\frac{\kappa_{k}^{2}(1+h_{k}\lambda)}{\lambda(\lambda+\kappa_{k})}W^{[j]}\Big)=0$, which implies that for $\lambda\neq-\kappa_{k}$, the equation
\begin{eqnarray}\label{eqn_v}
\Big(\Big(\frac{\mu_{k}}{\lambda}+\mu_{k} h_{k}+\eta_{k}\Big)(L_{k}\otimes I_{n})+\Big(\frac{\kappa_{k}}{\lambda}+\lambda+\kappa_{k} h_{k}\Big)I_{nq}-\frac{\kappa_{k}^{2}(1+h_{k}\lambda)}{\lambda(\lambda+\kappa_{k})}W^{[j]}\Big)\textbf{v}=\textbf{0}_{nq\times 1}
\end{eqnarray} has nontrivial solutions for $\textbf{v}\in\mathbb{C}^{nq}$. It follows from (\ref{hx2}) and (\ref{hx3}) that solving this $\textbf{v}$ is equivalent to solving $\textbf{x}_{2}$. Again, note that for every $j=1,\ldots,q$, $(L_{k}\otimes I_{n})W^{[j]}=\textbf{0}_{nq\times nq}$. Pre-multiplying $L_{k}\otimes I_{n}$ on both sides of (\ref{eqn_v}) yields $\Big(\Big(\frac{\mu_{k}}{\lambda}+\mu_{k} h_{k}+\eta_{k}\Big)(L_{k}\otimes I_{n})^{2}+\Big(\frac{\kappa_{k}}{\lambda}+\lambda+\kappa_{k} h_{k}\Big)(L_{k}\otimes I_{n})\Big)\textbf{v}=(L_{k}\otimes I_{n})\Big(\Big(\frac{\mu_{k}}{\lambda}+\mu_{k} h_{k}+\eta_{k}\Big)(L_{k}\otimes I_{n})+\Big(\frac{\kappa_{k}}{\lambda}+\lambda+\kappa_{k} h_{k}\Big)I_{nq}\Big)\textbf{v}=\textbf{0}_{nq\times 1}$, which implies that $\Big(\Big(\frac{\mu_{k}}{\lambda}+\mu_{k} h_{k}+\eta_{k}\Big)(L_{k}\otimes I_{n})+\Big(\frac{\kappa_{k}}{\lambda}+\lambda+\kappa_{k} h_{k}\Big)I_{nq}\Big)\textbf{v}\in\ker(L_{k}\otimes I_{n})$. Since $\ker(L_{k}\otimes I_{n})=\bigcup_{l=0}^{q-1-{\mathrm{rank}}(L_{k})}{\mathrm{span}}\{\textbf{w}_{l}\otimes\textbf{e}_{1},\ldots,\textbf{w}_{l}\otimes\textbf{e}_{n}\}$, it follows that
\begin{eqnarray}\label{Az=b}
\Big(\Big(\frac{\mu_{k}}{\lambda}+\mu_{k} h_{k}+\eta_{k}\Big)(L_{k}\otimes I_{n})+\Big(\frac{\kappa_{k}}{\lambda}+\lambda+\kappa_{k} h_{k}\Big)I_{nq}\Big)\textbf{v}=\sum_{i=1}^{n}\sum_{l=0}^{q-1-{\mathrm{rank}}(L_{k})}\omega_{li}\textbf{w}_{l}\otimes\textbf{e}_{i},
\end{eqnarray} where $\omega_{li}\in\mathbb{C}$. Now it follows from (\ref{eqn_v}) and (\ref{Az=b}) that 
\begin{eqnarray}\label{Wv_eqn}
\frac{\kappa_{k}^{2}(1+h_{k}\lambda)}{\lambda(\lambda+\kappa_{k})}W^{[j]}\textbf{v}=\sum_{i=1}^{n}\sum_{l=0}^{q-1-{\mathrm{rank}}(L_{k})}\omega_{li}\textbf{w}_{l}\otimes\textbf{e}_{i}.
\end{eqnarray}

If $\frac{\kappa_{k}}{\lambda}+\lambda+\kappa_{k} h_{k}\neq0$, then (\ref{Az=b}) has a particular solution $\textbf{v}=(\frac{\kappa_{k}}{\lambda}+\lambda+\kappa_{k} h_{k})^{-1}\sum_{i=1}^{n}\sum_{l=0}^{q-1-{\mathrm{rank}}(L_{k})}\omega_{li}\textbf{w}_{l}\otimes\textbf{e}_{i}$. Let $\textbf{w}_{l}=[w_{l1}^{*},\ldots,w_{lq}^{*}]^{*}$. Substituting this particular solution into (\ref{Wv_eqn}), together with $ii$) of Lemma~\ref{lemma_EW}, yields
\begin{eqnarray} &&\sum_{i=1}^{n}\sum_{l=0}^{q-1-{\mathrm{rank}}(L_{k})}\omega_{li}\textbf{w}_{l}\otimes\textbf{e}_{i}-\frac{\kappa_{k}^{2}(1+h_{k}\lambda)}{\lambda(\lambda+\kappa_{k})}W^{[j]}(\frac{\kappa_{k}}{\lambda}+\lambda+\kappa_{k} h_{k})^{-1}\sum_{i=1}^{n}\sum_{l=0}^{q-1-{\mathrm{rank}}(L_{k})}\omega_{li}\textbf{w}_{l}\otimes\textbf{e}_{i}\nonumber\\
&&=\sum_{i=1}^{n}\sum_{l=0}^{q-1-{\mathrm{rank}}(L_{k})}\omega_{li}\textbf{w}_{l}\otimes\textbf{e}_{i}-\frac{\kappa_{k}^{2}(1+h_{k}\lambda)}{(\lambda+\kappa_{k})(\lambda^{2}+\kappa_{k} h_{k}\lambda+\kappa_{k})}\sum_{i=1}^{n}\sum_{l=0}^{q-1-{\mathrm{rank}}(L_{k})}\omega_{li}w_{lj}\textbf{w}_{0}\otimes\textbf{e}_{i}\nonumber\\
&&=\sum_{i=1}^{n}\Big[\omega_{0i}-\frac{\kappa_{k}^{2}(1+h_{k}\lambda)}{(\lambda+\kappa_{k})(\lambda^{2}+\kappa_{k} h_{k}\lambda+\kappa_{k})}\sum_{l=0}^{q-1-{\mathrm{rank}}(L_{k})}\omega_{li}w_{lj}\Big]\textbf{w}_{0}\otimes\textbf{e}_{i}+\sum_{i=1}^{n}\sum_{l=1}^{q-1-{\mathrm{rank}}(L_{k})}\omega_{li}\textbf{w}_{l}\otimes\textbf{e}_{i}\nonumber\\
&&=\textbf{0}_{nq\times 1},
\end{eqnarray} which implies that 
\begin{eqnarray}\label{omega}
\omega_{0i}-\frac{\kappa_{k}^{2}(1+h_{k}\lambda)}{(\lambda+\kappa_{k})(\lambda^{2}+\kappa_{k} h_{k}\lambda+\kappa_{k})}\sum_{l=0}^{q-1-{\mathrm{rank}}(L_{k})}\omega_{li}w_{lj}=0
\end{eqnarray} and $\omega_{\ell i}=0$ for every $i=1,\ldots,n$ and every $\ell=1,\ldots,q-1-{\mathrm{rank}}(L_{k})$. Note that $w_{0j}=1$ for every $j=1,\ldots,q$. Substituting $\omega_{\ell i}=0$ into (\ref{omega}) yields
\begin{eqnarray}
\omega_{0i}-\frac{\kappa_{k}^{2}(1+h_{k}\lambda)}{(\lambda+\kappa_{k})(\lambda^{2}+\kappa_{k} h_{k}\lambda+\kappa_{k})}\omega_{0i}=0,\quad i=1,\ldots,n.
\end{eqnarray}
Then either $1-\frac{\kappa_{k}^{2}(1+h_{k}\lambda)}{(\lambda+\kappa_{k})(\lambda^{2}+\kappa_{k} h_{k}\lambda+\kappa_{k})}=0$ or $\omega_{0i}=0$ for every $i=1,\ldots,n$.

If $\frac{\kappa_{k}^{2}(1+h_{k}\lambda)}{(\lambda+\kappa_{k})(\lambda^{2}+\kappa_{k} h_{k}\lambda+\kappa_{k})}=1$, then $\lambda^{2}+\kappa_{k}(1+h_{k})\lambda+\kappa_{k}=0$. Hence, $\lambda=\lambda_{12}$ where $\lambda_{1,2}$ are given by (\ref{lambda1}). In this case, note that $\frac{\kappa_{k}}{\lambda_{1,2}}+\lambda_{1,2}+\kappa_{k} h_{k}=-\kappa_{k}\neq0$. Then it follows that 
(\ref{detlambda}) holds. Hence,
$\det\Big[\Big(\frac{\mu_{k}}{\lambda_{1,2}}+\mu_{k} h_{k}+\eta_{k}\Big)(L_{k}\otimes I_{n})+\Big(\frac{\kappa_{k}}{\lambda_{1,2}}+\lambda_{1,2}+\kappa_{k} h_{k}\Big)I_{nq}\Big]=0$ if and only if $1\in{\mathrm{spec}}((\frac{\mu_{k}}{\lambda_{1,2}\kappa_{k}}+\frac{\mu_{k} h_{k}}{\kappa_{k}}+\frac{\eta_{k}}{\kappa_{k}})L_{k})$. Furthermore, $\lambda_{1,2}\neq-\kappa_{k}$ if and only if $h_{k}\kappa_{k}\neq 1$. Thus, if 
$1\in{\mathrm{spec}}((\frac{\mu_{k}}{\lambda_{1,2}\kappa_{k}}+\frac{\mu_{k} h_{k}}{\kappa_{k}}+\frac{\eta_{k}}{\kappa_{k}})L_{k})$ and $h_{k}\kappa_{k}\neq 1$, then $\lambda_{1,2}$ given by (\ref{lambda1}) are indeed the eigenvalues of $A_{k}^{[j]}+h_{k}A_{{\mathrm{c}}k}$. In this case, 
(\ref{Az=b}) becomes
\begin{eqnarray}\label{Az=bx}
\Big(\Big(\frac{\mu_{k}}{\lambda_{1,2}}+\mu_{k} h_{k}+\eta_{k}\Big)(L_{k}\otimes I_{n})-\kappa_{k}I_{nq}\Big)\textbf{v}=\sum_{i=1}^{n}\omega_{0i}\textbf{w}_{0}\otimes\textbf{e}_{i}
\end{eqnarray} and a specific solution is given by $\textbf{v}=-\frac{1}{\kappa_{k}}\sum_{i=1}^{n}\omega_{0i}\textbf{w}_{0}\otimes\textbf{e}_{i}$. To find the general solution to (\ref{Az=bx}), let $G_{k}=(\frac{\mu_{k}}{\lambda_{1,2}}+\mu_{k} h_{k}+\eta_{k})L_{k}-\kappa_{k}I_{q}$ and consider 
\begin{eqnarray}\label{Gx0}
(G_{k}\otimes I_{n})\hat{\textbf{v}}=\textbf{0}_{nq\times 1}.
\end{eqnarray} It follows from $vi$) of Proposition 6.1.7 of \cite[p.~400]{Bernstein:2009} and $viii$) of Proposition 6.1.6 of \cite[p.~399]{Bernstein:2009} that the general solution $\hat{\textbf{v}}$ to (\ref{Gx0}) is given by the form
\begin{eqnarray}
\hat{\textbf{v}}&=&\Big[I_{nq}-(G_{k}\otimes I_{n})^{+}(G_{k}\otimes I_{n})\Big]\sum_{i=1}^{n}\sum_{l=1}^{q}\varpi_{li}\textbf{g}_{l}\otimes\textbf{e}_{i}\nonumber\\
&=&\Big[I_{nq}-(G_{k}^{+}\otimes I_{n})(G_{k}\otimes I_{n})\Big]\sum_{i=1}^{n}\sum_{l=1}^{q}\varpi_{li}\textbf{g}_{l}\otimes\textbf{e}_{i}\nonumber\\
&=&\Big[I_{q}\otimes I_{n}-((G_{k}^{+}G_{k})\otimes I_{n})\Big]\sum_{i=1}^{n}\sum_{l=1}^{q}\varpi_{li}\textbf{g}_{l}\otimes\textbf{e}_{i}\nonumber\\
&=&\Big[(I_{q}-G_{k}^{+}G_{k})\otimes I_{n}\Big]\sum_{i=1}^{n}\sum_{l=1}^{q}\varpi_{li}\textbf{g}_{l}\otimes\textbf{e}_{i}\nonumber\\
&=&\sum_{i=1}^{n}\sum_{l=1}^{q}\varpi_{li}(\textbf{g}_{l}-G_{k}^{+}G_{k}\textbf{g}_{l})\otimes \textbf{e}_{i},
\end{eqnarray} where $\varpi_{li}\in\mathbb{C}$, $j=1,\ldots,q$, and we used the facts that $(A\otimes B)^{+}=A^{+}\otimes B^{+}$, $A\otimes B-C\otimes B=(A-C)\otimes B$, and $(A\otimes B)(C\otimes D)=AC\otimes BD$ for compatible matrices $A,B,C,D$. Then the general solution to (\ref{Az=bx}) is given by
\begin{eqnarray}
\textbf{v}&=&\hat{\textbf{v}}-\frac{1}{\kappa_{k}}\sum_{i=1}^{n}\omega_{0i}\textbf{w}_{0}\otimes\textbf{e}_{i}\nonumber\\
&=&\sum_{i=1}^{n}\sum_{l=1}^{q}\varpi_{li}(\textbf{g}_{l}-G_{k}^{+}G_{k}\textbf{g}_{l})\otimes \textbf{e}_{i}-\frac{1}{\kappa_{k}}\sum_{i=1}^{n}\omega_{0i}\textbf{w}_{0}\otimes\textbf{e}_{i},
\end{eqnarray} and hence, $\textbf{x}_{2}=\textbf{v}\neq\textbf{0}_{nq\times 1}$ and  $\textbf{x}_{1}=\frac{1+h_{k}\lambda_{1,2}}{\lambda_{1,2}}\textbf{v}$. Furthermore,  note that $\textbf{g}_{j}^{\mathrm{T}}\textbf{w}_{0}=1$ for every $j=1,\ldots,q$, it follows that  
\begin{eqnarray}
\textbf{x}_{3}&=&\frac{\kappa_{k}+\kappa_{k}h_{k}\lambda_{1,2}}{\lambda_{1,2}(\lambda_{1,2}+\kappa_{k})}E_{n\times nq}^{[j]}\textbf{v}\nonumber\\
&=&\frac{\kappa_{k}+\kappa_{k}h_{k}\lambda_{1,2}}{\lambda_{1,2}(\lambda_{1,2}+\kappa_{k})}(\textbf{g}_{j}^{\mathrm{T}}\otimes I_{n})\textbf{v}\nonumber\\
&=&\frac{\kappa_{k}+\kappa_{k}h_{k}\lambda_{1,2}}{\lambda_{1,2}(\lambda_{1,2}+\kappa_{k})}\sum_{i=1}^{n}\sum_{l=1}^{q}\varpi_{li}(\textbf{g}_{j}^{\mathrm{T}}\otimes I_{n})((\textbf{g}_{l}-G_{k}^{+}G_{k}\textbf{g}_{l})\otimes \textbf{e}_{i})\nonumber\\
&&-\frac{1+h_{k}\lambda_{1,2}}{\lambda_{1,2}(\lambda_{1,2}+\kappa_{k})}\sum_{i=1}^{n}\omega_{0i}(\textbf{g}_{j}^{\mathrm{T}}\otimes I_{n})(\textbf{w}_{0}\otimes\textbf{e}_{i})\nonumber\\
&=&\frac{\kappa_{k}+\kappa_{k}h_{k}\lambda_{1,2}}{\lambda_{1,2}(\lambda_{1,2}+\kappa_{k})}\sum_{i=1}^{n}\sum_{l=1}^{q}\varpi_{li}(\textbf{g}_{j}^{\mathrm{T}}\textbf{g}_{l}-\textbf{g}_{j}^{\mathrm{T}}G_{k}^{+}G_{k}\textbf{g}_{l})\textbf{e}_{i}-\frac{1+h_{k}\lambda_{1,2}}{\lambda_{1,2}(\lambda_{1,2}+\kappa_{k})}\sum_{i=1}^{n}\omega_{0i}\textbf{e}_{i}.
\end{eqnarray} Hence,
the corresponding eigenvectors for $\lambda_{1,2}$ are given by 
\begin{eqnarray}
\textbf{x}&=&\Big[\frac{1+h_{k}\lambda_{1,2}^{*}}{\lambda_{1,2}^{*}}\sum_{i=1}^{n}\sum_{l=1}^{q}\varpi_{li}((\textbf{g}_{l}-G_{k}^{+}G_{k}\textbf{g}_{l})\otimes \textbf{e}_{i})^{\mathrm{T}}-\frac{1+h_{k}\lambda_{1,2}^{*}}{\kappa_{k}\lambda_{1,2}^{*}}\sum_{i=1}^{n}\omega_{0i}(\textbf{w}_{0}\otimes\textbf{e}_{i})^{\mathrm{T}},\nonumber\\
&&\sum_{i=1}^{n}\sum_{l=1}^{q}\varpi_{li}((\textbf{g}_{l}-G_{k}^{+}G_{k}\textbf{g}_{l})\otimes \textbf{e}_{i})^{\mathrm{T}}-\frac{1}{\kappa_{k}}\sum_{i=1}^{n}\omega_{0i}(\textbf{w}_{0}\otimes\textbf{e}_{i})^{\mathrm{T}},\nonumber\\
&&\frac{\kappa_{k}+\kappa_{k}h_{k}\lambda_{1,2}^{*}}{\lambda_{1,2}^{*}(\lambda_{1,2}^{*}+\kappa_{k})}\sum_{i=1}^{n}\sum_{l=1}^{q}\varpi_{li}(\textbf{g}_{j}^{\mathrm{T}}\textbf{g}_{l}-\textbf{g}_{j}^{\mathrm{T}}G_{k}^{+}G_{k}\textbf{g}_{l})\textbf{e}_{i}^{\mathrm{T}}-\frac{1+h_{k}\lambda_{1,2}^{*}}{\lambda_{1,2}^{*}(\lambda_{1,2}^{*}+\kappa_{k})}\sum_{i=1}^{n}\omega_{0i}\textbf{e}_{i}^{\mathrm{T}}\Big]^{*},
\end{eqnarray} where $\varpi_{li}\in\mathbb{C}$, $\omega_{0i}\in\mathbb{C}$, and not all of them are zero. Therefore,
$\ker\Big(A_{k}^{[j]}+h_{k}A_{{\mathrm{c}}k}-\lambda_{1,2} I_{2nq+n}\Big)$ is given by (\ref{egns2}).

If $\omega_{0i}=0$ for every $i=1,\ldots,n$, then it follows from (\ref{eqn_v}) and (\ref{Az=b}) that
\begin{eqnarray} 
\frac{\kappa_{k}^{2}(1+h_{k}\lambda)}{\lambda(\lambda+\kappa_{k})}W^{[j]}\textbf{v}=\textbf{0}_{nq\times 1},\label{v-1}\\
\Big(\Big(\frac{\mu_{k}}{\lambda}+\mu_{k} h_{k}+\eta_{k}\Big)(L_{k}\otimes I_{n})+\Big(\frac{\kappa_{k}}{\lambda}+\lambda+\kappa_{k} h_{k}\Big)I_{nq}\Big)\textbf{v}=\textbf{0}_{nq\times 1}.\label{v-2}
\end{eqnarray} In this case, since $\frac{\kappa_{k}}{\lambda}+\lambda+\kappa_{k} h_{k}\neq0$ and $\lambda\neq-\kappa_{k}$,  $\det\Big[\Big(\frac{\mu_{k}}{\lambda}+\mu_{k} h_{k}+\eta_{k}\Big)(L_{k}\otimes I_{n})+\Big(\frac{\kappa_{k}}{\lambda}+\lambda+\kappa_{k} h_{k}\Big)I_{nq}\Big]=0$ if and only if $\frac{\mu_{k}}{\lambda}+\mu_{k} h_{k}+\eta_{k}\neq0$ and $\frac{\lambda^{2}+\kappa_{k} h_{k}\lambda+\kappa_{k}}{\eta_{k}\lambda+\mu_{k} h_{k}\lambda+\mu_{k}}\in{\mathrm{spec}}(-L_{k})$. Thus, if $\frac{\kappa_{k}}{\lambda}+\lambda+\kappa_{k} h_{k}\neq0$, $\lambda\neq-\kappa_{k}$, $\frac{\mu_{k}}{\lambda}+\mu_{k} h_{k}+\eta_{k}\neq0$, and $\frac{\lambda^{2}+\kappa_{k} h_{k}\lambda+\kappa_{k}}{\eta_{k}\lambda+\mu_{k} h_{k}\lambda+\mu_{k}}\in{\mathrm{spec}}(-L_{k})$, then $\lambda=\lambda_{4}$, where
\begin{eqnarray}\label{lambda4}
\frac{\lambda_{4}^{2}+\kappa_{k} h_{k}\lambda_{4}+\kappa_{k}}{\eta_{k}\lambda_{4}+\mu_{k} h_{k}\lambda_{4}+\mu_{k}}\in{\mathrm{spec}}(-L_{k}),
\end{eqnarray} are the eigenvalues of $A_{k}^{[j]}+h_{k}A_{{\mathrm{c}}k}$. To find their corresponding eigenvectors,
let $F_{k}=\Big(\frac{\mu_{k}}{\lambda_{4}}+\mu_{k} h_{k}+\eta_{k}\Big)L_{k}+\Big(\frac{\kappa_{k}}{\lambda_{4}}+\lambda_{4}+\kappa_{k} h_{k}\Big)I_{q}$. We first show that (\ref{v-1}) is equivalent to 
\begin{eqnarray}\label{v-3}
\frac{\kappa_{k}^{2}(1+h_{k}\lambda)}{\lambda(\lambda+\kappa_{k})}E_{n\times nq}^{[j]}\textbf{v}=\textbf{0}_{n\times 1}
\end{eqnarray} for every $j=1,\ldots,q$. To see this, let $\textbf{v}=[\textbf{v}_{1}^{*},\ldots,\textbf{v}_{q}^{*}]^{*}$. Then it follows from (\ref{Wj2}) that $W^{[j]}\textbf{v}=[\textbf{v}_{j}^{*},\ldots,\textbf{v}_{j}^{*}]^{*}$. Hence (\ref{v-1}) holds if and only if $\frac{\kappa_{k}^{2}(1+h_{k}\lambda)}{\lambda(\lambda+\kappa_{k})}\textbf{v}_{j}=\textbf{0}_{n\times 1}$. On the other hand, note that $E_{n\times nq}^{[j]}\textbf{v}=\textbf{v}_{j}$. Hence, (\ref{v-1}) is equivalent to (\ref{v-3}).
Then by noting that $E_{n\times nq}^{[j]}=\textbf{g}_{j}^{\mathrm{T}}\otimes I_{n}$ for every $j=1,\ldots,q$, it follows from (\ref{v-2}) and (\ref{v-3}) that
\begin{eqnarray}\label{Fv}
\small\left[\begin{array}{c}
F_{k}\otimes I_{n}\\
\frac{\kappa_{k}^{2}(1+h_{k}\lambda_{4})}{\lambda_{4}(\lambda_{4}+\kappa_{k})}(\textbf{g}_{j}^{\mathrm{T}}\otimes I_{n})
\end{array}\right]\textbf{v}=\Big(\small\left[\begin{array}{c}
F_{k}\\
\frac{\kappa_{k}^{2}(1+h_{k}\lambda_{4})}{\lambda_{4}(\lambda_{4}+\kappa_{k})}\textbf{g}_{j}^{\mathrm{T}}
\end{array}\right]\otimes I_{n}\Big)\textbf{v}=\textbf{0}_{(nq+n)\times 1}.
\end{eqnarray} Next, it follows from $vi$) of Proposition 6.1.7 of \cite[p.~400]{Bernstein:2009} and $viii$) of Proposition 6.1.6 of \cite[p.~399]{Bernstein:2009} that the general solution $\textbf{v}$ to (\ref{Fv}) is given by the form
\begin{eqnarray}
\textbf{v}&=&\Big[I_{nq}-\Big(\small\left[\begin{array}{c}
F_{k}\\
\frac{\kappa_{k}^{2}(1+h_{k}\lambda_{4})}{\lambda_{4}(\lambda_{4}+\kappa_{k})}\textbf{g}_{j}^{\mathrm{T}}
\end{array}\right]\otimes I_{n}\Big)^{+}\Big(\small\left[\begin{array}{c}
F_{k}\\
\frac{\kappa_{k}^{2}(1+h_{k}\lambda_{4})}{\lambda_{4}(\lambda_{4}+\kappa_{k})}\textbf{g}_{j}^{\mathrm{T}}
\end{array}\right]\otimes I_{n}\Big)\Big]\sum_{i=1}^{n}\sum_{l=1}^{q}\varpi_{li}\textbf{g}_{l}\otimes\textbf{e}_{i}\nonumber\\
&=&\Big[I_{nq}-\Big(\small\left[\begin{array}{c}
F_{k}\\
\frac{\kappa_{k}^{2}(1+h_{k}\lambda_{4})}{\lambda_{4}(\lambda_{4}+\kappa_{k})}\textbf{g}_{j}^{\mathrm{T}}
\end{array}\right]^{+}\otimes I_{n}\Big)\Big(\small\left[\begin{array}{c}
F_{k}\\
\frac{\kappa_{k}^{2}(1+h_{k}\lambda_{4})}{\lambda_{4}(\lambda_{4}+\kappa_{k})}\textbf{g}_{j}^{\mathrm{T}}
\end{array}\right]\otimes I_{n}\Big)\Big]\sum_{i=1}^{n}\sum_{l=1}^{q}\varpi_{li}\textbf{g}_{l}\otimes\textbf{e}_{i}\nonumber\\
&=&\Big[I_{q}\otimes I_{n}-\Big(\small\left[\begin{array}{c}
F_{k}\\
\frac{\kappa_{k}^{2}(1+h_{k}\lambda_{4})}{\lambda_{4}(\lambda_{4}+\kappa_{k})}\textbf{g}_{j}^{\mathrm{T}}
\end{array}\right]^{+}\small\left[\begin{array}{c}
F_{k}\\
\frac{\kappa_{k}^{2}(1+h_{k}\lambda_{4})}{\lambda_{4}(\lambda_{4}+\kappa_{k})}\textbf{g}_{j}^{\mathrm{T}}
\end{array}\right]\otimes I_{n}\Big)\Big]\sum_{i=1}^{n}\sum_{l=1}^{q}\varpi_{li}\textbf{g}_{l}\otimes\textbf{e}_{i}\nonumber\\
&=&\Big[\Big(I_{q}-\small\left[\begin{array}{c}
F_{k}\\
\frac{\kappa_{k}^{2}(1+h_{k}\lambda_{4})}{\lambda_{4}(\lambda_{4}+\kappa_{k})}\textbf{g}_{j}^{\mathrm{T}}
\end{array}\right]^{+}\small\left[\begin{array}{c}
F_{k}\\
\frac{\kappa_{k}^{2}(1+h_{k}\lambda_{4})}{\lambda_{4}(\lambda_{4}+\kappa_{k})}\textbf{g}_{j}^{\mathrm{T}}
\end{array}\right]\Big)\otimes I_{n}\Big]\sum_{i=1}^{n}\sum_{l=1}^{q}\varpi_{li}\textbf{g}_{l}\otimes\textbf{e}_{i}\nonumber\\
&=&\sum_{i=1}^{n}\sum_{l=1}^{q}\varpi_{li}\Big(\textbf{g}_{l}-\small\left[\begin{array}{c}
F_{k}\\
\frac{\kappa_{k}^{2}(1+h_{k}\lambda_{4})}{\lambda_{4}(\lambda_{4}+\kappa_{k})}\textbf{g}_{j}^{\mathrm{T}}
\end{array}\right]^{+}\small\left[\begin{array}{c}
F_{k}\\
\frac{\kappa_{k}^{2}(1+h_{k}\lambda_{4})}{\lambda_{4}(\lambda_{4}+\kappa_{k})}\textbf{g}_{j}^{\mathrm{T}}
\end{array}\right]\textbf{g}_{l}\Big)\otimes \textbf{e}_{i},\label{gsolution1a}
\end{eqnarray} where $\varpi_{li}\in\mathbb{C}$ and $j=1,\ldots,q$. Note that by Proposition 6.1.6 of \cite[p.~399]{Bernstein:2009}, $F_{k}^{\mathrm{T}}(F_{k}^{\mathrm{T}})^{+}=F_{k}^{\mathrm{T}}(F_{k}^{+})^{\mathrm{T}}=(F_{k}^{+}F_{k})^{\mathrm{T}}=F_{k}^{+}F_{k}$. It follows from Fact 6.5.17 of \cite[p.~427]{Bernstein:2009} that 
\begin{eqnarray}
\small\left[\begin{array}{c}
F_{k}\\
\frac{\kappa_{k}^{2}(1+h_{k}\lambda_{4})}{\lambda_{4}(\lambda_{4}+\kappa_{k})}\textbf{g}_{j}^{\mathrm{T}}
\end{array}\right]^{+}=\left[\begin{array}{cc}
F_{k}^{+}(I_{q}-\frac{\kappa_{k}^{2}(1+h_{k}\lambda_{4})}{\lambda_{4}(\lambda_{4}+\kappa_{k})}\psi_{k}\textbf{g}_{j}^{\mathrm{T}}) & \psi_{k}
\end{array}\right],
\end{eqnarray} where $\psi_{k}$ is given by (\ref{psik}). Hence, it follows that for every $j,l=1,\ldots,q$,
\begin{eqnarray}
\textbf{g}_{l}-\small\left[\begin{array}{c}
F_{k}\\
\frac{\kappa_{k}^{2}(1+h_{k}\lambda_{4})}{\lambda_{4}(\lambda_{4}+\kappa_{k})}\textbf{g}_{j}^{\mathrm{T}}
\end{array}\right]^{+}\small\left[\begin{array}{c}
F_{k}\\
\frac{\kappa_{k}^{2}(1+h_{k}\lambda_{4})}{\lambda_{4}(\lambda_{4}+\kappa_{k})}\textbf{g}_{j}^{\mathrm{T}}
\end{array}\right]\textbf{g}_{l}&=&\textbf{g}_{l}-\left[\begin{array}{cc}
F_{k}^{+}(I_{q}-\frac{\kappa_{k}^{2}(1+h_{k}\lambda_{4})}{\lambda_{4}(\lambda_{4}+\kappa_{k})}\psi_{k}\textbf{g}_{j}^{\mathrm{T}}) & \psi_{k}
\end{array}\right]\small\left[\begin{array}{c}
F_{k}\\
\frac{\kappa_{k}^{2}(1+h_{k}\lambda_{4})}{\lambda_{4}(\lambda_{4}+\kappa_{k})}\textbf{g}_{j}^{\mathrm{T}}
\end{array}\right]\textbf{g}_{l}\nonumber\\
&=&\textbf{g}_{l}-\small\left[\begin{array}{cc}
F_{k}^{+}(I_{q}-\frac{\kappa_{k}^{2}(1+h_{k}\lambda_{4})}{\lambda_{4}(\lambda_{4}+\kappa_{k})}\psi_{k}\textbf{g}_{j}^{\mathrm{T}}) & \psi_{k}
\end{array}\right]\small\left[\begin{array}{c}
F_{k}\textbf{g}_{l}\\
\frac{\kappa_{k}^{2}(1+h_{k}\lambda_{4})}{\lambda_{4}(\lambda_{4}+\kappa_{k})}\textbf{g}_{j}^{\mathrm{T}}\textbf{g}_{l}
\end{array}\right]\nonumber\\
&=&\textbf{g}_{l}-F_{k}^{+}\Big(I_{q}-\frac{\kappa_{k}^{2}(1+h_{k}\lambda_{4})}{\lambda_{4}(\lambda_{4}+\kappa_{k})}\psi_{k}\textbf{g}_{j}^{\mathrm{T}}\Big)F_{k}\textbf{g}_{l}\nonumber\\
&&-\frac{\kappa_{k}^{2}(1+h_{k}\lambda_{4})}{\lambda_{4}(\lambda_{4}+\kappa_{k})}(\textbf{g}_{j}^{\mathrm{T}}\textbf{g}_{l})\psi_{k}\nonumber\\
&=&\textbf{g}_{l}-F_{k}^{+}F_{k}\textbf{g}_{l}+\frac{\kappa_{k}^{2}(1+h_{k}\lambda_{4})}{\lambda_{4}(\lambda_{4}+\kappa_{k})}(\textbf{g}_{j}^{\mathrm{T}}F_{k}\textbf{g}_{l})F_{k}^{+}\psi_{k}\nonumber\\
&&-\frac{\kappa_{k}^{2}(1+h_{k}\lambda_{4})}{\lambda_{4}(\lambda_{4}+\kappa_{k})}(\textbf{g}_{j}^{\mathrm{T}}\textbf{g}_{l})\psi_{k}.
\end{eqnarray} Thus, (\ref{gsolution1a}) becomes 
\begin{eqnarray}\label{gsoa}
\textbf{v}=\sum_{i=1}^{n}\sum_{l=1}^{q}\varpi_{li}\Big(\textbf{g}_{l}-F_{k}^{+}F_{k}\textbf{g}_{l}+\frac{\kappa_{k}^{2}(1+h_{k}\lambda_{4})}{\lambda_{4}(\lambda_{4}+\kappa_{k})}(\textbf{g}_{j}^{\mathrm{T}}F_{k}\textbf{g}_{l})F_{k}^{+}\psi_{k}-\frac{\kappa_{k}^{2}(1+h_{k}\lambda_{4})}{\lambda_{4}(\lambda_{4}+\kappa_{k})}(\textbf{g}_{j}^{\mathrm{T}}\textbf{g}_{l})\psi_{k}\Big)\otimes \textbf{e}_{i}.
\end{eqnarray} Hence, $\textbf{x}_{1}=\frac{1+h_{k}\lambda_{4}}{\lambda_{4}}\textbf{v}$, $\textbf{x}_{2}=\textbf{v}\neq\textbf{0}_{nq\times 1}$ given by (\ref{gsoa}), and 
\begin{eqnarray}
\textbf{x}_{3}&=&\frac{\kappa_{k}+\kappa_{k}h_{k}\lambda_{4}}{\lambda_{4}(\lambda_{4}+\kappa_{k})}E_{n\times nq}^{[j]}\textbf{v}\nonumber\\
&=&\frac{\kappa_{k}+\kappa_{k}h_{k}\lambda_{4}}{\lambda_{4}(\lambda_{4}+\kappa_{k})}(\textbf{g}_{j}^{\mathrm{T}}\otimes I_{n})\textbf{v}\nonumber\\
&=&\frac{\kappa_{k}+\kappa_{k}h_{k}\lambda_{4}}{\lambda_{4}(\lambda_{4}+\kappa_{k})}\sum_{i=1}^{n}\sum_{l=1}^{q}\varpi_{li}(\textbf{g}_{j}^{\mathrm{T}}\otimes I_{n})\Big(\Big(\textbf{g}_{l}-F_{k}^{+}F_{k}\textbf{g}_{l}+\frac{\kappa_{k}^{2}(1+h_{k}\lambda_{4})}{\lambda_{4}(\lambda_{4}+\kappa_{k})}(\textbf{g}_{j}^{\mathrm{T}}F_{k}\textbf{g}_{l})F_{k}^{+}\psi_{k}\nonumber\\
&&-\frac{\kappa_{k}^{2}(1+h_{k}\lambda_{4})}{\lambda_{4}(\lambda_{4}+\kappa_{k})}(\textbf{g}_{j}^{\mathrm{T}}\textbf{g}_{l})\psi_{k}\Big)\otimes \textbf{e}_{i}\Big)\nonumber\\
&=&\frac{\kappa_{k}+\kappa_{k}h_{k}\lambda_{4}}{\lambda_{4}(\lambda_{4}+\kappa_{k})}\sum_{i=1}^{n}\sum_{l=1}^{q}\varpi_{li}\Big(\textbf{g}_{j}^{\mathrm{T}}\textbf{g}_{l}-\textbf{g}_{j}^{\mathrm{T}}F_{k}^{+}F_{k}\textbf{g}_{l}+\frac{\kappa_{k}^{2}(1+h_{k}\lambda_{4})}{\lambda_{4}(\lambda_{4}+\kappa_{k})}(\textbf{g}_{j}^{\mathrm{T}}F_{k}\textbf{g}_{l})\textbf{g}_{j}^{\mathrm{T}}F_{k}^{+}\psi_{k}\nonumber\\
&&-\frac{\kappa_{k}^{2}(1+h_{k}\lambda_{4})}{\lambda_{4}(\lambda_{4}+\kappa_{k})}(\textbf{g}_{j}^{\mathrm{T}}\textbf{g}_{l})\textbf{g}_{j}^{\mathrm{T}}\psi_{k}\Big)\otimes \textbf{e}_{i},\label{gsox3}
\end{eqnarray}
where not all of $\omega_{\ell i}$ and $\varpi_{li}$ are zero. The corresponding eigenvectors for $\lambda_{4}$ are given by
\begin{eqnarray}\label{eigv4}
&&\textbf{x}=\nonumber\\
&&\Big[\frac{1+h_{k}\lambda_{4}^{*}}{\lambda_{4}^{*}}\sum_{i=1}^{n}\sum_{l=1}^{q}\varpi_{li}\Big(\textbf{g}_{l}-F_{k}^{+}F_{k}\textbf{g}_{l}+\frac{\kappa_{k}^{2}(1+h_{k}\lambda_{4})}{\lambda_{4}(\lambda_{4}+\kappa_{k})}(\textbf{g}_{j}^{\mathrm{T}}F_{k}\textbf{g}_{l})F_{k}^{+}\psi_{k}-\frac{\kappa_{k}^{2}(1+h_{k}\lambda_{4})}{\lambda_{4}(\lambda_{4}+\kappa_{k})}(\textbf{g}_{j}^{\mathrm{T}}\textbf{g}_{l})\psi_{k}\Big)^{*}\otimes \textbf{e}_{i}^{\mathrm{T}},\nonumber\\
&&\sum_{i=1}^{n}\sum_{l=1}^{q}\varpi_{li}\Big(\textbf{g}_{l}-F_{k}^{+}F_{k}\textbf{g}_{l}+\frac{\kappa_{k}^{2}(1+h_{k}\lambda_{4})}{\lambda_{4}(\lambda_{4}+\kappa_{k})}(\textbf{g}_{j}^{\mathrm{T}}F_{k}\textbf{g}_{l})F_{k}^{+}\psi_{k}-\frac{\kappa_{k}^{2}(1+h_{k}\lambda_{4})}{\lambda_{4}(\lambda_{4}+\kappa_{k})}(\textbf{g}_{j}^{\mathrm{T}}\textbf{g}_{l})\psi_{k}\Big)^{*}\otimes \textbf{e}_{i}^{\mathrm{T}},\nonumber\\
&&\frac{\kappa_{k}+\kappa_{k}h_{k}\lambda_{4}^{*}}{\lambda_{4}^{*}(\lambda_{4}^{*}+\kappa_{k})}\sum_{i=1}^{n}\sum_{l=1}^{q}\varpi_{li}\Big(\textbf{g}_{j}^{\mathrm{T}}\textbf{g}_{l}-\textbf{g}_{j}^{\mathrm{T}}F_{k}^{+}F_{k}\textbf{g}_{l}+\frac{\kappa_{k}^{2}(1+h_{k}\lambda_{4})}{\lambda_{4}(\lambda_{4}+\kappa_{k})}(\textbf{g}_{j}^{\mathrm{T}}F_{k}\textbf{g}_{l})\textbf{g}_{j}^{\mathrm{T}}F_{k}^{+}\psi_{k}\nonumber\\
&&-\frac{\kappa_{k}^{2}(1+h_{k}\lambda_{4})}{\lambda_{4}(\lambda_{4}+\kappa_{k})}(\textbf{g}_{j}^{\mathrm{T}}\textbf{g}_{l})\textbf{g}_{j}^{\mathrm{T}}\psi_{k}\Big)^{*}\otimes \textbf{e}_{i}^{\mathrm{T}}\Big]^{*},
\end{eqnarray} where $\varpi_{li}\in\mathbb{C}$ and not all of them are zero. Therefore,
$\ker\Big(A_{k}^{[j]}+h_{k}A_{{\mathrm{c}}k}-\lambda_{4} I_{2nq+n}\Big)$ is given by (\ref{egns3}).

If $\frac{\kappa_{k}}{\lambda}+\lambda+\kappa_{k} h_{k}=0$, then $\frac{\kappa_{k}^{2}(1+h_{k}\lambda)}{\lambda(\lambda+\kappa_{k})}=-\frac{\kappa_{k}\lambda}{\lambda+\kappa_{k}}\neq0$ since $\lambda\neq0$ and $\kappa_{k}\neq0$. In this case, it follows from (\ref{eqn_v}) and (\ref{Az=b}) that
\begin{eqnarray} &&\frac{\kappa_{k}^{2}(1+h_{k}\lambda)}{\lambda(\lambda+\kappa_{k})}W^{[j]}\textbf{v}=\sum_{i=1}^{n}\sum_{l=0}^{q-1-{\mathrm{rank}}(L_{k})}\omega_{li}\textbf{w}_{l}\otimes\textbf{e}_{i},\label{id_1}\\
&&\Big(\frac{\mu_{k}}{\lambda}+\mu_{k} h_{k}+\eta_{k}\Big)(L_{k}\otimes I_{n})\textbf{v}=\sum_{i=1}^{n}\sum_{l=0}^{q-1-{\mathrm{rank}}(L_{k})}\omega_{li}\textbf{w}_{l}\otimes\textbf{e}_{i}.\label{id_2}
\end{eqnarray} Since $W^{[j]}$ is idempotent by $i$) of Lemma~\ref{lemma_EW}, it follows from (\ref{id_1}) and $ii$) of Lemma~\ref{lemma_EW} that 
 \begin{eqnarray}
 \sum_{i=1}^{n}\sum_{l=0}^{q-1-{\mathrm{rank}}(L_{k})}\omega_{li}\textbf{w}_{l}\otimes\textbf{e}_{i}=\sum_{i=1}^{n}\sum_{l=0}^{q-1-{\mathrm{rank}}(L_{k})}\omega_{li}w_{lj}\textbf{w}_{0}\otimes\textbf{e}_{i},
 \end{eqnarray}
 and hence, 
\begin{eqnarray}
\sum_{i=1}^{n}\Big(\omega_{0i}-\sum_{l=0}^{q-1-{\mathrm{rank}}(L_{k})}\omega_{li}w_{lj}\Big)\textbf{w}_{0}\otimes\textbf{e}_{i}+\sum_{i=1}^{n}\sum_{l=1}^{q-1-{\mathrm{rank}}(L_{k})}\omega_{li}\textbf{w}_{l}\otimes\textbf{e}_{i}=\textbf{0}_{nq\times 1},
\end{eqnarray} which implies that $\omega_{0i}-\sum_{l=0}^{q-1-{\mathrm{rank}}(L_{k})}\omega_{li}w_{lj}=0$ and $\omega_{\ell i}=0$ for every $i=1,\ldots,n$, $j=1,\ldots,q$, and $\ell=1,\ldots,q-1-{\mathrm{rank}}(L_{k})$. Consequently, (\ref{id_1}) and (\ref{id_2}) can be simplified as
\begin{eqnarray} &&\frac{\kappa_{k}^{2}(1+h_{k}\lambda)}{\lambda(\lambda+\kappa_{k})}W^{[j]}\textbf{v}=\sum_{i=1}^{n}\omega_{0i}\textbf{w}_{0}\otimes\textbf{e}_{i},\label{new_1}\\
&&\Big(\frac{\mu_{k}}{\lambda}+\mu_{k} h_{k}+\eta_{k}\Big)(L_{k}\otimes I_{n})\textbf{v}=\sum_{i=1}^{n}\omega_{0i}\textbf{w}_{0}\otimes\textbf{e}_{i}.\label{new_2}
\end{eqnarray} It follows from $ii$) of Lemma~\ref{lemma_EW} that (\ref{new_1}) has a specific solution 
\begin{eqnarray}\label{vss}
\textbf{v}=\Big(\frac{\kappa_{k}^{2}(1+h_{k}\lambda)}{\lambda(\lambda+\kappa_{k})}\Big)^{-1}\sum_{i=1}^{n}\omega_{0i}\textbf{w}_{0}\otimes\textbf{e}_{i}.
\end{eqnarray} Substituting (\ref{vss}) into (\ref{new_2}) yields $\sum_{i=1}^{n}\omega_{0i}\textbf{w}_{0}\otimes\textbf{e}_{i}=\textbf{0}_{nq\times 1}$, which implies that $\omega_{0i}=0$ for every $i=1,\ldots,n$. Hence, (\ref{new_1}) and (\ref{new_2}) can be further simplified as 
\begin{eqnarray} 
&&\frac{\kappa_{k}^{2}(1+h_{k}\lambda)}{\lambda(\lambda+\kappa_{k})}W^{[j]}\textbf{v}=\textbf{0}_{nq\times 1},\label{gs0_1}\\
&&\Big(\frac{\mu_{k}}{\lambda}+\mu_{k} h_{k}+\eta_{k}\Big)(L_{k}\otimes I_{n})\textbf{v}=\textbf{0}_{nq\times 1}.\label{gs0_2}
\end{eqnarray} If $\frac{\mu_{k}}{\lambda}+\mu_{k} h_{k}+\eta_{k}\neq0$, note that for $\frac{\kappa_{k}}{\lambda}+\lambda+\kappa_{k} h_{k}=0$, $\det\Big[\Big(\frac{\mu_{k}}{\lambda}+\mu_{k} h_{k}+\eta_{k}\Big)(L_{k}\otimes I_{n})+\Big(\frac{\kappa_{k}}{\lambda}+\lambda+\kappa_{k} h_{k}\Big)I_{nq}\Big]=\det\Big[\Big(\frac{\mu_{k}}{\lambda}+\mu_{k} h_{k}+\eta_{k}\Big)(L_{k}\otimes I_{n})\Big]=0$. Hence, the general solution $\textbf{v}$ to (\ref{gs0_1}) and (\ref{gs0_2}) is given by the form of (\ref{gsoa}) in which $\lambda_{4}$ is replaced by $\lambda_{5,6}$ satisfying $\frac{\kappa_{k}}{\lambda_{5,6}}+\lambda_{5,6}+\kappa_{k} h_{k}=0$. Thus, this case is similar to the previous case where (\ref{lambda4}) still holds for $\lambda_{4}$ being replaced by $\lambda_{5,6}$, where
\begin{eqnarray}\label{lambda5}
\lambda_{5,6}=-\frac{\kappa_{k}h_{k}}{2}\pm\frac{1}{2}\sqrt{\kappa_{k}^{2}h_{k}^{2}-4\kappa_{k}}.
\end{eqnarray} Thus, $\lambda=\lambda_{5,6}$ are indeed the eigenvalues of $A_{k}^{[j]}+h_{k}A_{{\mathrm{c}}k}$ and the corresponding eigenvectors are given by the form (\ref{eigv4}) with $\lambda_{4}$ being replaced by $\lambda_{5,6}$. 

Otherwise, if $\frac{\mu_{k}}{\lambda}+\mu_{k} h_{k}+\eta_{k}=0$ and $\frac{\kappa_{k}}{\lambda}+\lambda+\kappa_{k} h_{k}=0$, then $\mu_{k}(\frac{1}{\lambda}+h_{k})=-\eta_{k}$ and $\kappa_{k}(\frac{1}{\lambda}+h_{k})=-\lambda$. Again, since $\lambda\neq0$, it follows from $\frac{\kappa_{k}}{\lambda}+\lambda+\kappa_{k} h_{k}=0$ that $\kappa_{k}\neq0$. If $\mu_{k}=0$, then it follows from $\mu_{k}(\frac{1}{\lambda}+h_{k})=-\eta_{k}$ that $\eta_{k}=0$. In this case, $\lambda=\lambda_{5,6}$ are the eigenvalues of $A_{k}^{[j]}+h_{k}A_{{\mathrm{c}}k}$. Furthermore, (\ref{gs0_2}) becomes trivial and (\ref{gs0_1}) is equivalent to $E_{n\times nq}^{[j]}\textbf{v}=\textbf{0}_{n\times 1}$, that is, $(\textbf{g}_{j}^{\mathrm{T}}\otimes I_{n})\textbf{v}=\textbf{0}_{n\times 1}$. It follows from $vi$) of Proposition 6.1.7 of \cite[p.~400]{Bernstein:2009} and $viii$) of Proposition 6.1.6 of \cite[p.~399]{Bernstein:2009} that the general solution $\textbf{v}$ to $(\textbf{g}_{j}^{\mathrm{T}}\otimes I_{n})\textbf{v}=\textbf{0}_{n\times 1}$ is given by the form
\begin{eqnarray}
\textbf{v}&=&\Big[I_{nq}-(\textbf{g}_{j}^{\mathrm{T}}\otimes I_{n})^{+}(\textbf{g}_{j}^{\mathrm{T}}\otimes I_{n})\Big]\sum_{i=1}^{n}\sum_{l=1}^{q}\varpi_{li}\textbf{g}_{l}\otimes\textbf{e}_{i}\nonumber\\
&=&\Big[I_{nq}-((\textbf{g}_{j}^{\mathrm{T}})^{+}\otimes I_{n})(\textbf{g}_{j}^{\mathrm{T}}\otimes I_{n})\Big]\sum_{i=1}^{n}\sum_{l=1}^{q}\varpi_{li}\textbf{g}_{l}\otimes\textbf{e}_{i}\nonumber\\
&=&\Big[I_{q}\otimes I_{n}-(((\textbf{g}_{j}^{\mathrm{T}})^{+}\textbf{g}_{j}^{\mathrm{T}})\otimes I_{n})\Big]\sum_{i=1}^{n}\sum_{l=1}^{q}\varpi_{li}\textbf{g}_{l}\otimes\textbf{e}_{i}\nonumber\\
&=&\Big[(I_{q}-((\textbf{g}_{j}^{\mathrm{T}})^{+}\textbf{g}_{j}^{\mathrm{T}}))\otimes I_{n}\Big]\sum_{i=1}^{n}\sum_{l=1}^{q}\varpi_{li}\textbf{g}_{l}\otimes\textbf{e}_{i}\nonumber\\
&=&\sum_{i=1}^{n}\sum_{l=1}^{q}\varpi_{li}(\textbf{g}_{l}-((\textbf{g}_{j}^{\mathrm{T}})^{+}\textbf{g}_{j}^{\mathrm{T}})\textbf{g}_{l})\otimes \textbf{e}_{i},
\end{eqnarray} where $\varpi_{li}\in\mathbb{C}$ and $j=1,\ldots,q$. Note that it follows from Fact 6.3.2 of \cite[p.~404]{Bernstein:2009} that $\textbf{g}_{j}^{+}=\textbf{g}_{j}^{\mathrm{T}}$, and hence, $(\textbf{g}_{j}^{\mathrm{T}})^{+}=\textbf{g}_{j}$ for every $j=1,\ldots,q$. Then we have
\begin{eqnarray}\label{newv}
\textbf{v}&=&\sum_{i=1}^{n}\sum_{l=1}^{q}\varpi_{li}(\textbf{g}_{l}-(\textbf{g}_{j}\textbf{g}_{j}^{\mathrm{T}})\textbf{g}_{l})\otimes \textbf{e}_{i}\nonumber\\
&=&\sum_{i=1}^{n}\sum_{l=1}^{q}\varpi_{li}(\textbf{g}_{l}-(\textbf{g}_{j}^{\mathrm{T}}\textbf{g}_{l})\textbf{g}_{j})\otimes \textbf{e}_{i}.
\end{eqnarray} Hence, $\textbf{x}_{1}=\frac{1+h_{k}\lambda_{5,6}}{\lambda_{5,6}}\textbf{v}$, $\textbf{x}_{2}=\textbf{v}\neq\textbf{0}_{nq\times 1}$ where $\textbf{v}$ is given by (\ref{newv}), and $\textbf{x}_{3}=\textbf{0}_{n\times 1}$. The corresponding eigenvectors for $\lambda_{5,6}$ in this case are given by 
\begin{eqnarray}
\textbf{x}=\Big[\frac{1+h_{k}\lambda_{5,6}^{*}}{\lambda_{5,6}^{*}}\sum_{i=1}^{n}\sum_{l=1}^{q}\varpi_{li}(\textbf{g}_{l}-(\textbf{g}_{j}^{\mathrm{T}}\textbf{g}_{l})\textbf{g}_{j})^{\mathrm{T}}\otimes \textbf{e}_{i}^{\mathrm{T}},\sum_{i=1}^{n}\sum_{l=1}^{q}\varpi_{li}(\textbf{g}_{l}-(\textbf{g}_{j}^{\mathrm{T}}\textbf{g}_{l})\textbf{g}_{j})^{\mathrm{T}}\otimes \textbf{e}_{i}^{\mathrm{T}},\textbf{0}_{1\times n}\Big]^{*},
\end{eqnarray} where $\varpi_{li}\in\mathbb{C}$ and not all of them are zero. Consequently, in this case
$\ker\Big(A_{k}^{[j]}+h_{k}A_{{\mathrm{c}}k}-\lambda_{5,6} I_{2nq+n}\Big)$ is given by (\ref{egns4}).

Finally, if $\mu_{k}\neq0$, then it follows from $\mu_{k}(\frac{1}{\lambda}+h_{k})=-\eta_{k}$ that $\frac{1}{\lambda}+h_{k}=-\frac{\eta_{k}}{\mu_{k}}$. Together with $\kappa_{k}(\frac{1}{\lambda}+h_{k})=-\lambda$, we have $\lambda=\frac{\kappa_{k}\eta_{k}}{\mu_{k}}$. Since $\lambda\neq0$, it follows that $\eta_{k}\neq0$. Substituting this $\lambda$ into $\frac{1}{\lambda}+h_{k}=-\frac{\eta_{k}}{\mu_{k}}$ yields $h_{k}=-\frac{\eta_{k}}{\mu_{k}}-\frac{\mu_{k}}{\kappa_{k}\eta_{k}}<0$, which is a contradiction since $h_{k}>0$. Hence, this case is impossible.

\textit{Case 2.} If $\lambda=-\kappa_{k}$, then $\kappa_{k}\neq0$ and (\ref{detcon}) becomes
\begin{eqnarray}\label{detcase2}
\det\small\left[\begin{array}{cc}
\Big(\frac{\mu_{k}}{\kappa_{k}}(\kappa_{k} h_{k}-1)+\eta_{k}\Big)(L_{k}\otimes I_{n})+(\kappa_{k} h_{k}-1-\kappa_{k})I_{nq} & -\kappa_{k}(\textbf{1}_{q\times 1}\otimes I_{n})\\
(\kappa_{k} h_{k}-1)E_{n\times nq}^{[j]} & \textbf{0}_{n\times n}
\end{array}\right]=0.
\end{eqnarray} If $\kappa_{k}h_{k}=1$, then clearly (\ref{detcase2}) holds. In this case,
\begin{eqnarray*}
&&\det\Big[\Big(\frac{\mu_{k}}{\lambda}+\mu_{k} h_{k}+\eta_{k}\Big)(L_{k}\otimes I_{n})+\Big(\frac{\kappa_{k}}{\lambda}+\lambda+\kappa_{k} h_{k}\Big)I_{nq}\Big]\nonumber\\
&&=\det\Big[\Big(-\frac{\mu_{k}}{\kappa_{k}}+\mu_{k} h_{k}+\eta_{k}\Big)(L_{k}\otimes I_{n})-\kappa_{k}I_{nq}\Big]\nonumber\\
&&=\kappa_{k}^{nq}\det\Big[\Big(-\frac{\mu_{k}}{\kappa_{k}^{2}}+\frac{\mu_{k} h_{k}}{\kappa_{k}}+\frac{\eta_{k}}{\kappa_{k}}\Big)(L_{k}\otimes I_{n})-I_{nq}\Big]\nonumber\\
&&=\kappa_{k}^{nq}\det\Big[\frac{\eta_{k}}{\kappa_{k}}(L_{k}\otimes I_{n})-I_{nq}\Big].
\end{eqnarray*} Hence, $\det\Big[\Big(\frac{\mu_{k}}{\lambda}+\mu_{k} h_{k}+\eta_{k}\Big)(L_{k}\otimes I_{n})+\Big(\frac{\kappa_{k}}{\lambda}+\lambda+\kappa_{k} h_{k}\Big)I_{nq}\Big]=0$ if and only if $1\in{\mathrm{spec}}(\frac{\eta_{k}}{\kappa_{k}}L_{k})$. Thus, if 
$1\in{\mathrm{spec}}(\frac{\eta_{k}}{\kappa_{k}}L_{k})$ and $\kappa_{k}h_{k}=1$, then $\lambda=-\kappa_{k}$ is indeed an eigenvalue of $A_{k}^{[j]}+h_{k}A_{{\mathrm{c}}k}$. Clearly when $\kappa_{k}h_{k}=1$ and $\lambda=-\kappa_{k}$, $\textbf{x}_{1}=\frac{1+h_{k}\lambda}{\lambda}\textbf{x}_{2}=\textbf{0}_{nq\times 1}$, (\ref{hx3}) becomes trivial, and (\ref{hx2}) becomes
\begin{eqnarray}\label{x2x3}
(\eta_{k}(L_{k}\otimes I_{n})-\kappa_{k}I_{nq})\textbf{x}_{2}-\kappa_{k}(\textbf{1}_{q\times 1}\otimes I_{n})\textbf{x}_{3}=\textbf{0}_{nq\times1}.
\end{eqnarray} Pre-multiplying $E_{n\times nq}^{[j]}$ on both sides of (\ref{x2x3}) yields
\begin{eqnarray}\label{x3x2}
\textbf{x}_{3}=\Big[\frac{\eta_{k}}{\kappa_{k}}(L_{k}\otimes I_{n})-I_{nq}\Big]\textbf{x}_{2}.
\end{eqnarray} Note that $\textbf{x}_{2}$ can be chosen arbitrarily in $\mathbb{C}^{nq}$ other than $\textbf{0}_{nq\times 1}$. Then $\textbf{x}_{2}$ can be represented as $\textbf{x}_{2}=\sum_{i=1}^{n}\sum_{l=1}^{q}\\\alpha_{li}(\textbf{g}_{l}\otimes\textbf{e}_{i})$, where $\alpha_{li}\in\mathbb{C}$, not all of $\alpha_{li}$ are zero, and $[\textbf{g}_{1},\ldots,\textbf{g}_{q}]=I_{q}$. Then it follows from (\ref{x3x2}) that $\textbf{x}_{3}=\sum_{i=1}^{n}\sum_{l=1}^{q}\frac{\eta_{k}}{\kappa_{k}}\alpha_{li}(L_{k}\otimes I_{n})(\textbf{g}_{l}\otimes\textbf{e}_{i})-\sum_{i=1}^{n}\sum_{l=1}^{q}\alpha_{li}(\textbf{g}_{l}\otimes\textbf{e}_{i})=\sum_{i=1}^{n}\sum_{l=1}^{q}\frac{\eta_{k}}{\kappa_{k}}\alpha_{li}(L_{k}\textbf{g}_{l}\otimes\textbf{e}_{i})-\sum_{i=1}^{n}\sum_{l=1}^{q}\alpha_{li}(\textbf{g}_{l}\otimes\textbf{e}_{i})$, where $\alpha_{li}\in\mathbb{C}$ and not all of $\alpha_{il}$ are zero. Clearly such $\textbf{x}_{i}$, $i=1,2,3$, satisfy (\ref{Aeig_1})--(\ref{Aeig_3}). Thus, the corresponding eigenvectors for the eigenvalue $\lambda=\lambda_{3}$ are given by
\begin{eqnarray}
\textbf{x}=\Big[\textbf{0}_{1\times nq},\sum_{i=1}^{n}\sum_{l=1}^{q}\alpha_{li}(\textbf{g}_{l}\otimes\textbf{e}_{i})^{\mathrm{T}},\sum_{i=1}^{n}\sum_{l=1}^{q}\frac{\eta_{k}}{\kappa_{k}}\alpha_{li}(L_{k}\textbf{g}_{l}\otimes\textbf{e}_{i})^{\mathrm{T}}-\sum_{i=1}^{n}\sum_{l=1}^{q}\alpha_{li}(\textbf{g}_{l}\otimes\textbf{e}_{i})^{\mathrm{T}}\Big]^{*},
\end{eqnarray} where $\alpha_{li}\in\mathbb{C}$, not all of $\alpha_{il}$ are zero, and 
\begin{eqnarray}\label{lambda3}
\lambda_{3}=-\kappa_{k}.
\end{eqnarray} Therefore,
$\ker\Big(A_{k}^{[j]}+h_{k}A_{{\mathrm{c}}k}-\lambda_{3} I_{2nq+n}\Big)$ is given by (\ref{egns5}).

Now we consider the case where $\kappa_{k}h_{k}\neq1$. Then in this case (\ref{detcase2}) holds if and only if the equation
\begin{eqnarray}
\small\left[\begin{array}{cc}
\Big(\frac{\mu_{k}}{\kappa_{k}}(\kappa_{k} h_{k}-1)+\eta_{k}\Big)(L_{k}\otimes I_{n})+(\kappa_{k} h_{k}-1-\kappa_{k})I_{nq} & -\kappa_{k}(\textbf{1}_{q\times 1}\otimes I_{n})\\
(\kappa_{k} h_{k}-1)E_{n\times nq}^{[j]} & \textbf{0}_{n\times n}
\end{array}\right]\textbf{u}=\textbf{0}_{(nq+n)\times 1}\label{deteqn}
\end{eqnarray} has a nontrivial solution $\textbf{u}\in\mathbb{C}^{nq+n}$. Let $\textbf{u}=[\textbf{u}_{1}^{*},\ldots,\textbf{u}_{q}^{*},\textbf{u}_{0}^{*}]^{*}$, where $\textbf{u}_{i}\in\mathbb{C}^{n}$, $i=0,1,\ldots,q$. Then it follows from (\ref{deteqn}) that 
\begin{eqnarray}
\Big(\frac{\mu_{k}}{\kappa_{k}}(\kappa_{k} h_{k}-1)+\eta_{k}\Big)(L_{k}\otimes I_{n})[\textbf{u}_{1}^{*},\ldots,\textbf{u}_{q}^{*}]^{*}+(\kappa_{k} h_{k}-1-\kappa_{k})[\textbf{u}_{1}^{*},\ldots,\textbf{u}_{q}^{*}]^{*}\nonumber\\
-\kappa_{k}(\textbf{1}_{q\times 1}\otimes I_{n})\textbf{u}_{0}=\textbf{0}_{nq\times 1},\label{Enqu1}\\
(\kappa_{k} h_{k}-1)E_{n\times nq}^{[j]}[\textbf{u}_{1}^{*},\ldots,\textbf{u}_{q}^{*}]^{*}=\textbf{0}_{n\times 1}.\label{Enqu}
\end{eqnarray}   

If $\frac{\mu_{k}}{\kappa_{k}}(\kappa_{k} h_{k}-1)+\eta_{k}=0$, in this case, since $\lambda=-\kappa_{k}$, then it follows that
\begin{eqnarray*}
\det\Big[\Big(\frac{\mu_{k}}{\lambda}+\mu_{k} h_{k}+\eta_{k}\Big)(L_{k}\otimes I_{n})+\Big(\frac{\kappa_{k}}{\lambda}+\lambda+\kappa_{k} h_{k}\Big)I_{nq}\Big]&=&\det\Big[(\kappa_{k} h_{k}-1-\kappa_{k})I_{nq}\Big]\nonumber\\
&=&(\kappa_{k} h_{k}-1-\kappa_{k})^{nq}.
\end{eqnarray*} Hence, $\det\Big[\Big(\frac{\mu_{k}}{\lambda}+\mu_{k} h_{k}+\eta_{k}\Big)(L_{k}\otimes I_{n})+\Big(\frac{\kappa_{k}}{\lambda}+\lambda+\kappa_{k} h_{k}\Big)I_{nq}\Big]=0$ if and only if $\kappa_{k} h_{k}-1-\kappa_{k}=0$. If $\kappa_{k} h_{k}-1-\kappa_{k}=0$, eliminating $h_{k}$ in $\frac{\mu_{k}}{\kappa_{k}}(\kappa_{k} h_{k}-1)+\eta_{k}=0$ by using $\kappa_{k} h_{k}-1-\kappa_{k}=0$ yields $\mu_{k}+\eta_{k}=0$, and hence, $\mu_{k}=\eta_{k}=0$ since $\mu_{k},\eta_{k}\geq0$. Furthermore, $h_{k}\kappa_{k}=1+\kappa_{k}\neq 1$ due to $\kappa_{k}\neq0$. Next, since $\frac{\mu_{k}}{\kappa_{k}}(\kappa_{k} h_{k}-1)+\eta_{k}=0$ and $\kappa_{k} h_{k}-1-\kappa_{k}=0$, it follows from (\ref{Enqu1}) that $\textbf{u}_{0}=\textbf{0}_{n\times 1}$. Thus in this case, (\ref{Enqu}) becomes $E_{n\times nq}^{[j]}[\textbf{u}_{1}^{*},\ldots,\textbf{u}_{q}^{*}]^{*}=\textbf{0}_{n\times 1}$, that is, $(\textbf{g}_{j}^{\mathrm{T}}\otimes I_{n})[\textbf{u}_{1}^{*},\ldots,\textbf{u}_{q}^{*}]^{*}=\textbf{0}_{n\times 1}$. Now it follows from (\ref{newv}) that $[\textbf{u}_{1}^{*},\ldots,\textbf{u}_{q}^{*}]^{*}=\sum_{i=1}^{n}\sum_{l=1}^{q}\alpha_{li}(\textbf{g}_{l}-(\textbf{g}_{j}^{\mathrm{T}}\textbf{g}_{l})\textbf{g}_{j})\otimes \textbf{e}_{i}$, where $\alpha_{li}\in\mathbb{C}$ and not all of them are zero. Clearly $\textbf{x}_{1}=\textbf{0}_{nq\times 1}$, $\textbf{x}_{2}=\sum_{i=1}^{n}\sum_{l=1}^{q}\alpha_{li}(\textbf{g}_{l}-(\textbf{g}_{j}^{\mathrm{T}}\textbf{g}_{l})\textbf{g}_{j})\otimes \textbf{e}_{i}$, and $\textbf{x}_{3}=\textbf{0}_{n\times 1}$ satisfy (\ref{Aeig_1})--(\ref{Aeig_3}). Thus, if 
$\frac{\mu_{k}}{\kappa_{k}}(\kappa_{k} h_{k}-1)+\eta_{k}=0$ and $h_{k}=1+\frac{1}{\kappa_{k}}$, then $\lambda=-\kappa_{k}$ is indeed an eigenvalue of $A_{k}^{[j]}+h_{k}A_{{\mathrm{c}}k}$ and the corresponding eigenvectors for the eigenvalue $\lambda_{3}$ of the form (\ref{lambda3}) are given by
\begin{eqnarray}
\textbf{x}=\Big[\textbf{0}_{1\times nq},\sum_{i=1}^{n}\sum_{l=1}^{q}\alpha_{li}(\textbf{g}_{l}-(\textbf{g}_{j}^{\mathrm{T}}\textbf{g}_{l})\textbf{g}_{j})^{\mathrm{T}}\otimes \textbf{e}_{i}^{\mathrm{T}},\textbf{0}_{1\times n}\Big]^{*},
\end{eqnarray} where $\alpha_{li}\in\mathbb{C}$ and not all $\alpha_{li}$ are zero. Therefore,
$\ker\Big(A_{k}^{[j]}+h_{k}A_{{\mathrm{c}}k}-\lambda_{3} I_{2nq+n}\Big)$ is given by (\ref{egns6}).

If $\frac{\mu_{k}}{\kappa_{k}}(\kappa_{k} h_{k}-1)+\eta_{k}\neq 0$ and $\kappa_{k} h_{k}-1-\kappa_{k}=0$, then $h_{k}=1+\frac{1}{\kappa_{k}}$. Clearly $h_{k}\kappa_{k}\neq 1$. In this case, since $\lambda=-\kappa_{k}$, it follows that
\begin{eqnarray*}
&&\det\Big[\Big(\frac{\mu_{k}}{\lambda}+\mu_{k} h_{k}+\eta_{k}\Big)(L_{k}\otimes I_{n})+\Big(\frac{\kappa_{k}}{\lambda}+\lambda+\kappa_{k} h_{k}\Big)I_{nq}\Big]\nonumber\\
&&=\det\Big[\Big(-\frac{\mu_{k}}{\kappa_{k}}+\mu_{k} h_{k}+\eta_{k}\Big)(L_{k}\otimes I_{n})-\kappa_{k}I_{nq}\Big]\nonumber\\
&&=\kappa_{k}^{nq}\det\Big[\frac{\mu_{k}+\eta_{k}}{\kappa_{k}}(L_{k}\otimes I_{n})-I_{nq}\Big].
\end{eqnarray*} Hence, $\det\Big[\Big(\frac{\mu_{k}}{\lambda}+\mu_{k} h_{k}+\eta_{k}\Big)(L_{k}\otimes I_{n})+\Big(\frac{\kappa_{k}}{\lambda}+\lambda+\kappa_{k} h_{k}\Big)I_{nq}\Big]=0$ if and only if $1\in{\mathrm{spec}}(\frac{\mu_{k}+\eta_{k}}{\kappa_{k}}L_{k})$. Note that $1\in{\mathrm{spec}}(\frac{\mu_{k}+\eta_{k}}{\kappa_{k}}L_{k})$ implies that $\mu_{k}+\eta_{k}\neq0$ and hence, by using $\kappa_{k} h_{k}-1-\kappa_{k}=0$, $\frac{\mu_{k}}{\kappa_{k}}(\kappa_{k} h_{k}-1)+\eta_{k}=\mu_{k}+\eta_{k}\neq0$. Now we assume that
$1\in{\mathrm{spec}}(\frac{\mu_{k}+\eta_{k}}{\kappa_{k}}L_{k})$ and $h_{k}=1+\frac{1}{\kappa_{k}}$. Next, since $\kappa_{k} h_{k}-1-\kappa_{k}=0$ and $\mu_{k}+\eta_{k}\neq0$, it follows from (\ref{Enqu1}) that 
\begin{eqnarray}\label{u1u0uq}
(L_{k}\otimes I_{n})[\textbf{u}_{1}^{*},\ldots,\textbf{u}_{q}^{*}]^{*}=\frac{\kappa_{k}}{\mu_{k}+\eta_{k}}(\textbf{1}_{q\times 1}\otimes I_{n})\textbf{u}_{0}.
\end{eqnarray} Note that $(L_{k}\otimes I_{n})(\textbf{1}_{q\times 1}\otimes I_{n})=\textbf{0}_{nq\times n}$. Pre-multiplying $L_{k}\otimes I_{n}$ on both sides of (\ref{u1u0uq}) yields $(L_{k}\otimes I_{n})(L_{k}\otimes I_{n})[\textbf{u}_{1}^{*},\ldots,\textbf{u}_{q}^{*}]^{*}=\textbf{0}_{nq\times 1}$, which implies that $(L_{k}\otimes I_{n})[\textbf{u}_{1}^{*},\ldots,\textbf{u}_{q}^{*}]^{*}\in\ker(L_{k}\otimes I_{n})$. Using the similar arguments as in the proof of Case 2 of $ii$) in Lemma~\ref{lemma_Arank}, it follows that
\begin{eqnarray}\label{lku1}
(L_{k}\otimes I_{n})[\textbf{u}_{1}^{*},\ldots,\textbf{u}_{q}^{*}]^{*}=\sum_{l=0}^{q-1-{\mathrm{rank}}(L_{k})}\sum_{i=1}^{n}\alpha_{li}\textbf{w}_{l}\otimes\textbf{e}_{i},
\end{eqnarray} where $\alpha_{li}\in\mathbb{C}$. Let $\textbf{u}_{0}=\sum_{i=1}^{n}\beta_{i}\textbf{e}_{i}$, where $\beta_{i}\in\mathbb{C}$. Then it follows from $iii$) of Lemma~\ref{lemma_EW} that $(\textbf{1}_{q\times 1}\otimes I_{n})\textbf{u}_{0}=\sum_{i=1}^{n}\beta_{i}(\textbf{1}_{q\times 1}\otimes I_{n})\textbf{e}_{i}=\sum_{i=1}^{n}\beta_{i}(\textbf{w}_{0}\otimes\textbf{e}_{i})$. Now it follows from (\ref{u1u0uq}) and (\ref{lku1}) that 
\begin{eqnarray*}
\sum_{i=1}^{n}\Big(\alpha_{0i}-\beta_{i}\frac{\kappa_{k}}{\mu_{k}+\eta_{k}}\Big)\textbf{w}_{0}\otimes\textbf{e}_{i}+\sum_{l=1}^{q-1-{\mathrm{rank}}(L_{k})}\sum_{i=1}^{n}\alpha_{li}\textbf{w}_{l}\otimes\textbf{e}_{i}=\textbf{0}_{nq\times 1},
\end{eqnarray*} which implies that $\alpha_{0i}-\beta_{i}\frac{\kappa_{k}}{\mu_{k}+\eta_{k}}=0$ and $\alpha_{li}=0$ for every $i=1,\ldots,n$ and every $l=1,\ldots,q-1-{\mathrm{rank}}(L_{k})$. Hence, 
\begin{eqnarray}\label{Axbeqn}
(L_{k}\otimes I_{n})[\textbf{u}_{1}^{*},\ldots,\textbf{u}_{q}^{*}]^{*}=\frac{\kappa_{k}}{\mu_{k}+\eta_{k}}\sum_{i=1}^{n}\beta_{i}\textbf{w}_{0}\otimes\textbf{e}_{i}.
\end{eqnarray} Together with $E_{n\times nq}^{[j]}[\textbf{u}_{1}^{*},\ldots,\textbf{u}_{q}^{*}]^{*}=(\textbf{g}_{j}^{\mathrm{T}}\otimes I_{n})[\textbf{u}_{1}^{*},\ldots,\textbf{u}_{q}^{*}]^{*}=\textbf{0}_{n\times 1}$, we have
\begin{eqnarray}\label{Axbeqnab}
\small\left[\begin{array}{c}
L_{k}\otimes I_{n}\\
\textbf{g}_{j}^{\mathrm{T}}\otimes I_{n}\\
\end{array}\right][\textbf{u}_{1}^{*},\ldots,\textbf{u}_{q}^{*}]^{*}=\small\left[\begin{array}{c}
\frac{\kappa_{k}}{\mu_{k}+\eta_{k}}\sum_{i=1}^{n}\beta_{i}\textbf{w}_{0}\otimes\textbf{e}_{i} \\
\textbf{0}_{n\times 1}\\
\end{array}\right].
\end{eqnarray} Now it follows from $ii$) of Theorem 2.6.4 of \cite[p.~108]{Bernstein:2009} that (\ref{Axbeqnab}) has a solution $[\textbf{u}_{1}^{*},\ldots,\textbf{u}_{q}^{*}]^{*}$ if and only if
\begin{eqnarray}\label{rankconab}
{\mathrm{rank}}\small\left[\begin{array}{c}
L_{k}\otimes I_{n}\\
\textbf{g}_{j}^{\mathrm{T}}\otimes I_{n}\\
\end{array}\right]={\mathrm{rank}}\small\left[\begin{array}{cc}
L_{k}\otimes I_{n} & \frac{\kappa_{k}}{\mu_{k}+\eta_{k}}\sum_{i=1}^{n}\beta_{i}\textbf{w}_{0}\otimes\textbf{e}_{i}\\
\textbf{g}_{j}^{\mathrm{T}}\otimes I_{n} & \textbf{0}_{n\times 1}\\
\end{array}\right].
\end{eqnarray} We claim that (\ref{rankconab}) is indeed true. First, if $\beta_{i}=0$ for every $i=1,\ldots,n$, then it is clear that ${\mathrm{rank}}\small\left[\begin{array}{c}
L_{k}\otimes I_{n}\\
\textbf{g}_{j}^{\mathrm{T}}\otimes I_{n}\\
\end{array}\right]={\mathrm{rank}}\small\left[\begin{array}{cc}
L_{k}\otimes I_{n} & \textbf{0}_{nq\times 1}\\
\textbf{g}_{j}^{\mathrm{T}}\otimes I_{n} & \textbf{0}_{n\times 1}\\
\end{array}\right]$. Alternatively, assume that $\beta_{i}\neq0$ for some $i\in\{1,\ldots,n\}$. Note that it follows from Fact 2.11.8 of \cite[p.~132]{Bernstein:2009} that ${\mathrm{rank}}\small\left[\begin{array}{c}
L_{k}\otimes I_{n}\\
\textbf{g}_{j}^{\mathrm{T}}\otimes I_{n}\\
\end{array}\right]\leq{\mathrm{rank}}\small\left[\begin{array}{cc}
L_{k}\otimes I_{n} & \frac{\kappa_{k}}{\mu_{k}+\eta_{k}}\sum_{i=1}^{n}\beta_{i}\textbf{w}_{0}\otimes\textbf{e}_{i}\\
\textbf{g}_{j}^{\mathrm{T}}\otimes I_{n} & \textbf{0}_{n\times 1}\\
\end{array}\right]$. To show (\ref{rankconab}), it suffices to show that 
\begin{eqnarray*}
{\mathrm{def}}\small\left[\begin{array}{c}
L_{k}\otimes I_{n}\\
\textbf{g}_{j}^{\mathrm{T}}\otimes I_{n}\\
\end{array}\right]\leq{\mathrm{def}}\small\left[\begin{array}{cc}
L_{k}\otimes I_{n} & \frac{\kappa_{k}}{\mu_{k}+\eta_{k}}\sum_{i=1}^{n}\beta_{i}\textbf{w}_{0}\otimes\textbf{e}_{i}\\
\textbf{g}_{j}^{\mathrm{T}}\otimes I_{n} & \textbf{0}_{n\times 1}\\
\end{array}\right],
\end{eqnarray*} or, equivalently,
\begin{eqnarray*}
\dim\ker\small\left[\begin{array}{c}
L_{k}\otimes I_{n}\\
\textbf{g}_{j}^{\mathrm{T}}\otimes I_{n}\\
\end{array}\right]\leq\dim\ker\small\left[\begin{array}{cc}
L_{k}\otimes I_{n} & \frac{\kappa_{k}}{\mu_{k}+\eta_{k}}\sum_{i=1}^{n}\beta_{i}\textbf{w}_{0}\otimes\textbf{e}_{i}\\
\textbf{g}_{j}^{\mathrm{T}}\otimes I_{n} & \textbf{0}_{n\times 1}\\
\end{array}\right].
\end{eqnarray*}
 Let $s\in\mathbb{C}$ be such that $s\in\ker\small\left[\begin{array}{c}
 \frac{\kappa_{k}}{\mu_{k}+\eta_{k}}\sum_{i=1}^{n}\beta_{i}\textbf{w}_{0}\otimes\textbf{e}_{i}\\
 \textbf{0}_{n\times 1}\\
 \end{array}\right]$. Then $s\frac{\kappa_{k}}{\mu_{k}+\eta_{k}}\beta_{i}=0$ for some $i\in\{1,\ldots,n\}$, which implies that $s=0$. Thus, $\dim\ker\small\left[\begin{array}{c}
  \frac{\kappa_{k}}{\mu_{k}+\eta_{k}}\sum_{i=1}^{n}\beta_{i}\textbf{w}_{0}\otimes\textbf{e}_{i}\\
  \textbf{0}_{n\times 1}\\
  \end{array}\right]=0$. Consequently, it follows from Fact 2.11.8 of \cite[p.~132]{Bernstein:2009} that
 \begin{eqnarray*}
 \dim\ker\small\left[\begin{array}{c}
 L_{k}\otimes I_{n}\\
 \textbf{g}_{j}^{\mathrm{T}}\otimes I_{n}\\
 \end{array}\right]&=&\dim\ker\small\left[\begin{array}{c}
 L_{k}\otimes I_{n}\\
 \textbf{g}_{j}^{\mathrm{T}}\otimes I_{n}\\
 \end{array}\right]+\dim\ker\small\left[\begin{array}{c}
  \frac{\kappa_{k}}{\mu_{k}+\eta_{k}}\sum_{i=1}^{n}\beta_{i}\textbf{w}_{0}\otimes\textbf{e}_{i}\\
  \textbf{0}_{n\times 1}\\
  \end{array}\right]\nonumber\\
 &\leq&\dim\ker\small\left[\begin{array}{cc}
 L_{k}\otimes I_{n} & \frac{\kappa_{k}}{\mu_{k}+\eta_{k}}\sum_{i=1}^{n}\beta_{i}\textbf{w}_{0}\otimes\textbf{e}_{i}\\
 \textbf{g}_{j}^{\mathrm{T}}\otimes I_{n} & \textbf{0}_{n\times 1}\\
 \end{array}\right],
 \end{eqnarray*} which implies that  ${\mathrm{rank}}\small\left[\begin{array}{c}
 L_{k}\otimes I_{n}\\
 \textbf{g}_{j}^{\mathrm{T}}\otimes I_{n}\\
 \end{array}\right]\geq{\mathrm{rank}}\small\left[\begin{array}{cc}
 L_{k}\otimes I_{n} & \frac{\kappa_{k}}{\mu_{k}+\eta_{k}}\sum_{i=1}^{n}\beta_{i}\textbf{w}_{0}\otimes\textbf{e}_{i}\\
 \textbf{g}_{j}^{\mathrm{T}}\otimes I_{n} & \textbf{0}_{n\times 1}\\
 \end{array}\right]$. Hence, (\ref{rankconab}) holds.
Next, it follows from $vi$) of Proposition 6.1.7 of \cite[p.~400]{Bernstein:2009} and $viii)$ of Proposition 6.1.6 of \cite[p.~399]{Bernstein:2009} that the general solution to (\ref{Axbeqnab}) is given by the form
\begin{eqnarray}\label{usolutionab} 
[\textbf{u}_{1}^{*},\ldots,\textbf{u}_{q}^{*}]^{*}&=&\small\left[\begin{array}{c}
 L_{k}\otimes I_{n}\\
 \textbf{g}_{j}^{\mathrm{T}}\otimes I_{n}\\
 \end{array}\right]^{+}\small\left[\begin{array}{c}
   \frac{\kappa_{k}}{\mu_{k}+\eta_{k}}\sum_{i=1}^{n}\beta_{i}\textbf{w}_{0}\otimes\textbf{e}_{i}\\
   \textbf{0}_{n\times 1}\\
   \end{array}\right]+\sum_{l=1}^{q}\sum_{i=1}^{n}\gamma_{li}\Big(I_{nq}-\small\left[\begin{array}{c}
    L_{k}\otimes I_{n}\\
    \textbf{g}_{j}^{\mathrm{T}}\otimes I_{n}\\
    \end{array}\right]^{+}\small\left[\begin{array}{c}
     L_{k}\otimes I_{n}\\
     \textbf{g}_{j}^{\mathrm{T}}\otimes I_{n}\\
     \end{array}\right]\Big)\nonumber\\
     &&(\textbf{g}_{l}\otimes\textbf{e}_{i})\nonumber\\
&=&\Big(\small\left[\begin{array}{c}
 L_{k}\\
 \textbf{g}_{j}^{\mathrm{T}}\\
 \end{array}\right]\otimes I_{n}\Big)^{+}\small\left[\begin{array}{c}
   \frac{\kappa_{k}}{\mu_{k}+\eta_{k}}\sum_{i=1}^{n}\beta_{i}\textbf{w}_{0}\otimes\textbf{e}_{i}\\
   \sum_{i=1}^{n}0\otimes\textbf{e}_{i}\\
   \end{array}\right]+\sum_{l=1}^{q}\sum_{i=1}^{n}\gamma_{li}\Big(I_{nq}-\Big(\small\left[\begin{array}{c}
    L_{k}\\
    \textbf{g}_{j}^{\mathrm{T}}\\
    \end{array}\right]\otimes I_{n}\Big)^{+}\nonumber\\
    &&\Big(\small\left[\begin{array}{c}
     L_{k}\\
     \textbf{g}_{j}^{\mathrm{T}}\\
     \end{array}\right]\otimes I_{n}\Big)\Big)(\textbf{g}_{l}\otimes\textbf{e}_{i})\nonumber\\
&=&\Big(\small\left[\begin{array}{c}
      L_{k}\\
      \textbf{g}_{j}^{\mathrm{T}}\\
      \end{array}\right]^{+}\otimes I_{n}\Big)\Big(\sum_{i=1}^{n}\small\left[\begin{array}{c}
        \frac{\kappa_{k}}{\mu_{k}+\eta_{k}}\beta_{i}\textbf{w}_{0}\\
        0\\
        \end{array}\right]\otimes\textbf{e}_{i}\Big)+\sum_{l=1}^{q}\sum_{i=1}^{n}\gamma_{li}\Big(I_{q}\otimes I_{n}-\Big(\small\left[\begin{array}{c}
         L_{k}\\
         \textbf{g}_{j}^{\mathrm{T}}\\
         \end{array}\right]^{+}\otimes I_{n}\Big)\nonumber\\
         &&\Big(\small\left[\begin{array}{c}
          L_{k}\\
          \textbf{g}_{j}^{\mathrm{T}}\\
          \end{array}\right]\otimes I_{n}\Big)\Big)(\textbf{g}_{l}\otimes\textbf{e}_{i})\nonumber\\
          &=&\sum_{i=1}^{n}\Big(\small\left[\begin{array}{c}
                L_{k}\\
                \textbf{g}_{j}^{\mathrm{T}}\\
                \end{array}\right]^{+}\small\left[\begin{array}{c}
                        \frac{\kappa_{k}}{\mu_{k}+\eta_{k}}\beta_{i}\textbf{w}_{0}\\
                        0\\
                        \end{array}\right]\Big)\otimes\textbf{e}_{i}+\sum_{l=1}^{q}\sum_{i=1}^{n}\gamma_{li}\Big(I_{q}\otimes I_{n}-\Big(\small\left[\begin{array}{c}
                                 L_{k}\\
                                 \textbf{g}_{j}^{\mathrm{T}}\\
                                 \end{array}\right]^{+}\small\left[\begin{array}{c}
                                                                   L_{k}\\
                                                                   \textbf{g}_{j}^{\mathrm{T}}\\
                                                                   \end{array}\right]\otimes I_{n}\Big)\Big)\nonumber\\
                                 &&(\textbf{g}_{l}\otimes\textbf{e}_{i})\nonumber\\
          &=&\sum_{i=1}^{n}\Big(\small\left[\begin{array}{c}
                L_{k}\\
                \textbf{g}_{j}^{\mathrm{T}}\\
                \end{array}\right]^{+}\small\left[\begin{array}{c}
                        \frac{\kappa_{k}}{\mu_{k}+\eta_{k}}\beta_{i}\textbf{w}_{0}\\
                        0\\
                        \end{array}\right]\Big)\otimes\textbf{e}_{i}+\sum_{l=1}^{q}\sum_{i=1}^{n}\gamma_{li}\Big(\textbf{g}_{l}-\small\left[\begin{array}{c}
                                 L_{k}\\
                                 \textbf{g}_{j}^{\mathrm{T}}\\
                                 \end{array}\right]^{+}\small\left[\begin{array}{c}
                                                                   L_{k}\\
                                                                   \textbf{g}_{j}^{\mathrm{T}}\\
                                                                   \end{array}\right]\textbf{g}_{l}\Big)\otimes\textbf{e}_{i},
\end{eqnarray} where $\gamma_{li}\in\mathbb{C}$. Note that by Proposition 6.1.6 of \cite[p.~399]{Bernstein:2009}, $L_{k}^{\mathrm{T}}(L_{k}^{\mathrm{T}})^{+}=L_{k}^{\mathrm{T}}(L_{k}^{+})^{\mathrm{T}}=(L_{k}^{+}L_{k})^{\mathrm{T}}=L_{k}^{+}L_{k}$. It follows from Fact 6.5.17 of \cite[p.~427]{Bernstein:2009} that 
\begin{eqnarray}
\small\left[\begin{array}{c}
L_{k}\\
\textbf{g}_{j}^{\mathrm{T}}
\end{array}\right]^{+}=\small\left[\begin{array}{cc}
L_{k}^{+}(I_{q}-\varphi_{k}\textbf{g}_{j}^{\mathrm{T}}) & \varphi_{k}
\end{array}\right],
\end{eqnarray} where $\varphi_{k}$ is given by (\ref{varphik}). Note that $\textbf{g}_{j}^{\mathrm{T}}\textbf{w}_{0}=1$ for every $j=1,\ldots,q$. Hence, it follows that for every $i=1,\ldots,n$ and every $j,l=1,\ldots,q$,
\begin{eqnarray}
\small\left[\begin{array}{cc}
L_{k}^{+}(I_{q}-\varphi_{k}\textbf{g}_{j}^{\mathrm{T}}) & \varphi_{k}
\end{array}\right]\small\left[\begin{array}{c}
                        \frac{\kappa_{k}}{\mu_{k}+\eta_{k}}\beta_{i}\textbf{w}_{0}\\
                        0\\
                        \end{array}\right]&=&\frac{\kappa_{k}}{\mu_{k}+\eta_{k}}\beta_{i}L_{k}^{+}\textbf{w}_{0}-\frac{\kappa_{k}}{\mu_{k}+\eta_{k}}\beta_{i}L_{k}^{+}\varphi_{k},\\
\textbf{g}_{l}-\small\left[\begin{array}{c}
L_{k}\\
\textbf{g}_{j}^{\mathrm{T}}
\end{array}\right]^{+}\small\left[\begin{array}{c}
L_{k}\\
\textbf{g}_{j}^{\mathrm{T}}
\end{array}\right]\textbf{g}_{l}&=&\textbf{g}_{l}-\small\left[\begin{array}{cc}
L_{k}^{+}(I_{q}-\varphi_{k}\textbf{g}_{j}^{\mathrm{T}}) & \varphi_{k}
\end{array}\right]\small\left[\begin{array}{c}
L_{k}\\
\textbf{g}_{j}^{\mathrm{T}}
\end{array}\right]\textbf{g}_{l}\nonumber\\
&=&\textbf{g}_{l}-\small\left[\begin{array}{cc}
L_{k}^{+}(I_{q}-\varphi_{k}\textbf{g}_{j}^{\mathrm{T}}) & \varphi_{k}
\end{array}\right]\small\left[\begin{array}{c}
L_{k}\textbf{g}_{l}\\
\textbf{g}_{j}^{\mathrm{T}}\textbf{g}_{l}
\end{array}\right]\nonumber\\
&=&\textbf{g}_{l}-L_{k}^{+}(I_{q}-\varphi_{k}\textbf{g}_{j}^{\mathrm{T}})L_{k}\textbf{g}_{l}-(\textbf{g}_{j}^{\mathrm{T}}\textbf{g}_{l})\varphi_{k}\nonumber\\
&=&\textbf{g}_{l}-L_{k}^{+}L_{k}\textbf{g}_{l}+(\textbf{g}_{j}^{\mathrm{T}}L_{k}\textbf{g}_{l})L_{k}^{+}\varphi_{k}-(\textbf{g}_{j}^{\mathrm{T}}\textbf{g}_{l})\varphi_{k}.
\end{eqnarray} Then (\ref{usolutionab}) becomes
\begin{eqnarray}\label{gsolutionab}
[\textbf{u}_{1}^{*},\ldots,\textbf{u}_{q}^{*}]^{*}&=&\frac{\kappa_{k}}{\mu_{k}+\eta_{k}}\sum_{i=1}^{n}\beta_{i}L_{k}^{+}\textbf{w}_{0}\otimes\textbf{e}_{i}-\frac{\kappa_{k}}{\mu_{k}+\eta_{k}}\sum_{i=1}^{n}\beta_{i}L_{k}^{+}\varphi_{k}\otimes\textbf{e}_{i}\nonumber\\
&&+\sum_{l=1}^{q}\sum_{i=1}^{n}\gamma_{li}(\textbf{g}_{l}-L_{k}^{+}L_{k}\textbf{g}_{l}+(\textbf{g}_{j}^{\mathrm{T}}L_{k}\textbf{g}_{l})L_{k}^{+}\varphi_{k}-(\textbf{g}_{j}^{\mathrm{T}}\textbf{g}_{l})\varphi_{k})\otimes \textbf{e}_{i}.
\end{eqnarray} 

In summary, if $1\in{\mathrm{spec}}(\frac{\mu_{k}+\eta_{k}}{\kappa_{k}}L_{k})$ and $h_{k}=1+\frac{1}{\kappa_{k}}$, then $\lambda=-\kappa_{k}$ is indeed an eigenvalue of $A_{k}^{[j]}+h_{k}A_{{\mathrm{c}}k}$. In this case, $\textbf{x}_{1}=\textbf{0}_{nq\times 1}$, $\textbf{x}_{2}=[\textbf{u}_{1}^{*},\ldots,\textbf{u}_{q}^{*}]^{*}$ given by (\ref{gsolutionab}), and $\textbf{x}_{3}=\sum_{i=1}^{n}\beta_{i}\textbf{e}_{i}$, where not all of $\beta_{i}$ and $\gamma_{li}$ are zero. The corresponding eigenvectors for $\lambda_{3}$ are given by
\begin{eqnarray}\label{xw0ab}
\textbf{x}&=&\Big[\textbf{0}_{1\times nq},\frac{\kappa_{k}}{\mu_{k}+\eta_{k}}\sum_{i=1}^{n}\beta_{i}(L_{k}^{+}\textbf{w}_{0}\otimes\textbf{e}_{i})^{\mathrm{T}}-\frac{\kappa_{k}}{\mu_{k}+\eta_{k}}\sum_{i=1}^{n}\beta_{i}(L_{k}^{+}\varphi_{k}\otimes\textbf{e}_{i})^{\mathrm{T}}\nonumber\\
&&+\sum_{l=1}^{q}\sum_{i=1}^{n}\gamma_{li}(\textbf{g}_{l}-L_{k}^{+}L_{k}\textbf{g}_{l}+(\textbf{g}_{j}^{\mathrm{T}}L_{k}\textbf{g}_{l})L_{k}^{+}\varphi_{k}-(\textbf{g}_{j}^{\mathrm{T}}\textbf{g}_{l})\varphi_{k})^{\mathrm{T}}\otimes \textbf{e}_{i}^{\mathrm{T}},\sum_{i=1}^{n}\beta_{i}\textbf{e}_{i}^{\mathrm{T}}\Big]^{*},
\end{eqnarray} where $\beta_{i}\in\mathbb{C}$ and $\gamma_{li}\in\mathbb{C}$ and not all of them are zero. Therefore,
$\ker\Big(A_{k}^{[j]}+h_{k}A_{{\mathrm{c}}k}-\lambda_{3} I_{2nq+n}\Big)$ is given by (\ref{egns7}).

If $\frac{\mu_{k}}{\kappa_{k}}(\kappa_{k} h_{k}-1)+\eta_{k}\neq 0$, $\kappa_{k} h_{k}-1-\kappa_{k}\neq0$, and $\kappa_{k} h_{k}-1\neq0$, in this case, since $\lambda=-\kappa_{k}$, then it follows that
\begin{eqnarray*}
&&\det\Big[\Big(\frac{\mu_{k}}{\lambda}+\mu_{k} h_{k}+\eta_{k}\Big)(L_{k}\otimes I_{n})+\Big(\frac{\kappa_{k}}{\lambda}+\lambda+\kappa_{k} h_{k}\Big)I_{nq}\Big]\nonumber\\
&&=\det\Big[\Big(\frac{\mu_{k}}{\kappa_{k}}(\kappa_{k} h_{k}-1)+\eta_{k}\Big)(L_{k}\otimes I_{n})+(\kappa_{k} h_{k}-1-\kappa_{k})I_{nq}\Big]\nonumber\\
&&=(-\kappa_{k} h_{k}+1+\kappa_{k})^{nq}\det\Big[\frac{\mu_{k}(\kappa_{k}h_{k}-1)+\eta_{k}\kappa_{k}}{\kappa_{k}(-\kappa_{k} h_{k}+1+\kappa_{k})}(L_{k}\otimes I_{n})-I_{nq}\Big].
\end{eqnarray*} Hence, $\det\Big[\Big(\frac{\mu_{k}}{\lambda}+\mu_{k} h_{k}+\eta_{k}\Big)(L_{k}\otimes I_{n})+\Big(\frac{\kappa_{k}}{\lambda}+\lambda+\kappa_{k} h_{k}\Big)I_{nq}\Big]=0$ if and only if $1\in{\mathrm{spec}}(\frac{\mu_{k}(\kappa_{k}h_{k}-1)+\eta_{k}\kappa_{k}}{\kappa_{k}(-\kappa_{k} h_{k}+1+\kappa_{k})}L_{k})$. Again, note that $1\in{\mathrm{spec}}(\frac{\mu_{k}(\kappa_{k}h_{k}-1)+\eta_{k}\kappa_{k}}{\kappa_{k}(-\kappa_{k} h_{k}+1+\kappa_{k})}L_{k})$ implies that $\frac{\mu_{k}}{\kappa_{k}}(\kappa_{k} h_{k}-1)+\eta_{k}\neq 0$ and $\kappa_{k} h_{k}-1-\kappa_{k}\neq0$. Now we assume that
$1\in{\mathrm{spec}}(\frac{\mu_{k}(\kappa_{k}h_{k}-1)+\eta_{k}\kappa_{k}}{\kappa_{k}(-\kappa_{k} h_{k}+1+\kappa_{k})}L_{k})$ and $\kappa_{k}h_{k}\neq 1$. Next, let $\textbf{u}_{0}=\sum_{i=1}^{n}\beta_{i}\textbf{e}_{i}$, where $\beta_{i}\in\mathbb{C}$ and it follows from (\ref{Enqu1}) that 
\begin{eqnarray}\label{xabs}
\Big(\Big(\frac{\mu_{k}}{\kappa_{k}}(\kappa_{k} h_{k}-1)+\eta_{k}\Big)(L_{k}\otimes I_{n})+(\kappa_{k} h_{k}-1-\kappa_{k})I_{nq}\Big)[\textbf{u}_{1}^{*},\ldots,\textbf{u}_{q}^{*}]^{*}=\kappa_{k}\sum_{i=1}^{n}\beta_{i}\textbf{1}_{q\times 1}\otimes\textbf{e}_{i}.
\end{eqnarray} Note that a specific solution $[\textbf{u}_{1}^{*},\ldots,\textbf{u}_{q}^{*}]^{*}$ to (\ref{xabs}) is given by the form
\begin{eqnarray}\label{specific}
[\textbf{u}_{1}^{*},\ldots,\textbf{u}_{q}^{*}]^{*}=\frac{\kappa_{k}}{\kappa_{k}h_{k}-1-\kappa_{k}}\sum_{i=1}^{n}\beta_{i}\textbf{1}_{q\times 1}\otimes\textbf{e}_{i}.
\end{eqnarray} Substituting (\ref{specific}) into (\ref{Enqu}) by using $iii$) of Lemma~\ref{lemma_EW} yields $\frac{\kappa_{k}(\kappa_{k}h_{k}-1)}{\kappa_{k}h_{k}-1-\kappa_{k}}\sum_{i=1}^{n}\beta_{i}E_{n\times nq}^{[j]}(\textbf{1}_{q\times 1}\otimes\textbf{e}_{i})=\frac{\kappa_{k}(\kappa_{k}h_{k}-1)}{\kappa_{k}h_{k}-1-\kappa_{k}}\\\sum_{i=1}^{n}\beta_{i}\textbf{e}_{i}=\textbf{0}_{n\times 1}$, which implies that $\beta_{i}=0$ for every $i=1,\ldots,n$, and hence, $\textbf{u}_{0}=\textbf{0}_{n\times 1}$. Thus, (\ref{xabs}) becomes
\begin{eqnarray}\label{homo}
\Big(\Big(\frac{\mu_{k}}{\kappa_{k}}(\kappa_{k} h_{k}-1)+\eta_{k}\Big)(L_{k}\otimes I_{n})+(\kappa_{k} h_{k}-1-\kappa_{k})I_{nq}\Big)[\textbf{u}_{1}^{*},\ldots,\textbf{u}_{q}^{*}]^{*}=\textbf{0}_{nq\times 1}.
\end{eqnarray} Let $M_{k}=(\frac{\mu_{k}}{\kappa_{k}}(\kappa_{k} h_{k}-1)+\eta_{k})L_{k}+(\kappa_{k} h_{k}-1-\kappa_{k})I_{q}$. Again, note that $E_{n\times nq}^{[j]}=\textbf{g}_{j}^{\mathrm{T}}\otimes I_{n}$ for every $j=1,\ldots,q$. Then it follows from  (\ref{homo}) and (\ref{Enqu}) that
\begin{eqnarray}\label{uMg}
\small\left[\begin{array}{c}
M_{k}\otimes I_{n}\\
\textbf{g}_{j}^{\mathrm{T}}\otimes I_{n}
\end{array}\right][\textbf{u}_{1}^{*},\ldots,\textbf{u}_{q}^{*}]^{*}=\Big(\small\left[\begin{array}{c}
M_{k}\\
\textbf{g}_{j}^{\mathrm{T}}
\end{array}\right]\otimes I_{n}\Big)[\textbf{u}_{1}^{*},\ldots,\textbf{u}_{q}^{*}]^{*}=\textbf{0}_{(nq+n)\times 1}.
\end{eqnarray} Next, it follows from $vi$) of Proposition 6.1.7 of \cite[p.~400]{Bernstein:2009} and $viii$) of Proposition 6.1.6 of \cite[p.~399]{Bernstein:2009} that the general solution $[\textbf{u}_{1}^{*},\ldots,\textbf{u}_{q}^{*}]^{*}$ to (\ref{uMg}) is given by the form
\begin{eqnarray}
[\textbf{u}_{1}^{*},\ldots,\textbf{u}_{q}^{*}]^{*}&=&\Big[I_{nq}-\Big(\small\left[\begin{array}{c}
M_{k}\\
\textbf{g}_{j}^{\mathrm{T}}
\end{array}\right]\otimes I_{n}\Big)^{+}\Big(\small\left[\begin{array}{c}
M_{k}\\
\textbf{g}_{j}^{\mathrm{T}}
\end{array}\right]\otimes I_{n}\Big)\Big]\sum_{i=1}^{n}\sum_{l=1}^{q}\varpi_{li}\textbf{g}_{l}\otimes\textbf{e}_{i}\nonumber\\
&=&\Big[I_{nq}-\Big(\small\left[\begin{array}{c}
M_{k}\\
\textbf{g}_{j}^{\mathrm{T}}
\end{array}\right]^{+}\otimes I_{n}\Big)\Big(\small\left[\begin{array}{c}
M_{k}\\
\textbf{g}_{j}^{\mathrm{T}}
\end{array}\right]\otimes I_{n}\Big)\Big]\sum_{i=1}^{n}\sum_{l=1}^{q}\varpi_{li}\textbf{g}_{l}\otimes\textbf{e}_{i}\nonumber\\
&=&\Big[I_{q}\otimes I_{n}-\Big(\small\left[\begin{array}{c}
M_{k}\\
\textbf{g}_{j}^{\mathrm{T}}
\end{array}\right]^{+}\small\left[\begin{array}{c}
M_{k}\\
\textbf{g}_{j}^{\mathrm{T}}
\end{array}\right]\otimes I_{n}\Big)\Big]\sum_{i=1}^{n}\sum_{l=1}^{q}\varpi_{li}\textbf{g}_{l}\otimes\textbf{e}_{i}\nonumber\\
&=&\Big[\Big(I_{q}-\small\left[\begin{array}{c}
M_{k}\\
\textbf{g}_{j}^{\mathrm{T}}
\end{array}\right]^{+}\small\left[\begin{array}{c}
M_{k}\\
\textbf{g}_{j}^{\mathrm{T}}
\end{array}\right]\Big)\otimes I_{n}\Big]\sum_{i=1}^{n}\sum_{l=1}^{q}\varpi_{li}\textbf{g}_{l}\otimes\textbf{e}_{i}\nonumber\\
&=&\sum_{i=1}^{n}\sum_{l=1}^{q}\varpi_{li}\Big(\textbf{g}_{l}-\small\left[\begin{array}{c}
M_{k}\\
\textbf{g}_{j}^{\mathrm{T}}
\end{array}\right]^{+}\small\left[\begin{array}{c}
M_{k}\\
\textbf{g}_{j}^{\mathrm{T}}
\end{array}\right]\textbf{g}_{l}\Big)\otimes \textbf{e}_{i},\label{gsolution1}
\end{eqnarray} where $\varpi_{li}\in\mathbb{C}$ and $j=1,\ldots,q$. Note that by Proposition 6.1.6 of \cite[p.~399]{Bernstein:2009}, $M_{k}^{\mathrm{T}}(M_{k}^{\mathrm{T}})^{+}=M_{k}^{\mathrm{T}}(M_{k}^{+})^{\mathrm{T}}=(M_{k}^{+}M_{k})^{\mathrm{T}}=M_{k}^{+}M_{k}$. It follows from Fact 6.5.17 of \cite[p.~427]{Bernstein:2009} that 
\begin{eqnarray}
\small\left[\begin{array}{c}
M_{k}\\
\textbf{g}_{j}^{\mathrm{T}}
\end{array}\right]^{+}=\small\left[\begin{array}{cc}
M_{k}^{+}(I_{q}-\phi_{k}\textbf{g}_{j}^{\mathrm{T}}) & \phi_{k}
\end{array}\right],
\end{eqnarray} where $\phi_{k}$ is given by (\ref{phik}). Hence, it follows that for every $j,l=1,\ldots,q$,
\begin{eqnarray}
\textbf{g}_{l}-\small\left[\begin{array}{c}
M_{k}\\
\textbf{g}_{j}^{\mathrm{T}}
\end{array}\right]^{+}\small\left[\begin{array}{c}
M_{k}\\
\textbf{g}_{j}^{\mathrm{T}}
\end{array}\right]\textbf{g}_{l}&=&\textbf{g}_{l}-\small\left[\begin{array}{cc}
M_{k}^{+}(I_{q}-\phi_{k}\textbf{g}_{j}^{\mathrm{T}}) & \phi_{k}
\end{array}\right]\small\left[\begin{array}{c}
M_{k}\\
\textbf{g}_{j}^{\mathrm{T}}
\end{array}\right]\textbf{g}_{l}\nonumber\\
&=&\textbf{g}_{l}-\small\left[\begin{array}{cc}
M_{k}^{+}(I_{q}-\phi_{k}\textbf{g}_{j}^{\mathrm{T}}) & \phi_{k}
\end{array}\right]\small\left[\begin{array}{c}
M_{k}\textbf{g}_{l}\\
\textbf{g}_{j}^{\mathrm{T}}\textbf{g}_{l}
\end{array}\right]\nonumber\\
&=&\textbf{g}_{l}-M_{k}^{+}(I_{q}-\phi_{k}\textbf{g}_{j}^{\mathrm{T}})M_{k}\textbf{g}_{l}-(\textbf{g}_{j}^{\mathrm{T}}\textbf{g}_{l})\phi_{k}\nonumber\\
&=&\textbf{g}_{l}-M_{k}^{+}M_{k}\textbf{g}_{l}+(\textbf{g}_{j}^{\mathrm{T}}M_{k}\textbf{g}_{l})M_{k}^{+}\phi_{k}-(\textbf{g}_{j}^{\mathrm{T}}\textbf{g}_{l})\phi_{k}.
\end{eqnarray} Thus, (\ref{gsolution1}) becomes 
\begin{eqnarray}\label{gso}
[\textbf{u}_{1}^{*},\ldots,\textbf{u}_{q}^{*}]^{*}=\sum_{i=1}^{n}\sum_{l=1}^{q}\varpi_{li}\Big(\textbf{g}_{l}-M_{k}^{+}M_{k}\textbf{g}_{l}+(\textbf{g}_{j}^{\mathrm{T}}M_{k}\textbf{g}_{l})M_{k}^{+}\phi_{k}-(\textbf{g}_{j}^{\mathrm{T}}\textbf{g}_{l})\phi_{k}\Big)\otimes \textbf{e}_{i}.
\end{eqnarray} 

In summary, if $1\in{\mathrm{spec}}(\frac{\mu_{k}(\kappa_{k}h_{k}-1)+\eta_{k}\kappa_{k}}{\kappa_{k}(-\kappa_{k} h_{k}+1+\kappa_{k})}L_{k})$ and $\kappa_{k}h_{k}\neq 1$, then $\lambda=-\kappa_{k}$ is indeed an eigenvalue of $A_{k}^{[j]}+h_{k}A_{{\mathrm{c}}k}$. In this case, $\textbf{x}_{1}=\textbf{0}_{nq\times 1}$, $\textbf{x}_{2}=[\textbf{u}_{1}^{*},\ldots,\textbf{u}_{q}^{*}]^{*}$ given by (\ref{gso}), and $\textbf{x}_{3}=\textbf{0}_{n\times 1}$, where not all of $\varpi_{li}$ are zero. The corresponding eigenvectors for $\lambda_{3}$ are given by
\begin{eqnarray}
\textbf{x}&=&\Big[\textbf{0}_{1\times nq},\sum_{i=1}^{n}\sum_{l=1}^{q}\varpi_{li}\Big(\textbf{g}_{l}-M_{k}^{+}M_{k}\textbf{g}_{l}+(\textbf{g}_{j}^{\mathrm{T}}M_{k}\textbf{g}_{l})M_{k}^{+}\phi_{k}-(\textbf{g}_{j}^{\mathrm{T}}\textbf{g}_{l})\phi_{k}\Big)^{\mathrm{T}}\otimes\textbf{e}_{i}^{\mathrm{T}},\textbf{0}_{1\times n}\Big]^{*},
\end{eqnarray} where $\varpi_{li}\in\mathbb{C}$ and not all of them are zero. Therefore, $\ker\Big(A_{k}^{[j]}+h_{k}A_{{\mathrm{c}}k}-\lambda_{3} I_{2nq+n}\Big)$ is given by (\ref{egns8}).
\end{IEEEproof}

\begin{remark}
One can obtain an alternative expression of (\ref{egns8}) by using the following method. First, it follows from $ii$) of Theorem 2.6.4 of \cite[p.~108]{Bernstein:2009} that (\ref{Axbeqn}) has a solution $[\textbf{u}_{1}^{*},\ldots,\textbf{u}_{q}^{*}]^{*}$ if and only if
\begin{eqnarray}\label{rankcon}
{\mathrm{rank}}(L_{k}\otimes I_{n})={\mathrm{rank}}\small\left[\begin{array}{cc}
L_{k}\otimes I_{n} & \frac{\kappa_{k}}{\mu_{k}+\eta_{k}}\sum_{i=1}^{n}\beta_{i}\textbf{w}_{0}\otimes\textbf{e}_{i}
\end{array}\right].
\end{eqnarray} We claim that (\ref{rankcon}) is indeed true. First, if $\beta_{i}=0$ for every $i=1,\ldots,n$, then it is clear that ${\mathrm{rank}}(L_{k}\otimes I_{n})={\mathrm{rank}}\small\left[\begin{array}{cc}
L_{k}\otimes I_{n} & \textbf{0}_{nq\times 1}
\end{array}\right]$. Alternatively, assume that $\beta_{i}\neq0$ for some $i\in\{1,\ldots,n\}$. Note that it follows from Fact 2.11.8 of \cite[p.~132]{Bernstein:2009} that ${\mathrm{rank}}(L_{k}\otimes I_{n})\leq{\mathrm{rank}}\small\left[\begin{array}{cc}
L_{k}\otimes I_{n} & \frac{\kappa_{k}}{\mu_{k}+\eta_{k}}\sum_{i=1}^{n}\beta_{i}\textbf{w}_{0}\otimes\textbf{e}_{i}
\end{array}\right]$. To show (\ref{rankcon}), it suffices to show that 
\begin{eqnarray*}
{\mathrm{def}}(L_{k}\otimes I_{n})\leq{\mathrm{def}}\small\left[\begin{array}{cc}
L_{k}\otimes I_{n} & \frac{\kappa_{k}}{\mu_{k}+\eta_{k}}\sum_{i=1}^{n}\beta_{i}\textbf{w}_{0}\otimes\textbf{e}_{i}
\end{array}\right],
\end{eqnarray*} or, equivalently,
\begin{eqnarray*}
\dim\ker(L_{k}\otimes I_{n})\leq\dim\ker\small\left[\begin{array}{cc}
L_{k}\otimes I_{n} & \frac{\kappa_{k}}{\mu_{k}+\eta_{k}}\sum_{i=1}^{n}\beta_{i}\textbf{w}_{0}\otimes\textbf{e}_{i}
\end{array}\right].
\end{eqnarray*}
 Let $s\in\mathbb{C}$ be such that $s\in\ker(\frac{\kappa_{k}}{\mu_{k}+\eta_{k}}\sum_{i=1}^{n}\beta_{i}\textbf{w}_{0}\otimes\textbf{e}_{i})$. Then $s\frac{\kappa_{k}}{\mu_{k}+\eta_{k}}\beta_{i}=0$ for some $i\in\{1,\ldots,n\}$, which implies that $s=0$. Thus, $\dim\ker\Big(\frac{\kappa_{k}}{\mu_{k}+\eta_{k}}\sum_{i=1}^{n}\beta_{i}\textbf{w}_{0}\otimes\textbf{e}_{i}\Big)=0$. Consequently, it follows from Fact 2.11.8 of \cite[p.~132]{Bernstein:2009} that
 \begin{eqnarray*}
 \dim\ker(L_{k}\otimes I_{n})&=&\dim\ker(L_{k}\otimes I_{n})+\dim\ker\Big(\frac{\kappa_{k}}{\mu_{k}+\eta_{k}}\sum_{i=1}^{n}\beta_{i}\textbf{w}_{0}\otimes\textbf{e}_{i}\Big)\nonumber\\
 &\leq&\dim\ker\small\left[\begin{array}{cc}
 L_{k}\otimes I_{n} & \frac{\kappa_{k}}{\mu_{k}+\eta_{k}}\sum_{i=1}^{n}\beta_{i}\textbf{w}_{0}\otimes\textbf{e}_{i}
 \end{array}\right],
 \end{eqnarray*} which implies that  ${\mathrm{rank}}(L_{k}\otimes I_{n})\geq{\mathrm{rank}}\small\left[\begin{array}{cc}
 L_{k}\otimes I_{n} & \frac{\kappa_{k}}{\mu_{k}+\eta_{k}}\sum_{i=1}^{n}\beta_{i}\textbf{w}_{0}\otimes\textbf{e}_{i}
 \end{array}\right]$. Hence, (\ref{rankcon}) holds.
Next, it follows from $vi$) of Proposition 6.1.7 of \cite[p.~400]{Bernstein:2009} and $viii$) of Proposition 6.1.6 of \cite[p.~399]{Bernstein:2009} that the general solution to (\ref{Axbeqn}) is given by the form
\begin{eqnarray}\label{usolution} 
[\textbf{u}_{1}^{*},\ldots,\textbf{u}_{q}^{*}]^{*}&=&\frac{\kappa_{k}}{\mu_{k}+\eta_{k}}\sum_{i=1}^{n}\beta_{i}(L_{k}\otimes I_{n})^{+}(\textbf{w}_{0}\otimes\textbf{e}_{i})+\sum_{l=1}^{q}\sum_{i=1}^{n}\gamma_{li}(I_{nq}-(L_{k}\otimes I_{n})^{+}(L_{k}\otimes I_{n}))(\textbf{g}_{l}\otimes\textbf{e}_{i})\nonumber\\
&=&\frac{\kappa_{k}}{\mu_{k}+\eta_{k}}\sum_{i=1}^{n}\beta_{i}(L_{k}^{+}\otimes I_{n})(\textbf{w}_{0}\otimes\textbf{e}_{i})+\sum_{l=1}^{q}\sum_{i=1}^{n}\gamma_{li}(I_{nq}-(L_{k}^{+}\otimes I_{n})(L_{k}\otimes I_{n}))(\textbf{g}_{l}\otimes\textbf{e}_{i})\nonumber\\
&=&\frac{\kappa_{k}}{\mu_{k}+\eta_{k}}\sum_{i=1}^{n}\beta_{i}L_{k}^{+}\textbf{w}_{0}\otimes\textbf{e}_{i}+\sum_{l=1}^{q}\sum_{i=1}^{n}\gamma_{li}(I_{q}\otimes I_{n}-(L_{k}^{+}L_{k}\otimes I_{n}))(\textbf{g}_{l}\otimes\textbf{e}_{i})\nonumber\\
&=&\frac{\kappa_{k}}{\mu_{k}+\eta_{k}}\sum_{i=1}^{n}\beta_{i}L_{k}^{+}\textbf{w}_{0}\otimes\textbf{e}_{i}+\sum_{l=1}^{q}\sum_{i=1}^{n}\gamma_{li}((I_{q}-L_{k}^{+}L_{k})\otimes I_{n}))(\textbf{g}_{l}\otimes\textbf{e}_{i})\nonumber\\
&=&\frac{\kappa_{k}}{\mu_{k}+\eta_{k}}\sum_{i=1}^{n}\beta_{i}L_{k}^{+}\textbf{w}_{0}\otimes\textbf{e}_{i}+\sum_{l=1}^{q}\sum_{i=1}^{n}\gamma_{li}(\textbf{g}_{l}-L_{k}^{+}L_{k}\textbf{g}_{l})\otimes \textbf{e}_{i},
\end{eqnarray} where $\gamma_{li}\in\mathbb{C}$. Hence, the general solution to (\ref{Axbeqn}) is given by the form
\begin{eqnarray}\label{gsolution}
[\textbf{u}_{1}^{*},\ldots,\textbf{u}_{q}^{*}]^{*}&=&\frac{\kappa_{k}}{\mu_{k}+\eta_{k}}\sum_{i=1}^{n}\beta_{i}L_{k}^{+}\textbf{w}_{0}\otimes\textbf{e}_{i}+\sum_{l=1}^{q}\sum_{i=1}^{n}\gamma_{li}(\textbf{g}_{l}-L_{k}^{+}L_{k}\textbf{g}_{l})\otimes \textbf{e}_{i}.
\end{eqnarray}
Note that it follows from (\ref{Enqu}) that $E_{n\times nq}^{[j]}[\textbf{u}_{1}^{*},\ldots,\textbf{u}_{q}^{*}]^{*}=\textbf{u}_{j}=\textbf{0}_{n\times 1}$. Then both the general solution (\ref{gsolution}) should satisfy this constraint. It now follows from (\ref{gsolution}) that $\beta_{i}\in\mathbb{C}$ and $\gamma_{li}\in\mathbb{C}$ in (\ref{gsolution}) should satisfy
\begin{eqnarray}
&&\frac{\kappa_{k}}{\mu_{k}+\eta_{k}}\sum_{i=1}^{n}\beta_{i}E_{n\times nq}^{[j]}(L_{k}^{+}\textbf{w}_{0}\otimes\textbf{e}_{i})+\sum_{l=1}^{q}\sum_{i=1}^{n}\gamma_{li}E_{n\times nq}^{[j]}((\textbf{g}_{l}-L_{k}^{+}L_{k}\textbf{g}_{l})\otimes \textbf{e}_{i})=\textbf{0}_{n\times 1}.\label{constraint}
\end{eqnarray} Note that $E_{n\times nq}^{[j]}=\textbf{g}_{j}^{\mathrm{T}}\otimes I_{n}$ for every $j=1,\ldots,q$. Then for every $i=1,\ldots,n$ and $j=1,\ldots,q$,
\begin{eqnarray*}
E_{n\times nq}^{[j]}(L_{k}^{+}\textbf{w}_{0}\otimes\textbf{e}_{i})=(\textbf{g}_{j}^{\mathrm{T}}\otimes I_{n})(L_{k}^{+}\textbf{w}_{0}\otimes\textbf{e}_{i})=\textbf{g}_{j}^{\mathrm{T}}L_{k}^{+}\textbf{w}_{0}\otimes \textbf{e}_{i}=(\textbf{g}_{j}^{\mathrm{T}}L_{k}^{+}\textbf{w}_{0})\textbf{e}_{i}.
\end{eqnarray*} Similarly, one can obtain that $E_{n\times nq}^{[j]}((\textbf{g}_{l}-L_{k}^{+}L_{k}\textbf{g}_{l})\otimes \textbf{e}_{i})=(\textbf{g}_{j}^{\mathrm{T}}\textbf{g}_{l}-\textbf{g}_{j}^{\mathrm{T}}L_{k}^{+}L_{k}\textbf{g}_{l})\textbf{e}_{i}$ for every $i=1,\ldots,n$, every $j=1,\ldots,q$, and every $l=1,\ldots,q$. Now using these relationships, (\ref{constraint}) can be simplified as 
\begin{eqnarray*}
\textbf{0}_{n\times 1}&=&\sum_{i=1}^{n}\frac{\kappa_{k}}{\mu_{k}+\eta_{k}}\beta_{i}(\textbf{g}_{j}^{\mathrm{T}}L_{k}^{+}\textbf{w}_{0})\textbf{e}_{i}+\sum_{i=1}^{n}\sum_{l=1}^{q}\gamma_{li}(\textbf{g}_{j}^{\mathrm{T}}\textbf{g}_{l}-\textbf{g}_{j}^{\mathrm{T}}L_{k}^{+}L_{k}\textbf{g}_{l})\textbf{e}_{i}\nonumber\\
&=&\sum_{i=1}^{n}\Big[\frac{\kappa_{k}}{\mu_{k}+\eta_{k}}\beta_{i}(\textbf{g}_{j}^{\mathrm{T}}L_{k}^{+}\textbf{w}_{0})+\sum_{l=1}^{q}\gamma_{li}(\textbf{g}_{j}^{\mathrm{T}}\textbf{g}_{l}-\textbf{g}_{j}^{\mathrm{T}}L_{k}^{+}L_{k}\textbf{g}_{l})\Big]\textbf{e}_{i},
\end{eqnarray*} which imply that $\beta_{i}\in\mathbb{C}$ and $\gamma_{li}\in\mathbb{C}$ in (\ref{gsolution}) satisfy
\begin{eqnarray}
\frac{\kappa_{k}}{\mu_{k}+\eta_{k}}\beta_{i}(\textbf{g}_{j}^{\mathrm{T}}L_{k}^{+}\textbf{w}_{0})+\sum_{l=1}^{q}\gamma_{li}(\textbf{g}_{j}^{\mathrm{T}}\textbf{g}_{l}-\textbf{g}_{j}^{\mathrm{T}}L_{k}^{+}L_{k}\textbf{g}_{l})=0,\label{consr1}
\end{eqnarray} for every $i=1,\ldots,n$ and every $j=1,\ldots,q$. Finally, since (\ref{Axbeqnab}) has infinitely many solutions due to (\ref{rankconab}), it follows that there exist infinitely many $\beta_{i}\in\mathbb{C}$ and $\gamma_{li}\in\mathbb{C}$ satisfying (\ref{consr1}).

In summary, if $1\in{\mathrm{spec}}(\frac{\mu_{k}+\eta_{k}}{\kappa_{k}}L_{k})$ and $h_{k}=1+\frac{1}{\kappa_{k}}$, then $\lambda=-\kappa_{k}$ is indeed an eigenvalue of $A_{k}^{[j]}+h_{k}A_{{\mathrm{c}}k}$. In this case, $\textbf{x}_{1}=\textbf{0}_{nq\times 1}$, $\textbf{x}_{2}=[\textbf{u}_{1}^{*},\ldots,\textbf{u}_{q}^{*}]^{*}$ given by (\ref{gsolution}) with $\textbf{u}_{j}=\textbf{0}_{n\times 1}$, and $\textbf{x}_{3}=\sum_{i=1}^{n}\beta_{i}\textbf{e}_{i}$, where not all of $\beta_{i}$ and $\gamma_{li}$ are zero. The corresponding eigenvectors for $\lambda_{3}$ are given by
\begin{eqnarray}\label{xw0}
\textbf{x}&=&\Big[\textbf{0}_{1\times nq},\frac{\kappa_{k}}{\mu_{k}+\eta_{k}}\sum_{i=1}^{n}\beta_{i}(L_{k}^{+}\textbf{w}_{0}\otimes\textbf{e}_{i})^{\mathrm{T}}+\sum_{l=1}^{q}\sum_{i=1}^{n}\gamma_{li}((\textbf{g}_{l}-L_{k}^{+}L_{k}\textbf{g}_{l})\otimes \textbf{e}_{i})^{\mathrm{T}},\sum_{i=1}^{n}\beta_{i}\textbf{e}_{i}^{\mathrm{T}}\Big]^{*},
\end{eqnarray} where $\beta_{i}\in\mathbb{C}$ and $\gamma_{li}\in\mathbb{C}$ satisfy (\ref{consr1}) and not all of them are zero. Therefore,
\begin{eqnarray}
&&\ker\Big(A_{k}^{[j]}+h_{k}A_{{\mathrm{c}}k}-\lambda_{3} I_{2nq+n}\Big)\nonumber\\
&&=\Big\{\Big[\textbf{0}_{1\times nq},\frac{\kappa_{k}}{\mu_{k}+\eta_{k}}\sum_{i=1}^{n}\beta_{i}(L_{k}^{+}\textbf{w}_{0}\otimes\textbf{e}_{i})^{\mathrm{T}}+\sum_{l=1}^{q}\sum_{i=1}^{n}\gamma_{li}((\textbf{g}_{l}-L_{k}^{+}L_{k}\textbf{g}_{l})\otimes \textbf{e}_{i})^{\mathrm{T}}),\sum_{i=1}^{n}\beta_{i}\textbf{e}_{i}^{\mathrm{T}}\Big]^{*}:(\ref{consr1})\,\,{\mathrm{holds}},\nonumber\\
&&\beta_{i}\in\mathbb{C},\gamma_{li}\in\mathbb{C},i=1,\ldots,n,l=1,\ldots,q\Big\}.
\end{eqnarray} This expression of $\ker\Big(A_{k}^{[j]}+h_{k}A_{{\mathrm{c}}k}-\lambda_{3} I_{2nq+n}\Big)$ is slightly different from the one in  (\ref{egns8}) since it involves the constraint (\ref{consr1}) for $\beta_{i}\in\mathbb{C}$ and $\gamma_{li}\in\mathbb{C}$. Nevertheless, they are equivalent to each other since both expressions are the general solution to the same form of linear equations. \hfill$\blacklozenge$
\end{remark}

\begin{remark}\label{rmk_q}
If ${\mathrm{rank}}(L_{k})=q-1$, then it follows that $\ker(L_{k})={\mathrm{span}}\{\textbf{w}_{0}\}$. In graph theory, this rank condition is implied by the strong connectivity condition on $\mathcal{G}_{k}$. In this case, if $\eta_{k}\lambda_{5,6}+\mu_{k} h_{k}\lambda_{5,6}+\mu_{k}\neq0$, where $\lambda_{5,6}$ are given by (\ref{lambda5}), then $\lambda_{5,6}$ are not the eigenvalues of $A_{k}^{[j]}+h_{k}A_{{\mathrm{c}}k}$ and $\{0\}\subseteq{\mathrm{spec}}(A_{k}^{[j]}+h_{k}A_{{\mathrm{c}}k})\subseteq\{0,-\kappa_{k},-\frac{\kappa_{k}(1+h_{k})}{2}\pm\frac{1}{2}\sqrt{\kappa^{2}(1+h_{k})^{2}-4\kappa_{k}},\lambda\in\mathbb{C}:\forall \frac{\lambda^{2}+\kappa_{k} h_{k}\lambda+\kappa_{k}}{\eta_{k}\lambda+\mu_{k} h_{k}\lambda+\mu_{k}}\in{\mathrm{spec}}(-L_{k})\backslash\{0\}\}$. This is because $\ker(L_{k}\otimes I_{n})={\mathrm{span}}\{\textbf{w}_{0}\otimes\textbf{e}_{1},\ldots,\textbf{w}_{0}\otimes\textbf{e}_{n}\}$, (\ref{gs0_1}) and (\ref{gs0_2}) only have the trivial solution $\textbf{v}=\textbf{0}_{nq\times 1}$, which contradicts the definition of eigenvectors for $\lambda=\lambda_{5,6}$. Hence, $\lambda\neq\lambda_{5,6}$. Furthermore, note that it follows from Lemma~\ref{lemma_argument} that if $\mathcal{G}_{k}$ is undirected, then $\lambda\in\mathbb{C}$ in ${\mathrm{spec}}(A_{k}^{[j]}+h_{k}A_{{\mathrm{c}}k})$ is such that $\frac{\lambda^{2}+\kappa_{k} h_{k}\lambda+\kappa_{k}}{\eta_{k}\lambda+\mu_{k} h_{k}\lambda+\mu_{k}}<0$. \hfill$\blacklozenge$
\end{remark}

\begin{lemma}\label{lemma_B}
Define a (possibly infinite) series of matrices $B^{[j]}_{k}$, $j=1,\ldots,q$, $k=0,1,2,\ldots$, as follows:
\begin{eqnarray}\label{Bmatrix}
B_{k}^{[j]}=\small\left[\begin{array}{ccc}
\textbf{0}_{nq\times nq} & h_{k}I_{nq} & \textbf{0}_{nq\times n} \\
-h_{k}\mu_{k} L_{k}\otimes I_{n}-h_{k}\kappa I_{nq} & -h_{k}\eta_{k} L_{k}\otimes I_{n} & h_{k}\kappa_{k} \textbf{1}_{q\times 1}\otimes I_{n} \\
E_{n\times nq}^{[j]} & \textbf{0}_{n\times nq} & -I_{n} \\
\end{array}\right],
\end{eqnarray} where $\mu_{k},\eta_{k},\kappa_{k}\geq0$, $h_{k}>0$, $k\in\overline{\mathbb{Z}}_{+}$, $L_{k}\in\mathbb{R}^{q\times q}$ denotes the Laplacian matrix of a node-fixed dynamic digraph $\mathcal{G}_{k}$, and $E_{n\times nq}^{[j]}\in\mathbb{R}^{n\times nq}$ is defined in Lemma~\ref{lemma_EW}. Then for every $j=1,\ldots,q$, $\{0\}\subseteq{\mathrm{spec}}(B_{k}^{[j]}+h_{k}^{2}A_{{\mathrm{c}}k})\subseteq\{0,-1,-\frac{h_{k}^{2}\kappa_{k}}{2}\pm\frac{1}{2}\sqrt{(h_{k}^{2}\kappa_{k})^{2}-4h_{k}^{2}\kappa_{k}},\lambda_{1},\lambda_{2}\in\mathbb{C}:\forall\frac{\lambda_{1}^{2}+\kappa_{k} h_{k}^{2}\lambda_{1}+\kappa_{k}h_{k}^{2}}{\eta_{k}h_{k}\lambda_{1}+\mu_{k} h_{k}^{2}\lambda_{1}+\mu_{k}h_{k}^{2}}\in{\mathrm{spec}}(-L_{k})\backslash\{0\},\lambda_{2}^{3}+(1+h_{k}^{2}\kappa_{k})\lambda_{2}^{2}+(2h_{k}^{2}\kappa_{k}-h_{k}\kappa_{k})\lambda_{2}+h_{k}^{2}\kappa_{k}=0\}$, where $A_{{\mathrm{c}}k}$ is defined by (\ref{Ac}) in Lemma~\ref{lemma_semisimple}. Furthermore, if $h_{k}\kappa_{k}\neq0$, then 0 is semisimple.  
\end{lemma}

\begin{IEEEproof}
For a fixed $j\in\{1,\ldots,q\}$, let $\lambda\in{\mathrm{spec}}(B_{k}^{[j]}+h_{k}^{2}A_{{\mathrm{c}}k})$ and $\textbf{x}=[\textbf{x}_{1}^{*},\textbf{x}_{2}^{*},\textbf{x}_{3}^{*}]^{*}\in\mathbb{C}^{2nq+n}$ be the corresponding eigenvector for $\lambda$, where $\textbf{x}_{1},\textbf{x}_{2}\in\mathbb{C}^{nq}$ and $\textbf{x}_{3}\in\mathbb{C}^{n}$. Then it follows from
$(B_{k}^{[j]}+h_{k}^{2}A_{{\mathrm{c}}k})\textbf{x}=\lambda\textbf{x}$ that 
\begin{eqnarray}
h_{k}\textbf{x}_{2}+h_{k}^{2}[-\mu_{k} (L_{k}\otimes I_{n})\textbf{x}_{1}-\kappa_{k}\textbf{x}_{1}-\eta_{k}(L_{k}\otimes I_{n})\textbf{x}_{2}+\kappa_{k}(\textbf{1}_{q\times 1}\otimes I_{n})\textbf{x}_{3}]=\lambda\textbf{x}_{1},\label{hx11}\\
h_{k}[-\mu_{k} (L_{k}\otimes I_{n})\textbf{x}_{1}-\kappa_{k}\textbf{x}_{1}-\eta_{k}(L_{k}\otimes I_{n})\textbf{x}_{2}+\kappa_{k}(\textbf{1}_{q\times 1}\otimes I_{n})\textbf{x}_{3}]=\lambda \textbf{x}_{2},\label{x21}\\
E_{n\times nq}^{[j]}\textbf{x}_{1}-\textbf{x}_{3}=\lambda \textbf{x}_{3}.\label{Aeig_31}
\end{eqnarray} Let $\textbf{x}_{3}\neq\textbf{0}_{n\times1}$ be arbitrary, $\textbf{x}_{1}=(\textbf{1}_{q\times 1}\otimes I_{n})\textbf{x}_{3}$, and $\textbf{x}_{2}=\textbf{0}_{nq\times 1}$. Clearly such $\textbf{x}_{i}$, $i=1,2,3$, satisfy (\ref{hx11})--(\ref{Aeig_31}) with $\lambda=0$. Hence, $\lambda=0$ is always an eigenvalue of $B_{k}^{[j]}+h_{k}^{2}A_{{\mathrm{c}}k}$. Next, we assume that $\lambda\neq0$.

Substituting (\ref{x21}) into (\ref{hx11}) yields
$\textbf{x}_{1}=\frac{h_{k}(1+\lambda)}{\lambda}\textbf{x}_{2}$. Replacing $\textbf{x}_{1}$ in (\ref{x21}) and (\ref{Aeig_31}) with $\textbf{x}_{1}=\frac{h_{k}(1+\lambda)}{\lambda}\textbf{x}_{2}$ yields
\begin{eqnarray}
\Big[\Big(\frac{h_{k}^{2}\mu_{k}}{\lambda}+\mu_{k} h_{k}^{2}+\eta_{k}h_{k}\Big)(L_{k}\otimes I_{n})+\Big(\frac{h_{k}^{2}\kappa_{k}}{\lambda}+\lambda+h_{k}^{2}\kappa_{k}\Big)I_{nq}\Big]\textbf{x}_{2}-h_{k}\kappa_{k}(\textbf{1}_{q\times 1}\otimes I_{n})\textbf{x}_{3}=\textbf{0}_{nq\times 1},\label{hx21}\\
E_{n\times nq}^{[j]}\textbf{x}_{2}-(1+\lambda)\textbf{x}_{3}=\textbf{0}_{n\times 1}.\label{hx31}
\end{eqnarray} Thus, (\ref{hx21}) and (\ref{hx31}) have nontrivial solutions if and only if
\begin{eqnarray}\label{detcon1}
\det\small\left[\begin{array}{cc}
\Big(\frac{h_{k}^{2}\mu_{k}}{\lambda}+\mu_{k} h_{k}^{2}+\eta_{k}h_{k}\Big)(L_{k}\otimes I_{n})+\Big(\frac{h_{k}^{2}\kappa_{k}}{\lambda}+\lambda+h_{k}^{2}\kappa_{k}\Big)I_{nq} & -h_{k}\kappa_{k}(\textbf{1}_{q\times 1}\otimes I_{n})\\
E_{n\times nq}^{[j]} & -(1+\lambda)I_{n}
\end{array}\right]=0.
\end{eqnarray} 

If $\det\Big[\Big(\frac{h_{k}^{2}\mu_{k}}{\lambda}+\mu_{k} h_{k}^{2}+\eta_{k}h_{k}\Big)(L_{k}\otimes I_{n})+\Big(\frac{h_{k}^{2}\kappa_{k}}{\lambda}+\lambda+h_{k}^{2}\kappa_{k}\Big)I_{nq}\Big]\neq0$, then pre-multiplying $-L_{k}\otimes I_{n}$ on both sides of (\ref{hx21}) and following the similar arguments as in the proof of $ii$) of Lemma~\ref{lemma_A}, we have $\textbf{x}_{2}=\sum_{l=0}^{q-1-{\mathrm{rank}}(L_{k})}\sum_{i=1}^{n}\varpi_{li}(\textbf{w}_{l}\otimes\textbf{e}_{i})$, where $\varpi_{li}\in\mathbb{C}$. Substituting this expression of $\textbf{x}_{2}$ into (\ref{hx21}) and (\ref{hx31}) by using $iii$) of Lemma~\ref{lemma_EW} yields
\begin{eqnarray}
\Big(\frac{h_{k}^{2}\kappa_{k}}{\lambda}+\lambda+h_{k}^{2}\kappa_{k}\Big)\sum_{l=0}^{q-1-{\mathrm{rank}}(L_{k})}\sum_{i=1}^{n}\varpi_{li}w_{lj}\textbf{e}_{i}-h_{k}\kappa_{k}\textbf{x}_{3}=\textbf{0}_{n\times 1},\label{x3_eqn1}\\
\sum_{l=0}^{q-1-{\mathrm{rank}}(L_{k})}\sum_{i=1}^{n}\varpi_{li}w_{lj}\textbf{e}_{i}-(1+\lambda)\textbf{x}_{3}=\textbf{0}_{n\times 1}.\label{x3_eqn2}
\end{eqnarray} Substituting (\ref{x3_eqn2}) into (\ref{x3_eqn1}) yields
\begin{eqnarray}
\Big[\Big(\frac{h_{k}^{2}\kappa_{k}}{\lambda}+\lambda+h_{k}^{2}\kappa_{k}\Big)(1+\lambda)-h_{k}\kappa_{k}\Big]\textbf{x}_{3}=\textbf{0}_{n\times 1}.
\end{eqnarray} If $\textbf{x}_{3}=\textbf{0}_{n\times 1}$, then it follows from (\ref{hx21}) that $\textbf{x}_{2}=\textbf{0}_{nq\times 1}$, and hence, $\textbf{x}_{1}=\textbf{0}_{nq\times 1}$, which is a contradiction since $\textbf{x}$ is an eigenvector. Thus, $\textbf{x}_{3}\neq\textbf{0}_{n\times 1}$ and consequently,
$\Big(\frac{h_{k}^{2}\kappa_{k}}{\lambda}+\lambda+h_{k}^{2}\kappa_{k}\Big)(1+\lambda)-h_{k}\kappa_{k}=0$,
i.e.,
\begin{eqnarray}\label{cubic3}
\lambda^{3}+(1+h_{k}^{2}\kappa_{k})\lambda^{2}+(2h_{k}^{2}\kappa_{k}-h_{k}\kappa_{k})\lambda+h_{k}^{2}\kappa_{k}=0.
\end{eqnarray}
Solving this cubic equation in terms of $\lambda$ gives the possible eigenvalues of $B_{k}^{[j]}+h_{k}^{2}A_{{\mathrm{c}}k}$. This can be done via Cardano's formula. If $h_{k}\kappa_{k}=0$, then $\lambda=-1$. Otherwise, if $h_{k}\kappa_{k}\neq0$, then it follows from Routh's Stability Criterion that ${\mathrm{Re}}\,\lambda<0$ if and only if $2h_{k}^{2}\kappa_{k}-h_{k}\kappa_{k}>0$ and $(1+h_{k}^{2}\kappa_{k})(2h_{k}^{2}\kappa_{k}-h_{k}\kappa_{k})>h_{k}^{2}\kappa_{k}$, that is, $h_{k}>1/2$ and $h_{k}+2h_{k}^{3}\kappa_{k}>1+h_{k}^{2}\kappa_{k}$. 

Alternatively, if $\det\Big[\Big(\frac{h_{k}^{2}\mu_{k}}{\lambda}+\mu_{k} h_{k}^{2}+\eta_{k}h_{k}\Big)(L_{k}\otimes I_{n})+\Big(\frac{h_{k}^{2}\kappa_{k}}{\lambda}+\lambda+h_{k}^{2}\kappa_{k}\Big)I_{nq}\Big]=0$, then in this case, (\ref{detcon1}) holds if $\lambda=-1$, or $\lambda\neq-1$ and by Proposition 2.8.4 of \cite[p.~116]{Bernstein:2009}, $\det\Big(\Big(\frac{h_{k}^{2}\mu_{k}}{\lambda}+\mu_{k} h_{k}^{2}+\eta_{k}h_{k}\Big)(L_{k}\otimes I_{n})+\Big(\frac{h_{k}^{2}\kappa_{k}}{\lambda}+\lambda+h_{k}^{2}\kappa_{k}\Big)I_{nq}-\frac{\kappa_{k}h_{k}}{1+\lambda}W^{[j]}\Big)=0$, which implies that for $\lambda\neq-1$, the equation
\begin{eqnarray}\label{eqn_v1}
\Big(\Big(\frac{h_{k}^{2}\mu_{k}}{\lambda}+\mu_{k} h_{k}^{2}+\eta_{k}h_{k}\Big)(L_{k}\otimes I_{n})+\Big(\frac{h_{k}^{2}\kappa_{k}}{\lambda}+\lambda+h_{k}^{2}\kappa_{k}\Big)I_{nq}-\frac{\kappa_{k}h_{k}}{1+\lambda}W^{[j]}\Big)\textbf{v}=\textbf{0}_{nq\times 1}
\end{eqnarray} has nontrivial solutions for $\textbf{v}\in\mathbb{C}^{nq}$. Again, note that for every $j=1,\ldots,q$, $(L_{k}\otimes I_{n})W^{[j]}=\textbf{0}_{nq\times nq}$. Pre-multiplying $L_{k}\otimes I_{n}$ on both sides of (\ref{eqn_v1}) yields $\Big(\Big(\frac{h_{k}^{2}\mu_{k}}{\lambda}+\mu_{k} h_{k}^{2}+\eta_{k}h_{k}\Big)(L_{k}\otimes I_{n})^{2}+\Big(\frac{h_{k}^{2}\kappa_{k}}{\lambda}+\lambda+h_{k}^{2}\kappa_{k}\Big)(L_{k}\otimes I_{n})\Big)\textbf{v}=(L_{k}\otimes I_{n})\Big(\Big(\frac{h_{k}^{2}\mu_{k}}{\lambda}+\mu_{k} h_{k}^{2}+\eta_{k}h_{k}\Big)(L_{k}\otimes I_{n})+\Big(\frac{h_{k}^{2}\kappa_{k}}{\lambda}+\lambda+h_{k}^{2}\kappa_{k}\Big)I_{nq}\Big)\textbf{v}=\textbf{0}_{nq\times 1}$, which implies that $\Big(\Big(\frac{h_{k}^{2}\mu_{k}}{\lambda}+\mu_{k} h_{k}^{2}+\eta_{k}h_{k}\Big)(L_{k}\otimes I_{n})+\Big(\frac{h_{k}^{2}\kappa_{k}}{\lambda}+\lambda+h_{k}^{2}\kappa_{k}\Big)I_{nq}\Big)\textbf{v}\in\ker(L_{k}\otimes I_{n})$. Since $\ker(L_{k}\otimes I_{n})=\bigcup_{l=0}^{q-1-{\mathrm{rank}}(L_{k})}{\mathrm{span}}\{\textbf{w}_{l}\otimes\textbf{e}_{1},\ldots,\textbf{w}_{l}\otimes\textbf{e}_{n}\}$, it follows that
\begin{eqnarray}\label{Az=b1}
\Big(\Big(\frac{h_{k}^{2}\mu_{k}}{\lambda}+\mu_{k} h_{k}^{2}+\eta_{k}h_{k}\Big)(L_{k}\otimes I_{n})+\Big(\frac{h_{k}^{2}\kappa_{k}}{\lambda}+\lambda+h_{k}^{2}\kappa_{k}\Big)I_{nq}\Big)\textbf{v}=\sum_{i=1}^{n}\sum_{l=0}^{q-1-{\mathrm{rank}}(L_{k})}\omega_{li}\textbf{w}_{l}\otimes\textbf{e}_{i},
\end{eqnarray} where $\omega_{li}\in\mathbb{C}$, which is similar to (\ref{Az=b}). Now it follows from (\ref{eqn_v1}) and (\ref{Az=b1}) that
\begin{eqnarray}\label{Wv_eqn1}
\frac{\kappa_{k}h_{k}}{1+\lambda}W^{[j]}\textbf{v}=\sum_{i=1}^{n}\sum_{l=0}^{q-1-{\mathrm{rank}}(L_{k})}\omega_{li}\textbf{w}_{l}\otimes\textbf{e}_{i}.
\end{eqnarray}

If $\frac{h_{k}^{2}\kappa_{k}}{\lambda}+\lambda+h_{k}^{2}\kappa_{k}\neq0$, then it follows from the similar arguments after (\ref{Wv_eqn}) that $\omega_{\ell i}=0$ for every $i=1,\ldots,n$ and every $\ell=1,\ldots,q-1-{\mathrm{rank}}(L_{k})$. Furthermore, 
\begin{eqnarray}
\omega_{0i}-\frac{\lambda\kappa_{k}h_{k}}{(1+\lambda)(\lambda^{2}+h_{k}^{2}\kappa_{k}\lambda+h_{k}^{2}\kappa_{k})}\omega_{0i}=0,\quad i=1,\ldots,n.
\end{eqnarray}
Then either $1-\frac{\lambda\kappa_{k}h_{k}}{(1+\lambda)(\lambda^{2}+h_{k}^{2}\kappa_{k}\lambda+h_{k}^{2}\kappa_{k})}=0$ or $\omega_{0i}=0$ for every $i=1,\ldots,n$.
If $\frac{\lambda\kappa_{k}h_{k}}{(1+\lambda)(\lambda^{2}+h_{k}^{2}\kappa_{k}\lambda+h_{k}^{2}\kappa_{k})}=1$, then
\begin{eqnarray} 
\lambda^{3}+(1+h_{k}^{2}\kappa_{k})\lambda^{2}+(2h_{k}^{2}\kappa_{k}-h_{k}\kappa_{k})\lambda+h_{k}^{2}\kappa_{k}=0,
\end{eqnarray} which is the same as (\ref{cubic3}). Since $\lambda\neq -1$, in this case $\kappa_{k}h_{k}\neq0$. Then it follows from Routh's Stability Criterion that ${\mathrm{Re}}\,\lambda<0$ if and only if $h_{k}>1/2$ and $h_{k}+2h_{k}^{3}\kappa_{k}>1+h_{k}^{2}\kappa_{k}$.
If $\omega_{0i}=0$ for every $i=1,\ldots,n$, then it follows from (\ref{eqn_v1}) and (\ref{Az=b1}) that $\frac{\kappa_{k}h_{k}}{1+\lambda}W^{[j]}\textbf{v}=\textbf{0}_{nq\times 1}$ and $\Big(\Big(\frac{h_{k}^{2}\mu_{k}}{\lambda}+\mu_{k} h_{k}^{2}+\eta_{k}h_{k}\Big)(L_{k}\otimes I_{n})+\Big(\frac{h_{k}^{2}\kappa_{k}}{\lambda}+\lambda+h_{k}^{2}\kappa_{k}\Big)I_{nq}\Big)\textbf{v}=\textbf{0}_{nq\times 1}$, which implies that $\textbf{v}\in\ker\Big(\Big(\frac{h_{k}^{2}\mu_{k}}{\lambda}+\mu_{k} h_{k}^{2}+\eta_{k}h_{k}\Big)(L_{k}\otimes I_{n})+\Big(\frac{h_{k}^{2}\kappa_{k}}{\lambda}+\lambda+h_{k}^{2}\kappa_{k}\Big)I_{nq}\Big)\cap\ker(\frac{\kappa_{k}h_{k}}{1+\lambda}W^{[j]})$. Clearly $\frac{h_{k}^{2}\mu_{k}}{\lambda}+\mu_{k} h_{k}^{2}+\eta_{k}h_{k}\neq0$. In this case, $\lambda\in\{\lambda_{1}\in\mathbb{C}:\forall\frac{\lambda_{1}^{2}+\kappa_{k} h_{k}^{2}\lambda_{1}+\kappa_{k}h_{k}^{2}}{\eta_{k}h_{k}\lambda_{1}+\mu_{k} h_{k}^{2}\lambda_{1}+\mu_{k}h_{k}^{2}}\in{\mathrm{spec}}(-L_{k})\backslash\{0\}\}$.

Alternatively, if $\frac{h_{k}^{2}\kappa_{k}}{\lambda}+\lambda+h_{k}^{2}\kappa_{k}=0$, then it follows from the similar arguments after (\ref{eigv4}) in Lemma~\ref{lemma_A} that
\begin{eqnarray}
\lambda=-\frac{h_{k}^{2}\kappa_{k}}{2}\pm\frac{1}{2}\sqrt{(h_{k}^{2}\kappa_{k})^{2}-4h_{k}^{2}\kappa_{k}}
\end{eqnarray} are the possible eigenvalues of $B_{k}^{[j]}+h_{k}^{2}A_{{\mathrm{c}}k}$. 

In summary,
\begin{eqnarray}\label{egspace} 
&&\{0\}\subseteq{\mathrm{spec}}(B_{k}^{[j]}+h_{k}^{2}A_{{\mathrm{c}}k})\subseteq\nonumber\\
&&\Big\{0,-1,-\frac{h_{k}^{2}\kappa_{k}}{2}\pm\frac{1}{2}\sqrt{(h_{k}^{2}\kappa_{k})^{2}-4h_{k}^{2}\kappa_{k}},\lambda_{1},\lambda_{2}\in\mathbb{C}:\forall\frac{\lambda_{1}^{2}+\kappa_{k} h_{k}^{2}\lambda_{1}+\kappa_{k}h_{k}^{2}}{\eta_{k}h_{k}\lambda_{1}+\mu_{k} h_{k}^{2}\lambda_{1}+\mu_{k}h_{k}^{2}}\in{\mathrm{spec}}(-L_{k})\backslash\{0\},\nonumber\\
&&\lambda_{2}^{3}+(1+h_{k}^{2}\kappa_{k})\lambda_{2}^{2}+(2h_{k}^{2}\kappa_{k}-h_{k}\kappa_{k})\lambda_{2}+h_{k}^{2}\kappa_{k}=0\Big\}. 
\end{eqnarray} Finally, the semisimplicity property of 0 can be proved by using the similar arguments as in the proof of Lemma~\ref{lemma_semisimple}. 
\end{IEEEproof}

\begin{remark}
Similar to Remark~\ref{rmk_q}, if ${\mathrm{rank}}(L_{k})=q-1$, then $-\frac{h_{k}^{2}\kappa_{k}}{2}\pm\frac{1}{2}\sqrt{(h_{k}^{2}\kappa_{k})^{2}-4h_{k}^{2}\kappa_{k}}$ are not the eigenvalues of $B_{k}^{[j]}+h_{k}^{2}A_{{\mathrm{c}}k}$ and $\{0\}\subseteq{\mathrm{spec}}(B_{k}^{[j]}+h_{k}^{2}A_{{\mathrm{c}}k})\subseteq\{0,-1,\lambda_{1},\lambda_{2}\in\mathbb{C}:\forall\frac{\lambda_{1}^{2}+\kappa_{k} h_{k}^{2}\lambda_{1}+\kappa_{k}h_{k}^{2}}{\eta_{k}h_{k}\lambda_{1}+\mu_{k} h_{k}^{2}\lambda_{1}+\mu_{k}h_{k}^{2}}\in{\mathrm{spec}}(-L_{k})\backslash\{0\},\lambda_{2}^{3}+(1+h_{k}^{2}\kappa_{k})\lambda_{2}^{2}+(2h_{k}^{2}\kappa_{k}-h_{k}\kappa_{k})\lambda_{2}+h_{k}^{2}\kappa_{k}=0\}$. Furthermore, note that it follows from Lemma~\ref{lemma_argument} that if $\mathcal{G}_{k}$ is undirected, then $\lambda_{1}\in\mathbb{C}$ in ${\mathrm{spec}}(B_{k}^{[j]}+h_{k}^{2}A_{{\mathrm{c}}k})$ is such that $\frac{\lambda_{1}^{2}+\kappa_{k} h_{k}^{2}\lambda_{1}+\kappa_{k}h_{k}^{2}}{\eta_{k}h_{k}\lambda_{1}+\mu_{k} h_{k}^{2}\lambda_{1}+\mu_{k}h_{k}^{2}}<0$.
Finally, one can also discuss the detailed eigenspace for each possible eigenvalue in (\ref{egspace}) by using the similar arguments in Lemmas \ref{lemma_Arank}--\ref{lemma_A}. \hfill$\blacklozenge$
\end{remark}

The following definition is due to \cite{HH:IJC:2009}.

\begin{definition}\label{def_so}
Let $A\in\mathbb{R}^{n\times n}$ and $C\in\mathbb{R}^{m\times n}$. The matrix pair $(A,C)$ is \textit{discrete-time semiobservable} if
\begin{eqnarray} 
\bigcap_{k=0}^{n-1}\ker(C(I_{n}-A)^{k})=\ker(I_{n}-A).
\end{eqnarray} 
\end{definition}

Next, we present an extended version of Definition~\ref{def_so} in \cite{HZ:CDC:2012}.

\begin{definition}\label{Defod}
Let $A\in\mathbb{R}^{n\times n}$ and $C\in\mathbb{R}^{m\times n}$.
The matrix pair $(A,C)$ is \textit{discrete-time $k$-semiobservable} if there exists a nonnegative integer $k$ such that
\begin{eqnarray}
k=\min\left\{l\in\overline{\mathbb{Z}}_{+}:\bigcap_{i=0}^{n-1}\ker\left(C(I_{n}-A)^{l+i}\right)=\ker(I_{n}-A)\right\}.\label{Nulld}
\end{eqnarray}
\end{definition}

An alternative extended version of Definition~\ref{def_so} to operator pairs can be found in \cite{HB:SCL:2013}. Define $\ell_{2}$ to be the collection of all sequences $\{x_{i}\}_{i=0}^{\infty}$ for which $\sum_{i=0}^{\infty}\|x_{i}\|^{2}<\infty$, where $\|\cdot\|$ denotes the 2-norm.

\begin{definition}\label{def_op}
Consider a Hilbert space $\ell_{2}$ and a linear system $\mathcal{G}_{\rm{a}}$ with a given infinitesimal generator $\mathcal{A}$ of the form
$\frac{{\rm{d}}}{{\rm{d}}t}\psi(t)=(\mathcal{A}\psi)(t)$
over $\ell_{2}$. 
Let $\mathcal{C}$ be a bounded operator on $\ell_{2}$.
The operator pair $(\mathcal{A},\mathcal{C})$ is \textit{discretely approximate semiobservable} if 
\begin{eqnarray}\label{DAS}
\bigcap_{k=0}^{\infty}\ker(\mathcal{C}\mathcal{A}^{k})=\ker(\mathcal{A}).
\end{eqnarray}  
\end{definition}  

Motivated by Definitions \ref{def_so} and \ref{def_op}, we propose a new notion of discrete-time approximate semiobservable for a (possibly infinite) set of matrix pairs.

\begin{definition}
Let $A_{k}\in\mathbb{R}^{n\times n}$, $k=0,1,2,\ldots$, and $C\in\mathbb{R}^{m\times n}$. The set of pairs $\{(A_{k},C)\}_{k\in\overline{\mathbb{Z}}_{+}}$ is called \textit{discrete-time approximate semiobservable with respect to some matrix $A\in\mathbb{R}^{n\times n}$} if
\begin{eqnarray}
\bigcap_{k=0}^{\infty}\ker(C(I_{n}-A_{k}))=\ker(I_{n}-A).
\end{eqnarray} 
\end{definition}

The following definition of paracontracting matrices is due to \cite{EKN:LAA:1990}.

\begin{definition}
Let $W\in\mathbb{R}^{n\times n}$. $W$ is called \textit{paracontracting} if for any $x\in\mathbb{R}^{n}$, $Wx\neq x$ is equivalent to $\|Wx\|<\|x\|$. 
\end{definition}

 Recall from \cite{Bernstein:2009,HCH:2009,Hui:TAC:2013} that a matrix $A\in\mathbb{R}^{n\times n}$ is called \textit{discrete-time semistable} if ${\rm{spec}}(A)\subseteq\{s\in\mathbb{C}:|s|<1\}\cup\{1\}$, and if $1\in{\rm{spec}}(A)$, then $1$ is semisimple.
$A\in\mathbb{R}^{n\times n}$ is called \textit{nontrivially discrete-time semistable} \cite{Hui:TAC:2013} if $A$ is discrete-time semistable and $A\neq I_{n}$. Finally, $A\in\mathbb{R}^{n\times n}$ is called \textit{normal} \cite[p.~179]{Bernstein:2009} if $AA^{\mathrm{T}}=A^{\mathrm{T}}A$. 

\begin{lemma}\label{normal}
Let $W\in\mathbb{R}^{q\times q}$. If $W$ is normal and nontrivially discrete-time semistable, then $W$ is paracontracting. Conversely, if $W$ is paracontracting, then $W$ is nontrivially discrete-time semistable.
\end{lemma}

\begin{IEEEproof}
Assume that $W$ is normal and nontrivially discrete-time semistable. Since $W$ is normal, it follows from Corollary 5.4.8 of \cite[p.~321]{Bernstein:2009} that $W$ has $q$ mutually orthogonal eigenvectors. In this case, for any $x\in\mathbb{R}^{q}$, we write $x$ as $x=\sum_{i=1}^{n}\alpha_{i}y_{i}$ where $\alpha_{i}$, $i=1,\ldots,n$, are either real or complex numbers, and  $\{y_{1},\ldots,y_{n}\}$ is an orthonormal set of eigenvectors of $W$ associated with the eigenvalues $\lambda_{i}\in\mathbb{C}$, $|\lambda_{i}|<1$ or $\lambda_{i}=1$, $i=1,\ldots,n$. 

Next, since $Wx=\sum_{i=1}^{n}\alpha_{i}\lambda_{i} y_{i}$, it follows that $\|Wx\|^{2}=\sum_{i=1}^{n}\|\alpha_{i}\lambda_{i}y_{i}\|^{2}$. Hence, $Wx=x$ if and only if $\alpha_{i}=\alpha_{i}\lambda_{i}$ for every $i=1,\ldots,n$, or, equivalently, $Wx\neq x$ if and only if $\alpha_{j}\neq\alpha_{j}\lambda_{j}$ for some $j\in\{1,\ldots,n\}$. Clearly if $Wx\neq x$, then
$\|\alpha_{j}\lambda_{j}y_{j}\|<\|\alpha_{j}y_{j}\|$ for some $j\in\{1,\ldots,n\}$. Thus, $\|Wx\|^{2}<\|\alpha_{1}\lambda_{1}y_{1}\|^{2}+\cdots+\|\alpha_{j}y_{j}\|^{2}+\cdots+\|\alpha_{n}\lambda_{n}y_{n}\|^{2}\leq\sum_{i=1}^{n}\|\alpha_{i}y_{i}\|^{2}=\|x\|^{2}$, which imply that $\|Wx\|^{2}<\|x\|^{2}$. Hence, $\|Wx\|<\|x\|$. On the other hand, if $\|Wx\|<\|x\|$ for any nonzero $x\in\mathbb{R}^{q}$, then it follows from the above expressions for $\|Wx\|^{2}$ and $\|x\|^{2}$ that there exists at least one integer $j\in\{1,\ldots,n\}$ such that $\|\alpha_{j}\lambda_{j}y_{j}\|<\|\alpha_{j}y_{j}\|$, which implies that $|\alpha_{j}\lambda_{j}|\neq|\alpha_{j}|$. Suppose there exists some nonzero $x\in\mathbb{R}^{q}$ such that $Wx=x$. Then it follows that $\sum_{i=1}^{n}\alpha_{i}\lambda_{i}y_{i}=\sum_{i=1}^{n}\alpha_{i}y_{i}$, which implies that $\alpha_{i}\lambda_{i}=\alpha_{i}$ for all $i=1,\ldots,n$. However, this contradicts $|\alpha_{j}\lambda_{j}|\neq|\alpha_{j}|$. Hence, $Wx\neq x$. 

Conversely, it follows from Proposition 3.2 of \cite{NN:NM:1987} that if $W$ is paracontracting, then $\lim_{k\to\infty}W^{k}$ exists, and hence, $W$ is discrete-time semistable by \cite[p.~735]{Bernstein:2009}. Clearly $W\neq I_{q}$.    
\end{IEEEproof}

A direct consequence from Lemma~\ref{normal} is that if $W\in\mathbb{R}^{q\times q}$ is symmetric and nontrivially discrete-time semistable, then $W$ is paracontracting. Next we generalize Lemma~\ref{normal} to have a necessary and sufficient condition to guarantee paracontraction of $W$. 

\begin{lemma}\label{null}
Let $W\in\mathbb{R}^{q\times q}$ and ${\mathrm{spec}}(W)=\{\lambda_{1},\ldots,\lambda_{r}\}$, where $r$ denotes the number of distinct eigenvalues for $W$. Then $W$ is nontrivially discrete-time semistable, $\|Wx\|\leq\|x\|$ for any $Wx\neq\lambda_{i}x$ and every $i=1,\ldots,r$, and $\ker(W^{\mathrm{T}}W-I_{q})=\ker((W-I_{q})^{\mathrm{T}}(W-I_{q})+(W-I_{q})^{2})$ 
if and only if $W$ is paracontracting.
\end{lemma}

\begin{IEEEproof}
First, note that $\|Wx\|^{2}-\|x\|^{2}=x^{\mathrm{T}}(W^{\mathrm{T}}W-I_{q})x=x^{\mathrm{T}}[(W-I_{q})^{\mathrm{T}}(W-I_{q})+W^{\mathrm{T}}-I_{q}+W-I_{q}]x$ for any $x\in\mathbb{R}^{q}$. Hence, $W$ is paracontracting if and only if $x^{\mathrm{T}}[(W-I_{q})^{\mathrm{T}}(W-I_{q})+W^{\mathrm{T}}-I_{q}+W-I_{q}]x<0$ is equivalent to $(W-I_{q})x\neq0$, $x\in\mathbb{R}^{q}$, or, equivalently speaking,  $x^{\mathrm{T}}[(W-I_{q})^{\mathrm{T}}(W-I_{q})+W^{\mathrm{T}}-I_{q}+W-I_{q}]x\leq0$ for any $x\in\mathbb{R}^{q}$, and $x^{\mathrm{T}}[(W-I_{q})^{\mathrm{T}}(W-I_{q})+W^{\mathrm{T}}-I_{q}+W-I_{q}]x=0$ is equivalent to $(W-I_{q})x=0$. Furthermore, since $x^{\mathrm{T}}[(W-I_{q})^{\mathrm{T}}(W-I_{q})+W^{\mathrm{T}}-I_{q}+W-I_{q}]x\leq0$ for any $x\in\mathbb{R}^{q}$ is equivalent to $(W-I_{q})^{\mathrm{T}}(W-I_{q})+W^{\mathrm{T}}-I_{q}+W-I_{q}\leq0$, it follows that $W$ is paracontracting if and only if  $(W-I_{q})^{\mathrm{T}}(W-I_{q})+W^{\mathrm{T}}-I_{q}+W-I_{q}\leq0$, and $x^{\mathrm{T}}[(W-I_{q})^{\mathrm{T}}(W-I_{q})+W^{\mathrm{T}}-I_{q}+W-I_{q}]x=0$ is equivalent to $(W-I_{q})x=0$. Next, it follows from Fact 8.15.2 of \cite[p.~550]{Bernstein:2009} that the condition, $(W-I_{q})^{\mathrm{T}}(W-I_{q})+W^{\mathrm{T}}-I_{q}+W-I_{q}\leq0$ and $x^{\mathrm{T}}[(W-I_{q})^{\mathrm{T}}(W-I_{q})+W^{\mathrm{T}}-I_{q}+W-I_{q}]x=0$ if and only if $(W-I_{q})x=0$, is equivalent to a new condition,  $(W-I_{q})^{\mathrm{T}}(W-I_{q})+W^{\mathrm{T}}-I_{q}+W-I_{q}\leq0$ and $\ker((W-I_{q})^{\mathrm{T}}(W-I_{q})+W^{\mathrm{T}}-I_{q}+W-I_{q})=\ker(W-I_{q})$. Consequently, $W$ is paracontracting if and only if $(W-I_{q})^{\mathrm{T}}(W-I_{q})+W^{\mathrm{T}}-I_{q}+W-I_{q}\leq0$ and $\ker((W-I_{q})^{\mathrm{T}}(W-I_{q})+W^{\mathrm{T}}-I_{q}+W-I_{q})=\ker(W-I_{q})$.

Assume that $W$ is nontrivially discrete-time semistable and $\ker(W^{\mathrm{T}}W-I_{q})=\ker((W-I_{q})^{\mathrm{T}}(W-I_{q})+(W-I_{q})^{2})$. We first claim that if $W$ is discrete-time semistable, then $\ker(W-I_{q})=\ker((W-I_{q})^{2})$. Since $W$ is discrete-time semistable, it follows from Proposition 11.10.2 of \cite[p.~735]{Bernstein:2009} that $W-I_{q}$ is group invertible. Now it follows from Fact 3.6.1 of \cite[p.~191]{Bernstein:2009} that $\ker(W-I_{q})=\ker((W-I_{q})^{2})$. Since $\ker(W-I_{q})=\ker((W-I_{q})^{\mathrm{T}}(W-I_{q}))$, it follows that $\ker(W-I_{q})=\ker((W-I_{q})^{2})=\ker((W-I_{q})^{\mathrm{T}}(W-I_{q}))$.

We now claim that $(W-I_{q})^{\mathrm{T}}(W-I_{q})+W^{\mathrm{T}}-I_{q}+W-I_{q}\leq0$, or equivalently, $W^{\mathrm{T}}W\leq I_{q}$. Clearly by discrete-time semistability of $W$, $|\lambda_{i}|\leq 1$ for every $i=1,\ldots,r$. Next by definition $Wx_{i}=\lambda_{i}x_{i}$ for $x_{i}\in\ker(\lambda_{i}I_{q}-W)\backslash\{0\}$, $i=1,\ldots,r$. Hence, $x_{i}^{*}W^{\mathrm{T}}Wx_{i}=|\lambda_{i}|^{2}x_{i}^{*}x_{i}$, $i=1,\ldots,r$. By Proposition 4.5.4 of \cite[p.~268]{Bernstein:2009}, $x_{1},\ldots,x_{r}$ are linearly independent, and hence, $\ker(\lambda_{i}I_{q}-W)\cap\ker(\lambda_{j}I_{q}-W)=\{0\}$ for every $i,j=1,\ldots,r$, $i\neq j$. Then it follows from $|\lambda_{i}|\leq 1$, $i=1,\ldots,r$, that $x^{*}W^{\mathrm{T}}Wx\leq x^{*}x$ for every $x\in\bigcup_{i=1}^{r}\ker(\lambda_{i}I_{q}-W)$. 

Suppose that there exists $y\in\overline{\bigcup_{i=1}^{r}\ker(\lambda_{i}I_{q}-W)}$ such that $y^{*}W^{\mathrm{T}}Wy>y^{*}y$, where $\overline{S}$ denotes the complement of the set $S$. First note that $\overline{\bigcup_{i=1}^{r}\ker(\lambda_{i}I_{q}-W)}=\bigcap_{i=1}^{r}\overline{\ker(\lambda_{i}I_{q}-W)}$. Hence, $y\in\overline{\ker(\lambda_{i}I_{q}-W)}$ for every $i=1,\ldots,r$, or equivalently, $Wy\neq\lambda_{i}y$ for every $i=1,\ldots,r$. However, this contradicts the condition that $\|Wx\|\leq\|x\|$ for any $Wx\neq\lambda_{i}x$ and every $i=1,\ldots,r$. In summary, $x^{*}W^{\mathrm{T}}Wx\leq x^{*}x$ for every $x\in\mathbb{C}^{q}$. Thus, $W^{\mathrm{T}}W\leq I_{q}$.

Next, we show that $\ker((W-I_{q})^{\mathrm{T}}(W-I_{q})+W^{\mathrm{T}}-I_{q}+W-I_{q})=\ker(W-I_{q})$. If $x\in\ker(W-I_{q})$, then it follows from $\ker(W-I_{q})=\ker((W-I_{q})^{2})=\ker((W-I_{q})^{\mathrm{T}}(W-I_{q}))$ that  $((W-I_{q})^{\mathrm{T}}(W-I_{q})+(W-I_{q})^{2})x=0$, which implies that $x\in\ker((W-I_{q})^{\mathrm{T}}(W-I_{q})+(W-I_{q})^{2})$. Hence, $\ker(W-I_{q})\subseteq\ker((W-I_{q})^{\mathrm{T}}(W-I_{q})+(W-I_{q})^{2})$. On the other hand, if $x\in\ker((W-I_{q})^{\mathrm{T}}(W-I_{q})+(W-I_{q})^{2})$, then it follows from $(W-I_{q})^{\mathrm{T}}(W-I_{q})+W^{\mathrm{T}}-I_{q}+W-I_{q}\leq0$ that $0=x^{\mathrm{T}}(W-I_{q})^{\mathrm{T}}((W-I_{q})^{\mathrm{T}}(W-I_{q})+(W-I_{q})^{2})x=x^{\mathrm{T}}(W-I_{q})^{\mathrm{T}}(W^{\mathrm{T}}-I_{q}+W-I_{q})(W-I_{q})x\leq-x^{\mathrm{T}}((W-I_{q})^{2})^{\mathrm{T}}(W-I_{q})^{2}x\leq0$, which implies that $x^{\mathrm{T}}((W-I_{q})^{2})^{\mathrm{T}}(W-I_{q})^{2}x=0$, and hence, $x\in\ker((W-I_{q})^{2})=\ker(W-I_{q})$. Thus, $\ker((W-I_{q})^{\mathrm{T}}(W-I_{q})+(W-I_{q})^{2})\subseteq\ker(W-I_{q})$. Therefore, $\ker((W-I_{q})^{\mathrm{T}}(W-I_{q})+(W-I_{q})^{2})=\ker(W-I_{q})$. Now it follows from $\ker(W^{\mathrm{T}}W-I_{q})=\ker((W-I_{q})^{\mathrm{T}}(W-I_{q})+(W-I_{q})^{2})$ that $\ker((W-I_{q})^{\mathrm{T}}(W-I_{q})+W^{\mathrm{T}}-I_{q}+W-I_{q})=\ker(W-I_{q})$. Combining this kernel condition with $W^{\mathrm{T}}W\leq I_{q}$ yields paracontraction of $W$.

Alternatively, assume that $W$ is paracontracting. Then it follows from Lemma~\ref{normal} that $W$ is nontrivially discrete-time semistable. Hence, $W-I_{q}$ is group invertible. Moreover, it follows from Fact 3.6.1 of \cite[p.~191]{Bernstein:2009} that $\ker(W-I_{q})=\ker((W-I_{q})^{2})=\ker((W-I_{q})^{\mathrm{T}}(W-I_{q}))$. Next, since for any $x\in\mathbb{R}^{q}$, $Wx\neq x$ is equivalent to $\|Wx\|<\|x\|$, it follows that $\|Wx\|\leq\|x\|$ for every $x\in\mathbb{R}^{q}$, or equivalently, $W^{\mathrm{T}}W\leq I_{q}$. In particular, $\|Wx\|\leq\|x\|$ for any $Wx\neq\lambda_{i}x$ and every $i=1,\ldots,r$. Finally, to show that $\ker(W^{\mathrm{T}}W-I_{q})=\ker((W-I_{q})^{\mathrm{T}}(W-I_{q})+(W-I_{q})^{2})$, it suffices to show that $\ker((W-I_{q})^{\mathrm{T}}(W-I_{q})+(W-I_{q})^{2})=\ker(W-I_{q})$ since $\ker(W^{\mathrm{T}}W-I_{q})=\ker(W-I_{q})$ by paracontraction of $W$. This has actually been done in the above paragraph.
\end{IEEEproof}

Next, we replace $\|Wx\|\leq\|x\|$ for any $Wx\neq\lambda_{i}x$ and every $i=1,\ldots,r$, and  $\ker(W^{\mathrm{T}}W-I_{q})=\ker((W-I_{q})^{\mathrm{T}}(W-I_{q})+(W-I_{q})^{2})$ in Lemma~\ref{null} by new conditions which are easier to check practically. Recall from \cite[p.~608]{Bernstein:2009} that the H\"{o}lder-induced norm $\|\cdot\|$ for $W$ is defined by $\|W\|=\max_{x\in\mathbb{R}^{q}\backslash\{0\}}\|Ax\|/\|x\|$.

\begin{lemma}\label{pcm}
Let $W\in\mathbb{R}^{q\times q}$. Then $W$ is nontrivially discrete-time semistable, $\|W\|\leq 1$, and ${\mathrm{rank}}(W^{\mathrm{T}}W-I_{q})={\mathrm{rank}}((W-I_{q})^{\mathrm{T}}(W-I_{q})+(W-I_{q})^{2})={\mathrm{rank}}\left[\begin{array}{cc}
W^{\mathrm{T}}W-I_{q} & (W-I_{q})^{\mathrm{T}}(W-I_{q})+(W^{\mathrm{T}}-I_{q})^{2}
\end{array}\right]$ 
if and only if $W$ is paracontracting.
\end{lemma}

\begin{IEEEproof}
First, it follows from Proposition 9.4.9 of \cite[p.~609]{Bernstein:2009} that $\sigma_{\max}(W)=\|W\|\leq1$, where $\sigma_{\max}(W)$ denotes the maximum singular value of $W$. Next, it follows from Fact 5.11.35 of \cite[p.~358]{Bernstein:2009} that $\sigma_{\max}(W)\leq 1$ if and only if $W^{\mathrm{T}}W\leq I_{q}$. Thus, $\|W\|\leq 1$ if and only if $W^{\mathrm{T}}W\leq I_{q}$.

Second, it follows from Equation (2.4.13) of \cite[p.~103]{Bernstein:2009} that $\ker(W^{\mathrm{T}}W-I_{q})=\ker((W-I_{q})^{\mathrm{T}}(W-I_{q})+(W-I_{q})^{2})$ if and only if ${\mathrm{ran}}(W^{\mathrm{T}}W-I_{q})^{\bot}={\mathrm{ran}}((W-I_{q})^{\mathrm{T}}(W-I_{q})+(W^{\mathrm{T}}-I_{q})^{2})^{\bot}$, where ${\mathrm{ran}}(A)$ denotes the range of $A$ and $S^{\bot}$ denotes the orthogonal complement of $S$. Note that both ${\mathrm{ran}}(W^{\mathrm{T}}W-I_{q})$ and ${\mathrm{ran}}((W-I_{q})^{\mathrm{T}}(W-I_{q})+(W^{\mathrm{T}}-I_{q})^{2})$ are subspaces. Then it follows from Fact 2.9.14 of \cite[p.~121]{Bernstein:2009} that ${\mathrm{ran}}(W^{\mathrm{T}}W-I_{q})^{\bot}={\mathrm{ran}}((W-I_{q})^{\mathrm{T}}(W-I_{q})+(W^{\mathrm{T}}-I_{q})^{2})^{\bot}$ if and only if ${\mathrm{ran}}(W^{\mathrm{T}}W-I_{q})={\mathrm{ran}}((W-I_{q})^{\mathrm{T}}(W-I_{q})+(W^{\mathrm{T}}-I_{q})^{2})$. Now it follows from Fact 2.11.5 of \cite[p.~131]{Bernstein:2009} that ${\mathrm{ran}}(W^{\mathrm{T}}W-I_{q})={\mathrm{ran}}((W-I_{q})^{\mathrm{T}}(W-I_{q})+(W^{\mathrm{T}}-I_{q})^{2})$ if and only if ${\mathrm{rank}}(W^{\mathrm{T}}W-I_{q})={\mathrm{rank}}((W-I_{q})^{\mathrm{T}}(W-I_{q})+(W^{\mathrm{T}}-I_{q})^{2})={\mathrm{rank}}((W-I_{q})^{\mathrm{T}}(W-I_{q})+(W-I_{q})^{2})={\mathrm{rank}}\left[\begin{array}{cc}
W^{\mathrm{T}}W-I_{q} & (W-I_{q})^{\mathrm{T}}(W-I_{q})+(W^{\mathrm{T}}-I_{q})^{2}
\end{array}\right]$. 

Now the rest of the proof directly follows from the proof of Lemma~\ref{null}.   
\end{IEEEproof}

The following corollary is immediate based on Lemmas \ref{null} and \ref{pcm}.

\begin{corollary}\label{coro}
Let $W\in\mathbb{R}^{q\times q}$. Then $W$ is nontrivially discrete-time semistable, $\|W\|\leq 1$, and $\ker((W-I_{q})^{\mathrm{T}}(W-I_{q})+W^{\mathrm{T}}-I_{q}+W-I_{q})=\ker((W-I_{q})^{\mathrm{T}}(W-I_{q})+(W-I_{q})^{2})$ 
if and only if $W$ is paracontracting.
\end{corollary}

Motivated by Theorem 1 of \cite{EKN:LAA:1990} and Corollary 3.2 of \cite{BEN:ETNA:1994},  we have the following convergence results for a sequence of (possibly infinite) discrete-time semistable matrices. 

\begin{lemma}\label{lemma_DTSS}
Let $J$ be a (possibly infinite) countable index set and $P_{k}\in\mathbb{R}^{n\times n}$, $k\in J$, be discrete-time semistable, $\|P_{k}\|\leq1$, and $\ker(P_{k}^{\mathrm{T}}P_{k}-I_{n})=\ker((P_{k}-I_{n})^{\mathrm{T}}(P_{k}-I_{n})+(P_{k}-I_{n})^{2})$. Consider the sequence $\{x_{i}\}_{i=0}^{\infty}$ defined by the iterative process $x_{i+1}=Q_{i}x_{i}$, $i=0,1,2,\ldots$, where $Q_{i}\in\{P_{k}:\forall k\in J\}$. 
\begin{itemize}
\item[$i$)] If $|J|<\infty$, then $\lim_{i\to\infty}x_{i}$ exists. If in addition, $P_{k}\in\mathbb{R}^{n\times n}$ is nontrivially discrete-time semistable for every $k\in J$, then $\lim_{i\to\infty}x_{i}$ is in $\bigcap_{k\in\mathcal{I}}\ker(I_{n}-P_{k})$, where $\mathcal{I}$ is the set of all indexes $k$ for which $P_{k}$ appears infinitely often in $\{Q_{i}\}_{i=0}^{\infty}$.
\item[$ii$)] If there exists $s\in J$ such that $P_{s}$ is nontrivially discrete-time semistable, $\{(Q_{k},I_{n})\}_{k\in\overline{\mathbb{Z}}_{+}}$ is discrete-time approximate semiobservable with respect to some nontrivially discrete-time semistable matrix $Q_{r}$,  $r\in\overline{\mathbb{Z}}_{+}$, and for every positive integer $N$, there always exists $j\geq N$ such that $Q_{j}=Q_{r}$, then $\lim_{i\to\infty}x_{i}$ exists and the limit is in $\ker(I_{n}-Q_{r})$.
\end{itemize}
\end{lemma}

\begin{IEEEproof}
$i$) Since $P_{k}\in\mathbb{R}^{n\times n}$ is discrete-time semistable for every $k\in J$, it follows that either $P_{k}=I_{n}$ or $P_{k}$ is nontrivially discrete-time semistable. If there exists $N\geq 1$ such that $Q_{i}=I_{n}$ for all $i\geq N$, then $x_{i}=x_{N}$ for all $i\geq N$, which implies that $\lim_{i\to\infty}x_{i}$ exists. Otherwise, we select all the nontrivially discrete-time semistable matrices in $\{Q_{i}\}_{i=0}^{\infty}$ to form an infinite subsequence $\{Q_{i_{n}}\}_{n=0}^{\infty}$ of $\{Q_{i}\}_{i=0}^{\infty}$. Define $y_{n+1}=Q_{i_{n}}y_{n}$, $n=0,1,2,\ldots$. Then it follows from Corollary~\ref{coro} that $Q_{i_{n}}$ is paracontracting for every $n=0,1,2,\ldots$. Now by Theorem 1 of \cite{EKN:LAA:1990}, $\lim_{n\to\infty}y_{n}$ exists. Consequently, $\lim_{i\to\infty}x_{i}=\lim_{n\to\infty}y_{n}$ exists. The second assertion is a direct consequence of Corollary~\ref{coro} above and Theorem 1 of \cite{EKN:LAA:1990}.

$ii$) Again, it follows from Corollary~\ref{coro} that $Q_{i}$ is paracontracting for every $i=0,1,2,\ldots$. Then the assertion follows directly from Corollary 3.2 of \cite{BEN:ETNA:1994}.
\end{IEEEproof}

Now we have the main result for the global convergence of the iterative process in Algorithm~\ref{MCO}.

\begin{theorem}\label{thm_HMCO}
Consider the following discrete-time switched linear model to describe the iterative process for MCO:
\begin{eqnarray}
x_{i}[k+1]&=&x_{i}[k]+h_{k}v_{i}[k+1],\quad x_{i}[0]=x_{i0},\label{DE_1}\\
v_{i}[k+1]&=&v_{i}[k]+h_{k}\eta_{k}\sum_{j\in\mathcal{N}_{i}}(v_{j}[k]-v_{i}[k])+h_{k}\mu_{k}\sum_{j\in\mathcal{N}_{i}}(x_{j}[k]-x_{i}[k])\nonumber\\
&&+h_{k}\kappa_{k}(p[k]-x_{i}[k]),\quad v_{i}[0]=v_{i0},\\
p[k+1]&=&p[k]+h_{k}\kappa_{k}(x_{j}[k]-p[k]),\quad p[k]\not\in\mathcal{Z}_{p},\quad p[0]=p_{0},\\
p[k+1]&=&x_{j}[k],\quad p[k]\in\mathcal{Z}_{p},\quad k=0,1,2,\ldots,\quad i=1,\ldots,q,\label{MSO_4}
\end{eqnarray} where $x_{i}\in\mathbb{R}^{n}$, $v_{i}\in\mathbb{R}^{n}$, $p\in\mathbb{R}^{n}$, $\mu_{k},\eta_{k},\kappa_{k},h_{k}$ are randomly selected in $\Omega\subseteq[0,\infty)$, $\mathcal{Z}_{p}=\{p\in\mathbb{R}^{n}:f(x_{j})<f(p)\}$, and $x_{j}=\arg\min_{1\leq i\leq q}f(x_{i})$. Assume that for every $k\in\overline{\mathbb{Z}}_{+}$ and every $j=1,\ldots,q$:
\begin{itemize}
\item[H1)] $0<h_{k}<-\frac{\lambda+\bar{\lambda}}{|\lambda|^{2}}$ for every $\lambda\in\{-\kappa_{k},-\frac{\kappa_{k}(1+h_{k})}{2}\pm\frac{1}{2}\sqrt{\kappa_{k}^{2}(1+h_{k})^{2}-4\kappa_{k}},-\frac{\kappa_{k}h_{k}}{2}\pm\frac{1}{2}\sqrt{\kappa_{k}^{2}h_{k}^{2}-4\kappa_{k}},\lambda\in\mathbb{C}:\forall \frac{\lambda^{2}+\kappa_{k} h_{k}\lambda+\kappa_{k}}{\eta_{k}\lambda+\mu_{k} h_{k}\lambda+\mu_{k}}\in{\mathrm{spec}}(-L_{k})\backslash\{0\}\}$;
\item[H2)] $0<h_{k}<-\frac{\lambda+\bar{\lambda}}{|\lambda|^{2}}$ for every $\lambda\in\{-1,-\frac{h_{k}^{2}\kappa_{k}}{2}\pm\frac{1}{2}\sqrt{(h_{k}^{2}\kappa_{k})^{2}-4h_{k}^{2}\kappa_{k}},\lambda_{1},\lambda_{2}\in\mathbb{C}:\forall\frac{\lambda_{1}^{2}+\kappa_{k} h_{k}^{2}\lambda_{1}+\kappa_{k}h_{k}^{2}}{\eta_{k}h_{k}\lambda_{1}+\mu_{k} h_{k}^{2}\lambda_{1}+\mu_{k}h_{k}^{2}}\in{\mathrm{spec}}(-L_{k})\backslash\{0\},\lambda_{2}^{3}+(1+h_{k}^{2}\kappa_{k})\lambda_{2}^{2}+(2h_{k}^{2}\kappa_{k}-h_{k}\kappa_{k})\lambda_{2}+h_{k}^{2}\kappa_{k}=0\}$;
\item[H3)] $\|I_{2nq+n}+h_{k}A_{k}^{[j]}+h_{k}^{2}A_{{\mathrm{c}}k}\|\leq 1$ and $\|I_{2nq+n}+B_{k}^{[j]}+h_{k}^{2}A_{{\mathrm{c}}k}\|\leq1$.
\item[H4)] $\ker((h_{k}A_{k}^{[j]}+h_{k}^{2}A_{{\mathrm{c}}k})^{\mathrm{T}}(h_{k}A_{k}^{[j]}+h_{k}^{2}A_{{\mathrm{c}}k})+(h_{k}A_{k}^{[j]}+h_{k}^{2}A_{{\mathrm{c}}k})^{\mathrm{T}}+h_{k}A_{k}^{[j]}+h_{k}^{2}A_{{\mathrm{c}}k})=\ker((h_{k}A_{k}^{[j]}+h_{k}^{2}A_{{\mathrm{c}}k})^{\mathrm{T}}\\(h_{k}A_{k}^{[j]}+h_{k}^{2}A_{{\mathrm{c}}k})+(h_{k}A_{k}^{[j]}+h_{k}^{2}A_{{\mathrm{c}}k})^{2})$ and $\ker((B_{k}^{[j]}+h_{k}^{2}A_{{\mathrm{c}}k})^{\mathrm{T}}(B_{k}^{[j]}+h_{k}^{2}A_{{\mathrm{c}}k})+(B_{k}^{[j]}+h_{k}^{2}A_{{\mathrm{c}}k})^{\mathrm{T}}+B_{k}^{[j]}+h_{k}^{2}A_{{\mathrm{c}}k})=\ker((B_{k}^{[j]}+h_{k}^{2}A_{{\mathrm{c}}k})^{\mathrm{T}}(B_{k}^{[j]}+h_{k}^{2}A_{{\mathrm{c}}k})+(B_{k}^{[j]}+h_{k}^{2}A_{{\mathrm{c}}k})^{2})$.
\end{itemize}
Then the following conclusions hold: 
\begin{itemize}
\item[C1)] If $\Omega$ is a finite discrete set, then $x_{i}[k]\to p^{\dag}$, $v_{i}[k]\to \textbf{0}_{n\times 1}$, and $p[k]\to p^{\dag}$ as $k\to\infty$ for every $x_{i0}\in\mathbb{R}^{n}$, $v_{i0}\in\mathbb{R}^{n}$, $p_{0}\in\mathbb{R}^{n}$, and every $i=1,\ldots,q$, where $p^{\dag}\in\mathbb{R}^{n}$ is some constant vector.
\item[C2)] If for every positive integer $N$, there always exists $s\geq N$ such that  $h_{s}(A_{s}^{[j_{s}]}+h_{s}A_{{\mathrm{c}}s})=B_{s}^{[j_{s}]}+h_{s}^{2}A_{{\mathrm{c}}s}=h_{T}(A_{T}^{[j_{T}]}+h_{T}A_{{\mathrm{c}}T})=B_{T}^{[j_{T}]}+h_{T}^{2}A_{{\mathrm{c}}T}$ for some fixed $T\in\overline{\mathbb{Z}}_{+}$, where $j_{s},j_{T}\in\{1,\ldots,q\}$, then $x_{i}[k]\to p^{\dag}$, $v_{i}[k]\to \textbf{0}_{n\times 1}$, and $p[k]\to p^{\dag}$ as $k\to\infty$ for every $x_{i0}\in\mathbb{R}^{n}$, $v_{i0}\in\mathbb{R}^{n}$, $p_{0}\in\mathbb{R}^{n}$, and every $i=1,\ldots,q$, where $p^{\dag}\in\mathbb{R}^{n}$ is some constant vector.
\end{itemize}
\end{theorem}

\begin{IEEEproof}
Let $Z=[x_{1}^{\rm{T}},\ldots,x_{q}^{\rm{T}},v_{1}^{\rm{T}},\ldots,v_{q}^{\rm{T}},p^{\rm{T}}]^{\rm{T}}\in\mathbb{R}^{2nq+n}$. Note that (\ref{DE_1})--(\ref{MSO_4}) can be rewritten as the compact form $Z[k+1]=(I_{2nq+n}+h_{k}(A_{k}^{[j_{k}]}+h_{k}A_{{\mathrm{c}}k}))Z[k]$, $Z[k]\not\in\mathcal{S}$, and $Z[k+1]=(I_{2nq+n}+B_{k}^{[j_{k}]}+h_{k}^{2}A_{{\mathrm{c}}k})Z[k]$, $Z[k]\in\mathcal{S}$, $j_{k}\in\{1,\ldots,q\}$ is selected based on $\mathcal{Z}_{p}$. Let $h_{k}^{\dag}=\min\Big\{-\frac{\lambda+\bar{\lambda}}{|\lambda|^{2}}:\lambda\in\{-\kappa_{k},-\frac{\kappa_{k}(1+h_{k})}{2}\pm\frac{1}{2}\sqrt{\kappa_{k}^{2}(1+h_{k})^{2}-4\kappa_{k}},-\frac{\kappa_{k}h_{k}}{2}\pm\frac{1}{2}\sqrt{\kappa_{k}^{2}h_{k}^{2}-4\kappa_{k}},\lambda\in\mathbb{C}:\forall \frac{\lambda^{2}+\kappa_{k} h_{k}\lambda+\kappa_{k}}{\eta_{k}\lambda+\mu_{k} h_{k}\lambda+\mu_{k}}\in{\mathrm{spec}}(-L_{k})\backslash\{0\}\}\Big\}$. First, we show that if $h<h_{k}^{\dag}$, then $I_{2nq+n}+h_{k}(A_{k}^{[j]}+h_{k}A_{{\mathrm{c}}k})$ becomes discrete-time semistable for every $j=1,\ldots,q$ and every $k=0,1,2,\ldots$. Note that  ${\mathrm{spec}}(I_{2nq+n}+h_{k}(A_{k}^{[j]}+h_{k}A_{{\mathrm{c}}k}))=\{1+h\lambda:\forall\lambda\in{\mathrm{spec}}(A_{k}^{[j]}+h_{k}A_{{\mathrm{c}}k})\}$. Since by Lemma~\ref{lemma_A} and Assumption H1, $A_{k}^{[j]}+h_{k}A_{{\mathrm{c}}k}$ is semistable for every $j=1,\ldots,q$ and every $k=0,1,2,\ldots$, it follows that ${\mathrm{spec}}(I_{2nq+n}+h_{k}(A_{k}^{[j]}+h_{k}A_{{\mathrm{c}}k}))=\{1\}\cup\{1+h\lambda:\forall\lambda\in{\mathrm{spec}}(A_{k}^{[j]}+h_{k}A_{{\mathrm{c}}k}),{\mathrm{Re}}\,\lambda<0\}$. Hence, $I_{2nq+n}+h_{k}(A_{k}^{[j]}+h_{k}A_{{\mathrm{c}}k})$ is discrete-time semistable for every $j=1,\ldots,q$ and every $k=0,1,2,\ldots$ if $|1+h_{k}\lambda|<1$ for every $\lambda\in{\mathrm{spec}}(A_{k}^{[j]}+h_{k}A_{{\mathrm{c}}k})$ and ${\mathrm{Re}}\,\lambda<0$. Note that $|1+h_{k}\lambda|<1$ is equivalent to $(1+h_{k}\lambda)(1+h_{k}\bar{\lambda})=|1+h_{k}\lambda|^{2}<1$, i.e., $h_{k}<-(\lambda+\bar{\lambda})/|\lambda|^{2}$. By Lemma~\ref{lemma_A}, for any $h_{k}<h_{k}^{\dag}$, $I_{2nq+n}+h_{k}(A_{k}^{[j]}+h_{k}A_{{\mathrm{c}}k})$ is discrete-time semistable for every $j=1,\ldots,q$ and every $k=0,1,2,\ldots$. Similarly, it follows from Lemma~\ref{lemma_B} and Assumption H2 that $I_{2nq+n}+B_{k}^{[j]}+h_{k}^{2}A_{{\mathrm{c}}k}$ is discrete-time semistable for every $j=1,\ldots,q$ and every $k=0,1,2,\ldots$. And (\ref{DE_1})--(\ref{MSO_4}) can further be rewritten as an iteration $Z[k+1]=P_{k}Z[k]$, $k=0,1,2,\ldots$, where $P_{k}\in\{I_{2nq+n}+h_{k}(A_{k}^{[j]}+h_{k}A_{{\mathrm{c}}k}),I_{2nq+n}+B_{k}^{[j]}+h_{k}^{2}A_{{\mathrm{c}}k}:j=1,\ldots,q,k=0,1,2,\ldots\}=\{I_{2nq+n}+h_{k}(A_{k}^{[j]}+h_{k}A_{{\mathrm{c}}k}),I_{2nq+n}+B_{k}^{[j]}+h_{k}^{2}A_{{\mathrm{c}}k}:j=1,\ldots,q,\mu_{k},\eta_{k},\kappa_{k},h_{k}\in\Omega\}$. 

C1) By assumption, $\Omega$ is a finite discrete set. Hence, $\{I_{2nq+n}+h_{k}(A_{k}^{[j]}+h_{k}A_{{\mathrm{c}}k}),I_{2nq+n}+B_{k}^{[j]}+h_{k}^{2}A_{{\mathrm{c}}k}:j=1,\ldots,q,\mu_{k},\eta_{k},\kappa_{k},h_{k}\in\Omega\}$ is a finite discrete set. Now it follows from Assumptions H3 and H4 as well as $i$) of Lemma \ref{lemma_DTSS} that $\lim_{k\to\infty}Z[k]$ exists. The rest of the conclusion follows directly from (\ref{DE_1})--(\ref{MSO_4}).

C2) By assumption, either $h_{T}(A_{T}^{[j_{T}]}+h_{T}A_{{\mathrm{c}}T})$ or $B_{T}^{[j_{T}]}+h_{T}^{2}A_{{\mathrm{c}}T}$ appears infinitely many times in the sequence $\{P_{k}\}_{k=0}^{\infty}$. Next, it follows from Lemmas \ref{lemma_Arank} and \ref{lemma_Ah} as well as the assumption $h_{k}>0$ that $\ker(h_{k}(A_{k}^{[j_{k}]}+h_{k}A_{{\mathrm{c}}k}))=\ker(A_{k}^{[j_{k}]})=\ker(A_{s}^{[j_{s}]})=\ker(h_{s}(A_{s}^{[j_{s}]}+h_{s}A_{{\mathrm{c}}s}))$ for every $k,s\in\overline{\mathbb{Z}}_{+}$. Using the similar arguments, one can prove that $\ker(B_{k}^{[j_{k}]}+h_{k}^{2}A_{{\mathrm{c}}k})=\ker(B_{k}^{[j_{k}]})=\ker(B_{s}^{[j_{s}]})=\ker(B_{s}^{[j_{s}]}+h_{s}^{2}A_{{\mathrm{c}}s})$ for every $k,s\in\overline{\mathbb{Z}}_{+}$. Hence, it follows from Assumptions H3 and H4 as well as $ii$) of Lemma \ref{lemma_DTSS} that $\lim_{k\to\infty}Z[k]$ exists. The rest of the conclusion follows directly from (\ref{DE_1})--(\ref{MSO_4}). Note that in this case, $\Omega$ may be an infinite set.
\end{IEEEproof}

\begin{remark}
Since $\rho(A)\leq\|A\|$, where $\rho(A)$ denotes the spectrum abscissa of $A$, it follows from Lemmas \ref{lemma_A}
and \ref{lemma_B} that $\|I_{2nq+n}+h_{k}A_{k}^{[j]}+h_{k}^{2}A_{{\mathrm{c}}k}\|\geq 1$ and $\|I_{2nq+n}+B_{k}^{[j]}+h_{k}^{2}A_{{\mathrm{c}}k}\|\geq1$. Hence, to guarantee H3, one only needs to assume that $\|I_{2nq+n}+h_{k}A_{k}^{[j]}+h_{k}^{2}A_{{\mathrm{c}}k}\|=1$ and $\|I_{2nq+n}+B_{k}^{[j]}+h_{k}^{2}A_{{\mathrm{c}}k}\|=1$. \hfill$\blacklozenge$
\end{remark}

\section{Numerical Evaluation}\label{ne}

\subsection{Test Function Review}\label{tf}

In order to show the performance of the parallel MCO, we conduct a numerical comparison evaluation between the standard PSO, serial MCO, and parallel MCO. In particular, we use the following eight test functions chosen from
\cite{IoanCristian2003317,ZSQ:DSCC:2011} to evaluate the proposed algorithm.
\begin{itemize}
\item Sphere function:
$f(x)=\sum_{i=1}^{n}x_i^2$.
The test area is usually restricted to the hypercube $-30\leq x_{i}\leq 30$, $i=1, \ldots, n$. The global minimum of $f(x)$ is 0 at $x_i=0$.
\item Rosenbrock's valley:
$f(x)=\sum_{i=1}^{n-1}[100(x_{i+1}-x_i^2)^2+(1-x_i)^2]$.
The test area is usually restricted to the hypercube $-30 \leq x_i \leq
30$, $i=1, \ldots, n$. The global minimum of $f(x)$ is 0 at $x_i=1$.
\item Rastrigin function:
$f(x)=10n+\sum_{i=1}^n[x_i^2-10\cos(2 \pi x_i )]$.
The test area is usually restricted to the hypercube $-30 \leq x_i \leq
30$, $i=1, \ldots, n$. The global minimum of $f(x)$ is 0.
\item Griewank function:
$f(x)=\frac{1}{4000}\sum_{i=1}^nx_i^2-\prod_{i=1}^{n}
\cos(\frac{x_i}{\sqrt{i}})+1$.
The test area is usually restricted to the hypercube $-600 \leq x_i \leq
600$, $i=1, \ldots, n$. The global minimum of $f(x)$ is 0 at $x_i=0$.
\item Ackley function:
$f(x)=-20 \exp(-0.2\times
\sqrt{\frac{1}{n}\sum_{i=1}^nx_i^2})-\exp(\frac{1}{n}\sum_{i=1}^n\cos{2\pi x_i})+20+e$.
The test area is usually restricted to the hypercube $-32.768 \leq x_i \leq
32.768$, $i=1,\ldots,n$. The global minimum of $f(x)$ is 0 at $x_i=0$.
\item De Jong's f4 function:
$f(x)=\sum_{i=1}^n(ix_i^4)$.
The test area is usually restricted to the hypercube $-20 \leq x_i \leq 20$,
$i=1, \ldots, n$. The global minimum of $f(x)$ is 0 at $x_i=0$.
\item Zakharov function:
$f(x)=\sum_{i=1}^{n}x_i^2+(0.5ix_i)^2+(0.5ix_i)^4$. The test area is usually restricted to the hypercube $-10 \leq x_i \leq 10$,
$i=1, \ldots, n$. The global minimum of $f(x)$ is 0 at $x_i=0$.
\item Levy function:
$f(x)=\sin^2(\pi x_1)+(x_n-1)^2(1+\sin^2(2 \pi x_n))-\sum_{i=1}^{n-1}(x_i-1)^2(1+10\sin^2(\pi x_i+1))$. The test area is usually restricted to the hypercube $-10 \leq x_i \leq 10$,
$i=1, \ldots, n$. The global minimum of $f(x)$ is 0 at $x_i=1$.
\end{itemize}

\subsection{Evaluation of Computational Time for the Parallel MCO}

\begin{figure} \centering
\subfigure[MCO saved time]{
\includegraphics[height=0.35\textwidth, width=0.45\textwidth]{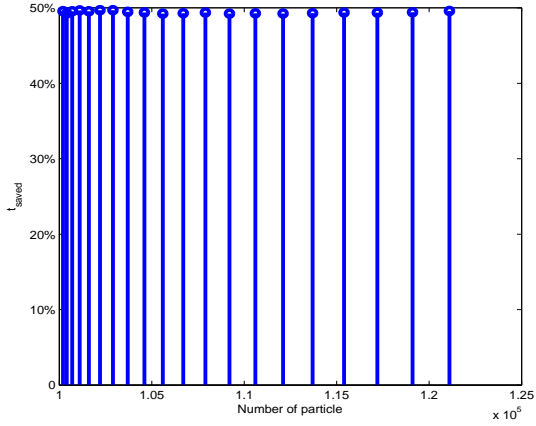}
 }
 \subfigure[MCO speedup]{
\includegraphics[height=0.35\textwidth, width=0.45\textwidth]{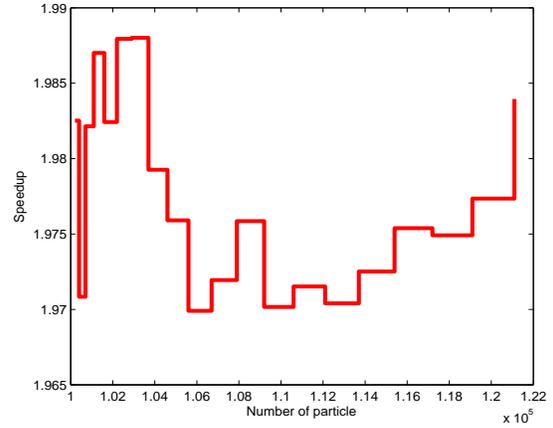}
}
\subfigure[PSO saved time]{
\includegraphics[height=0.35\textwidth, width=0.45\textwidth]{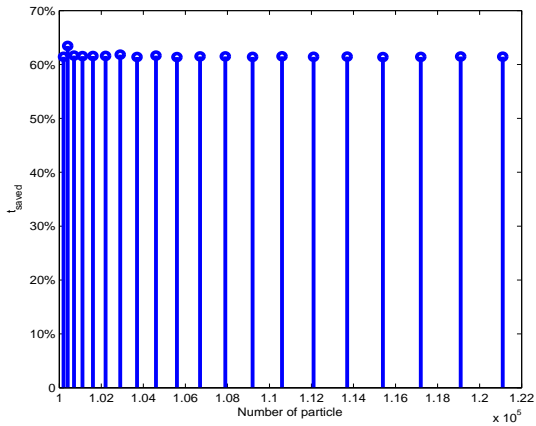}
 }
 \subfigure[PSO speedup]{
\includegraphics[height=0.35\textwidth, width=0.45\textwidth]{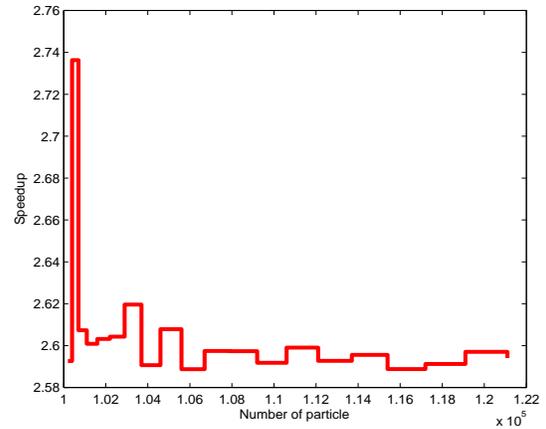}
}
\caption[]{Parallel MCO and PSO solving Sphere function.}
\label{fsphere}
\end{figure}

\begin{figure} \centering
\subfigure[MCO saved time]{
\includegraphics[height=0.35\textwidth, width=0.45\textwidth]{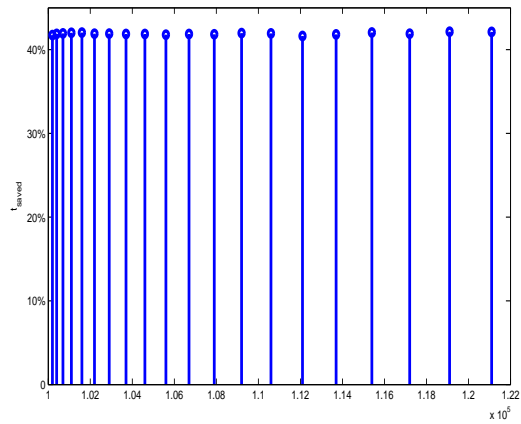}
 }
 \subfigure[MCO speedup]{
\includegraphics[height=0.35\textwidth, width=0.45\textwidth]{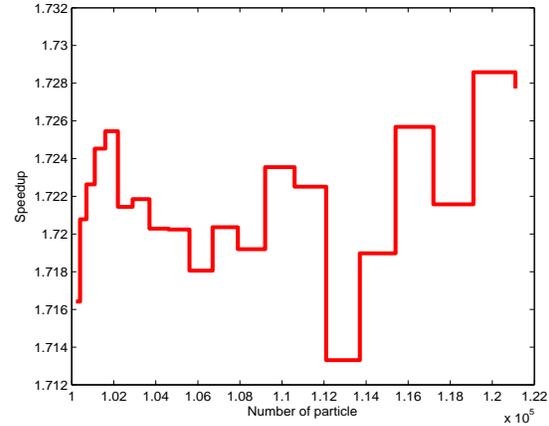}
}
\subfigure[PSO saved time]{
\includegraphics[height=0.35\textwidth, width=0.45\textwidth]{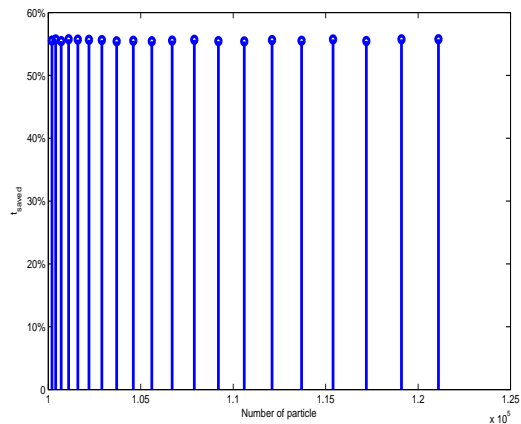}
 }
 \subfigure[PSO speedup]{
\includegraphics[height=0.35\textwidth, width=0.45\textwidth]{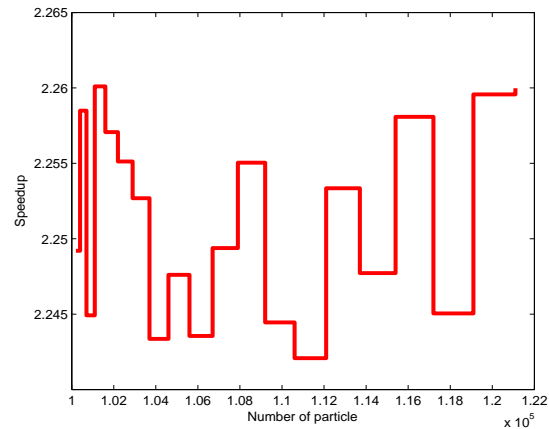}
}
\caption[]{Parallel MCO and PSO solving Rosenbrock function.}
\label{fRosenbrock}
\end{figure}

\begin{figure} \centering
\subfigure[MCO saved time]{
\includegraphics[height=0.35\textwidth, width=0.45\textwidth]{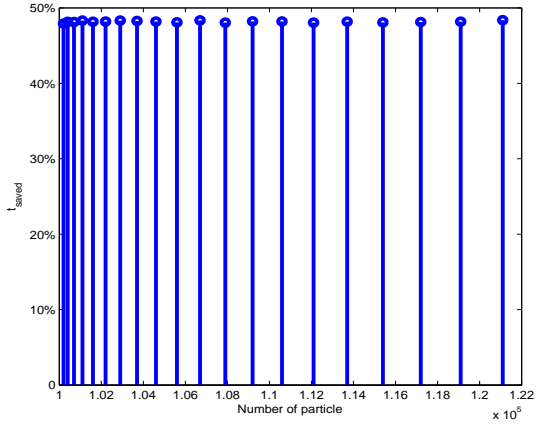}
 }
 \subfigure[MCO speedup]{
\includegraphics[height=0.35\textwidth, width=0.45\textwidth]{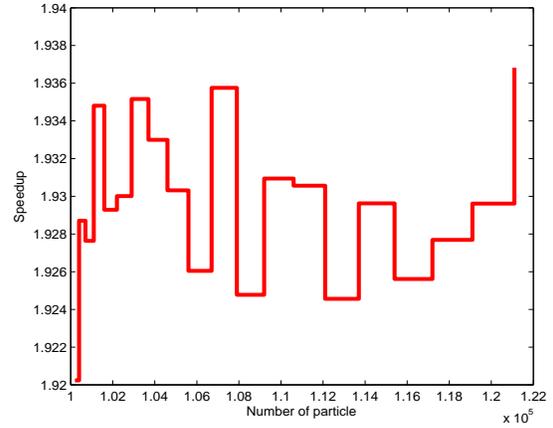}
}
\subfigure[PSO saved time]{
\includegraphics[height=0.35\textwidth, width=0.45\textwidth]{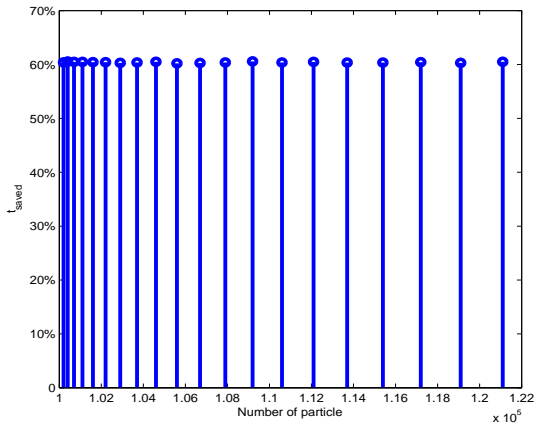}
 }
 \subfigure[PSO speedup]{
\includegraphics[height=0.35\textwidth, width=0.45\textwidth]{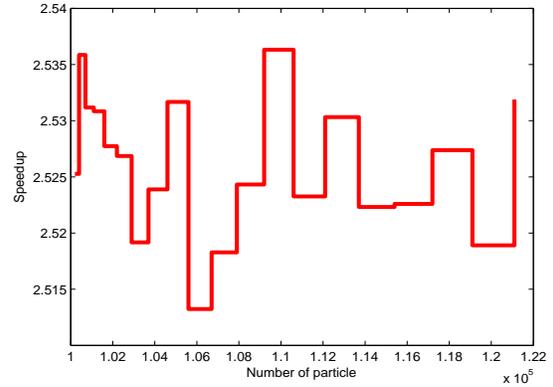}
}
\caption[]{Parallel MCO and PSO solving Rastrigin function.}
\label{fRRastrigin}
\end{figure}

\begin{figure} \centering
\subfigure[MCO saved time]{
\includegraphics[height=0.35\textwidth, width=0.45\textwidth]{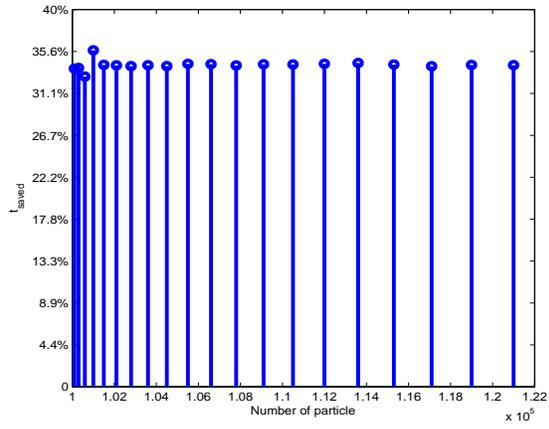}
 }
 \subfigure[MCO speedup]{
\includegraphics[height=0.35\textwidth, width=0.45\textwidth]{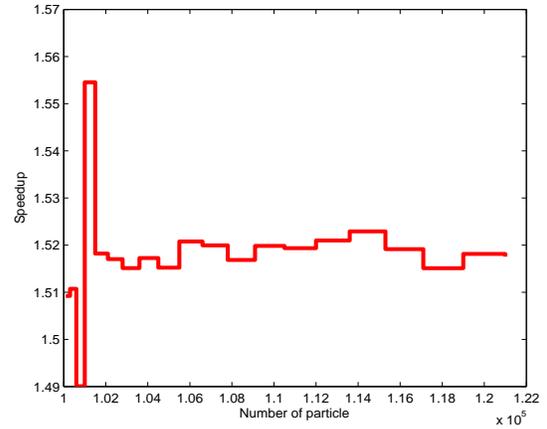}
}
\subfigure[PSO saved time]{
\includegraphics[height=0.35\textwidth, width=0.45\textwidth]{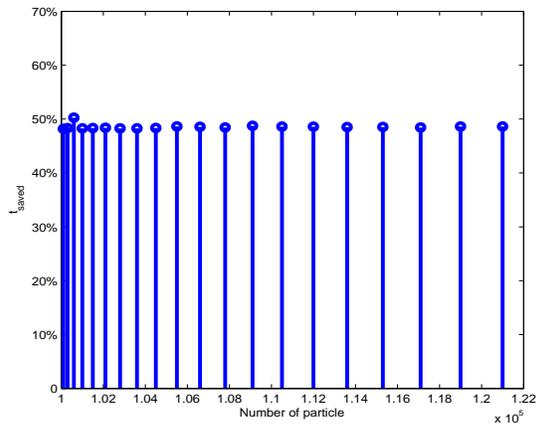}
 }
 \subfigure[PSO speedup]{
\includegraphics[height=0.35\textwidth, width=0.45\textwidth]{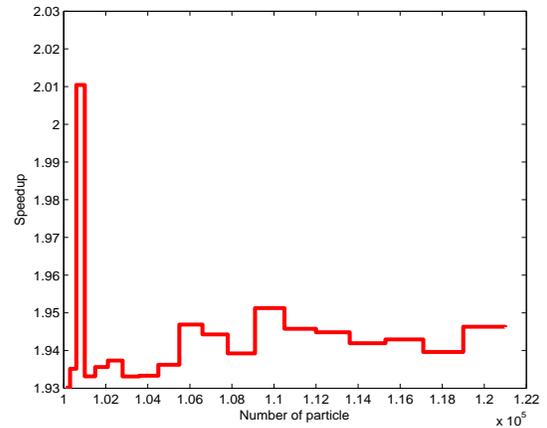}
}
\caption[]{Parallel MCO and PSO solving Griewank function.}
\label{fGriewank}
\end{figure}

\begin{figure} \centering
\subfigure[MCO saved time]{
\includegraphics[height=0.35\textwidth, width=0.45\textwidth]{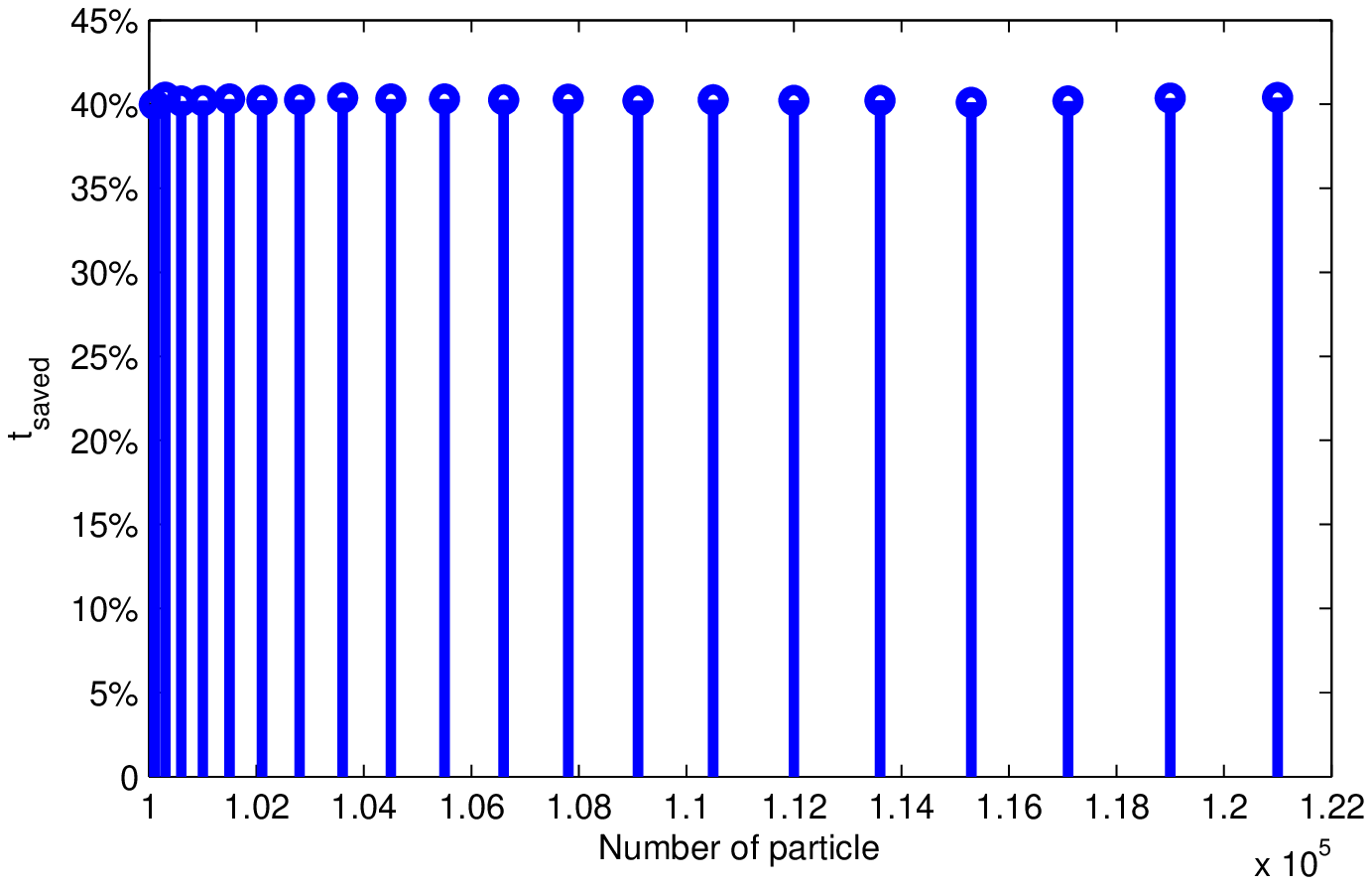}
 }
 \subfigure[MCO speedup]{
\includegraphics[height=0.35\textwidth, width=0.45\textwidth]{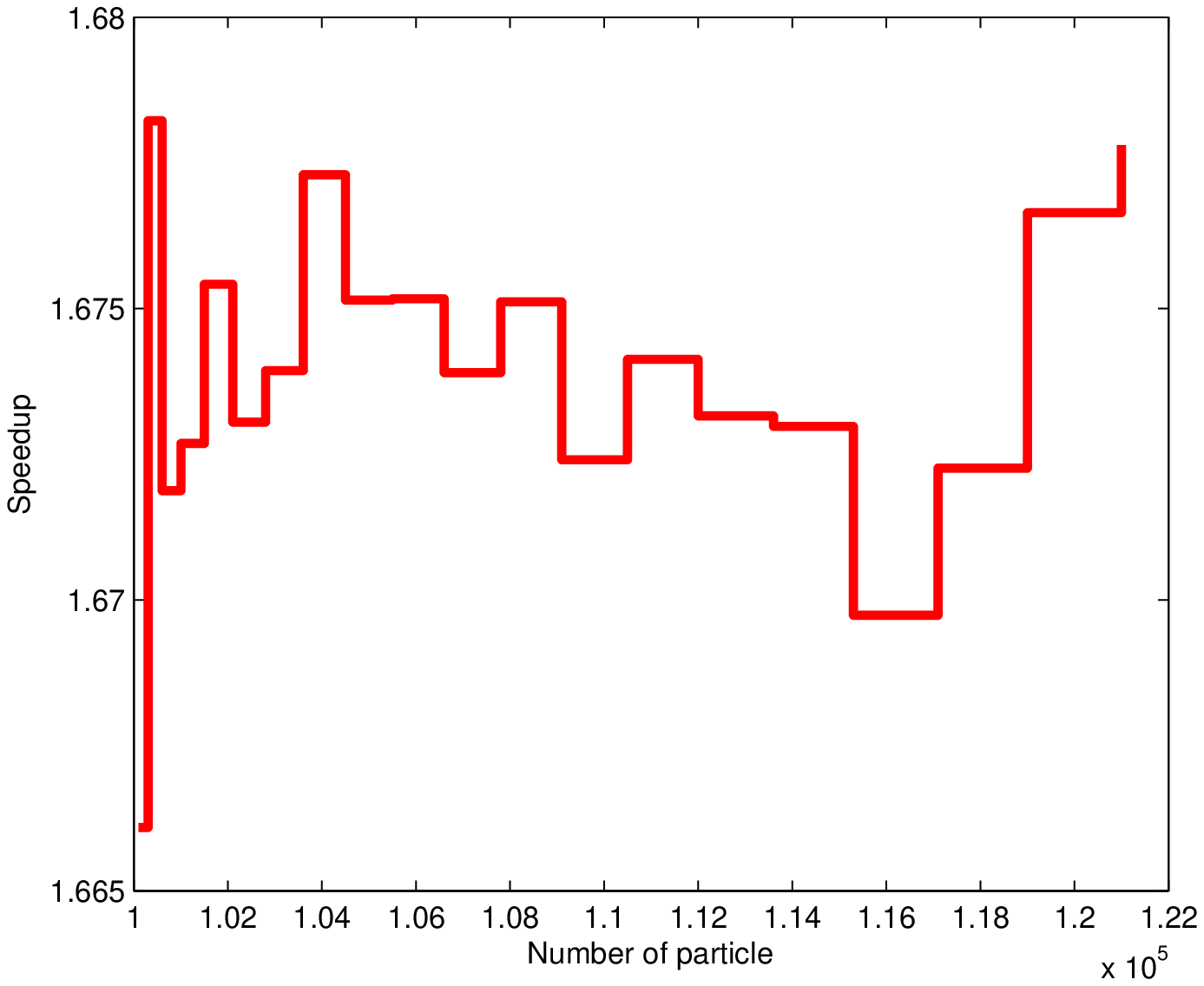}
}
\subfigure[PSO saved time]{
\includegraphics[height=0.35\textwidth, width=0.45\textwidth]{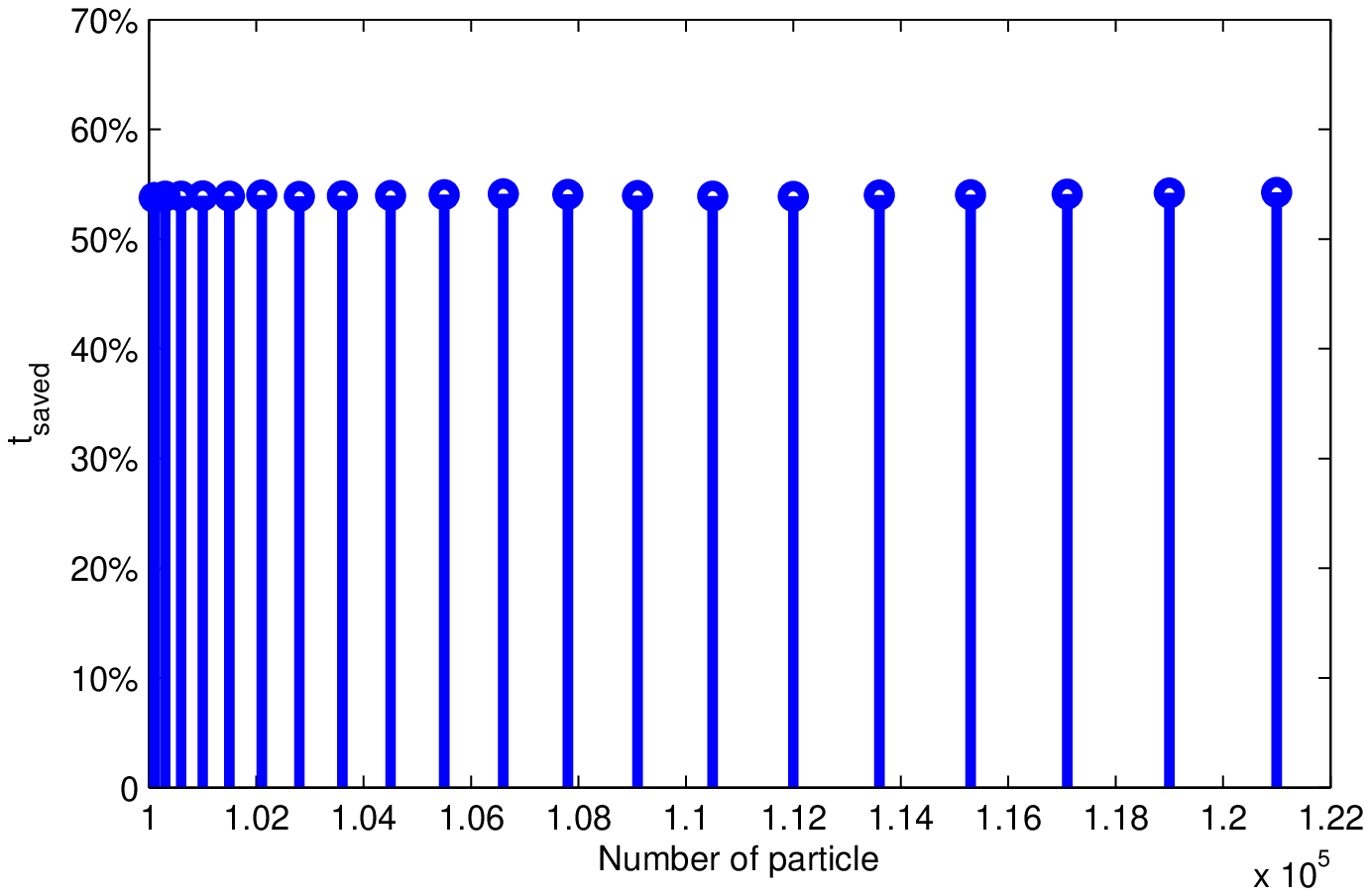}
 }
 \subfigure[PSO speedup]{
\includegraphics[height=0.35\textwidth, width=0.45\textwidth]{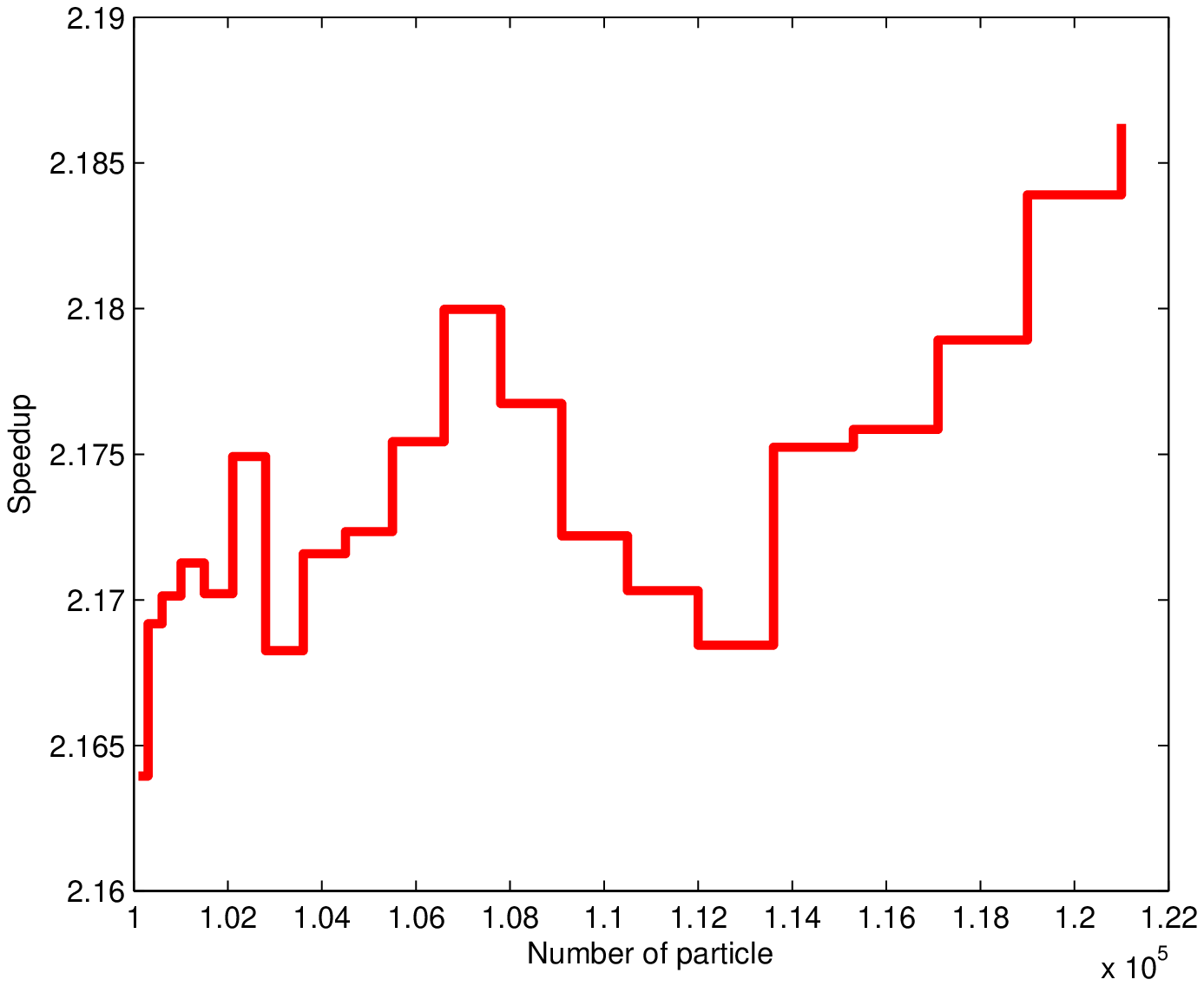}
}
\caption[]{Parallel MCO and PSO solving Ackley function.}
\label{fAckley}
\end{figure}

\begin{figure} \centering
\subfigure[MCO saved time]{
\includegraphics[height=0.35\textwidth, width=0.45\textwidth]{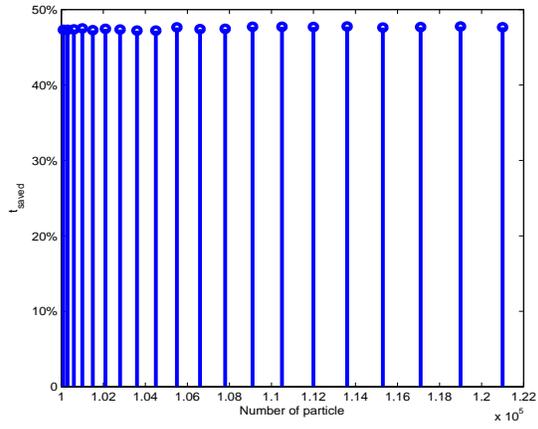}
 }
 \subfigure[MCO speedup]{
\includegraphics[height=0.35\textwidth, width=0.45\textwidth]{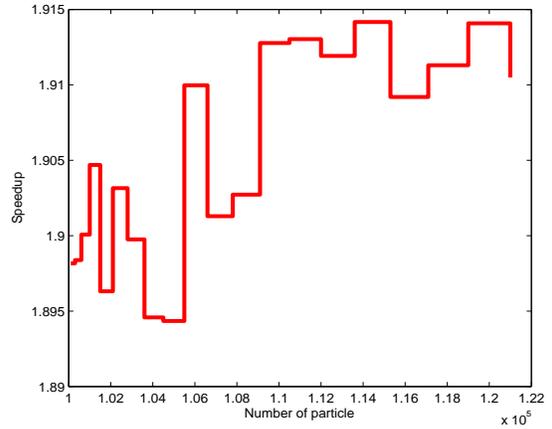}
}
\subfigure[PSO saved time]{
\includegraphics[height=0.35\textwidth, width=0.45\textwidth]{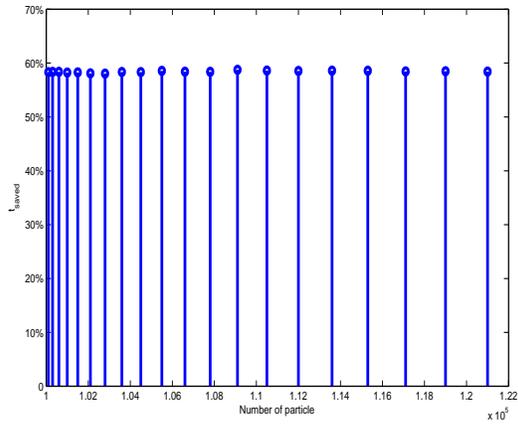}
 }
 \subfigure[PSO speedup]{
\includegraphics[height=0.35\textwidth, width=0.45\textwidth]{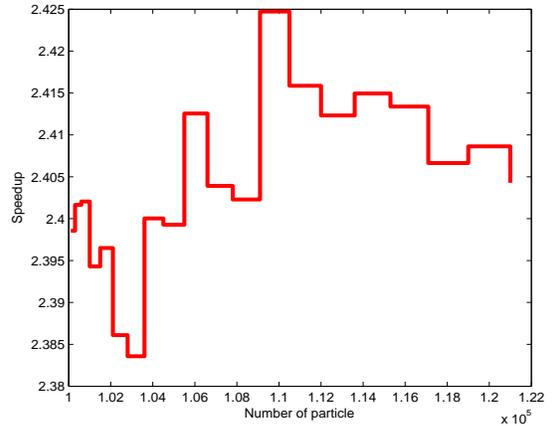}
}
\caption[]{Parallel MCO and PSO solving De Jong's f4 function.}
\label{fDe Jong's f4}
\end{figure}

\begin{figure} \centering
\subfigure[MCO saved time]{
\includegraphics[height=0.35\textwidth, width=0.45\textwidth]{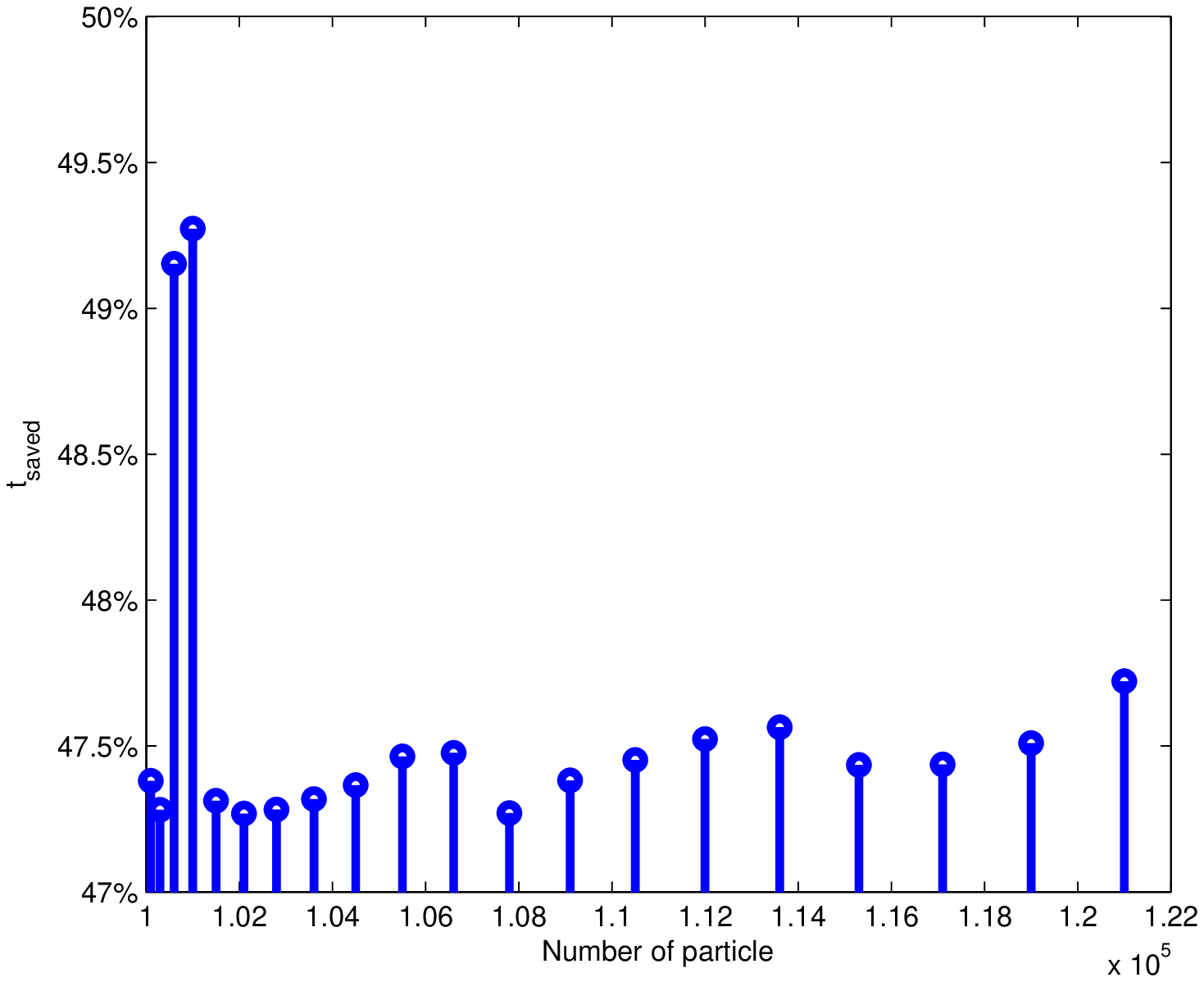}
 }
 \subfigure[MCO speedup]{
\includegraphics[height=0.35\textwidth, width=0.45\textwidth]{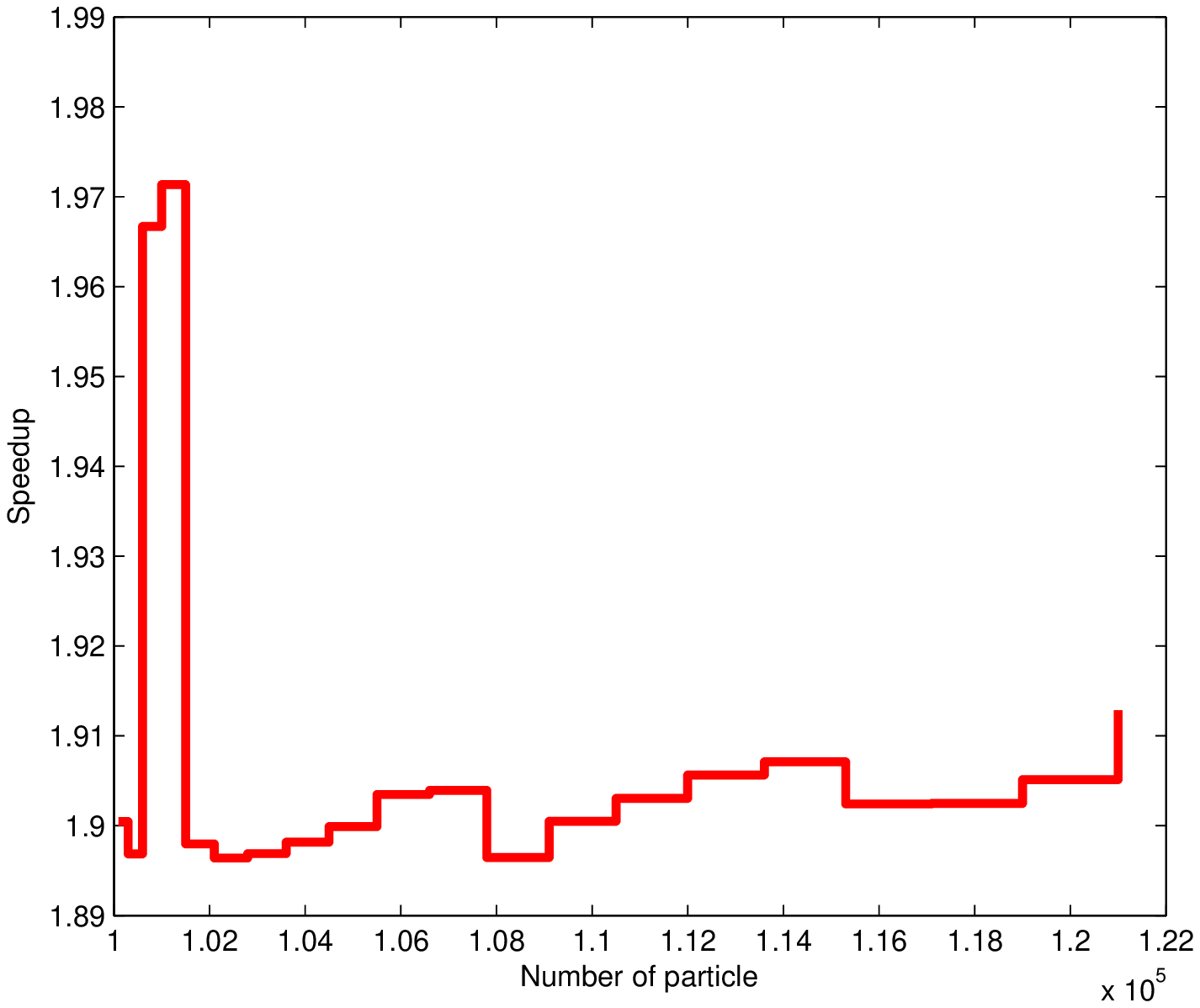}
}
\subfigure[PSO saved time]{
\includegraphics[height=0.35\textwidth, width=0.45\textwidth]{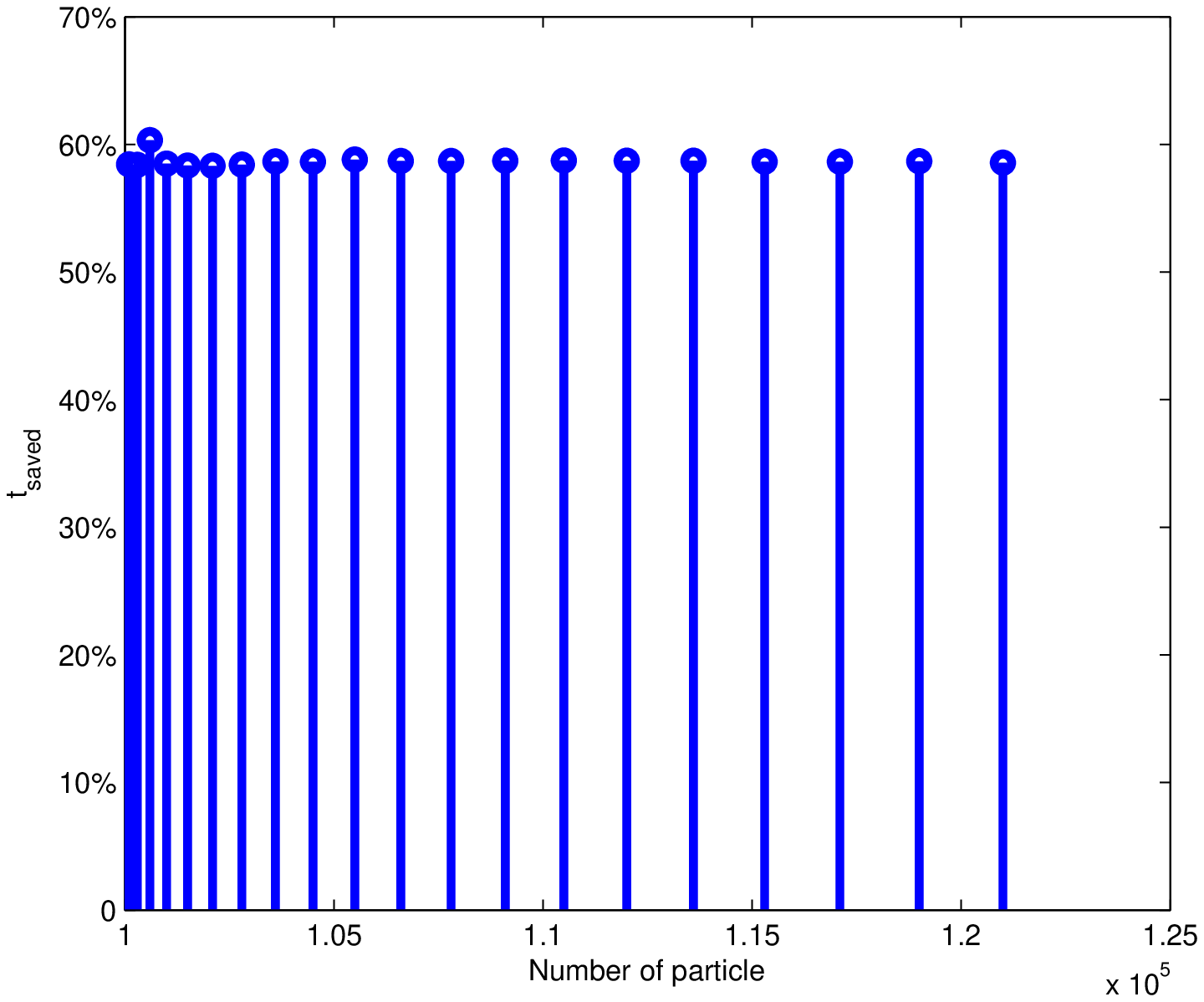}
 }
 \subfigure[PSO speedup]{
\includegraphics[height=0.35\textwidth, width=0.45\textwidth]{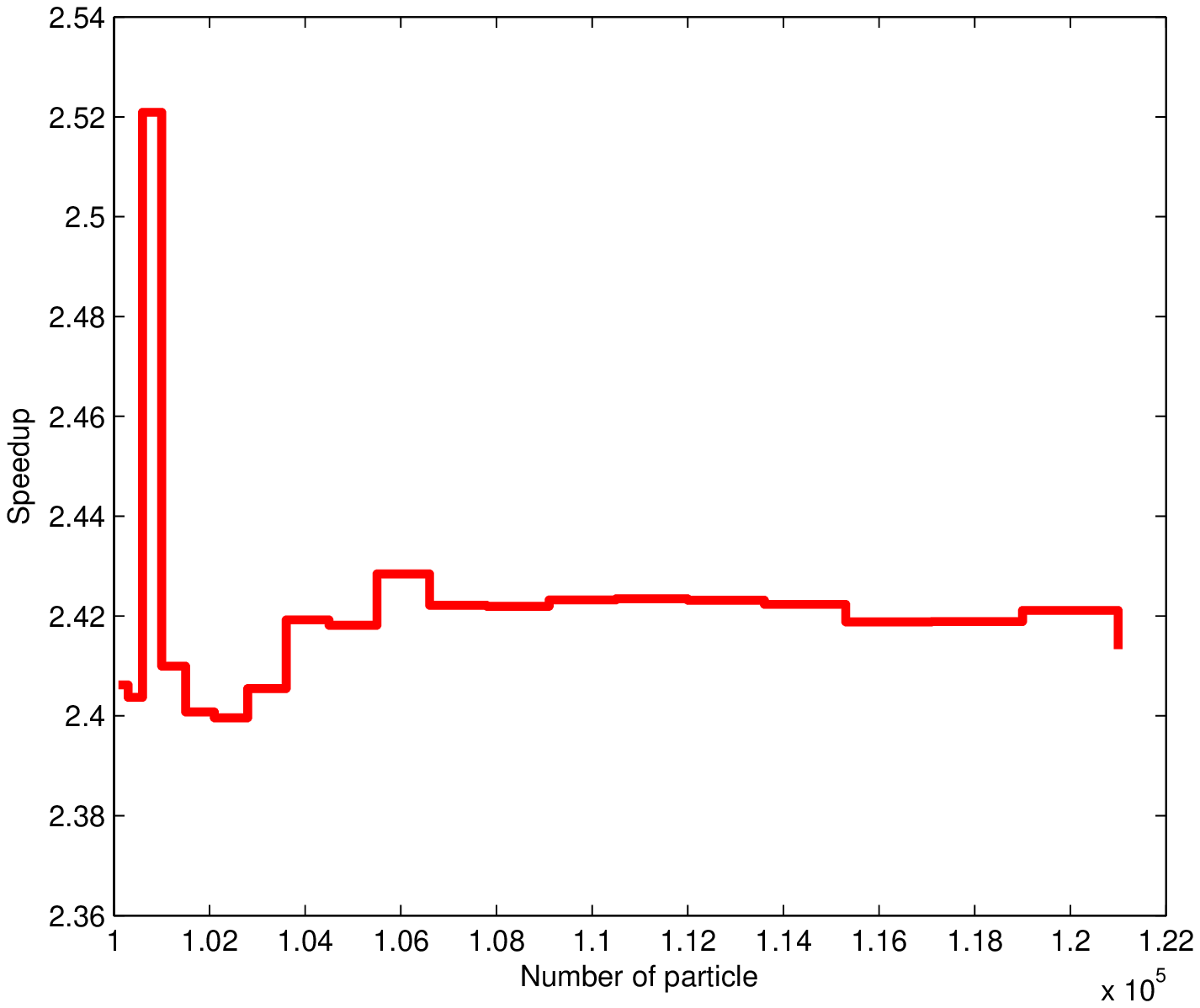}
}
\caption[]{Parallel MCO and PSO solving Levy function.}
\label{fLevy}
\end{figure}

\begin{figure} \centering
\subfigure[MCO saved time]{
\includegraphics[height=0.35\textwidth, width=0.45\textwidth]{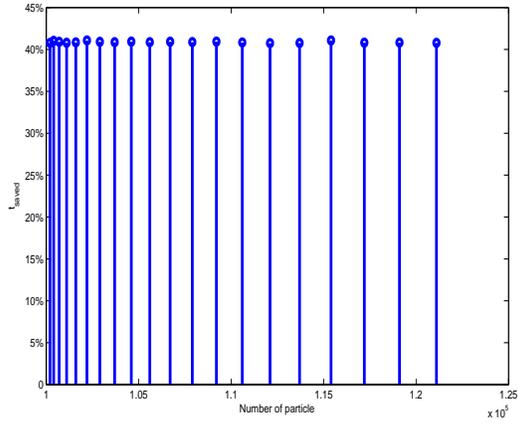}
 }
 \subfigure[MCO speedup]{
\includegraphics[height=0.35\textwidth, width=0.45\textwidth]{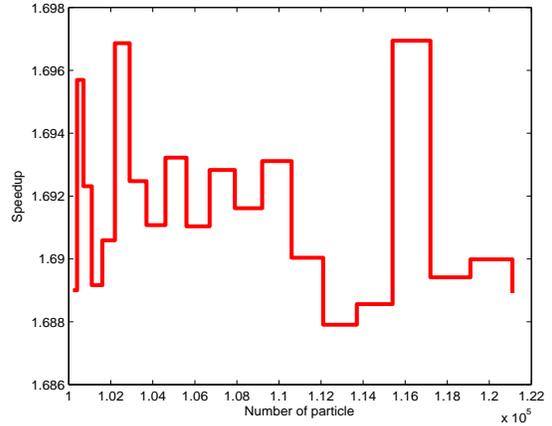}
}
\subfigure[PSO saved time]{
\includegraphics[height=0.35\textwidth, width=0.45\textwidth]{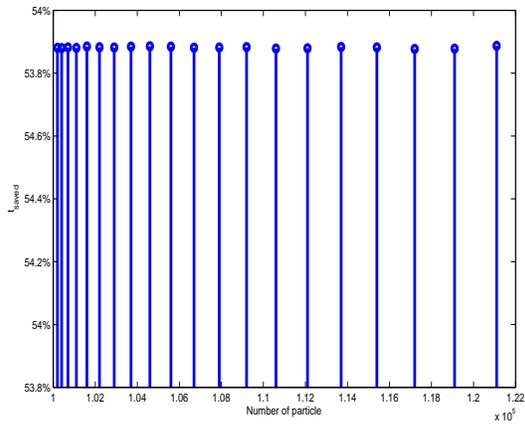}
 }
 \subfigure[PSO speedup]{
\includegraphics[height=0.35\textwidth, width=0.45\textwidth]{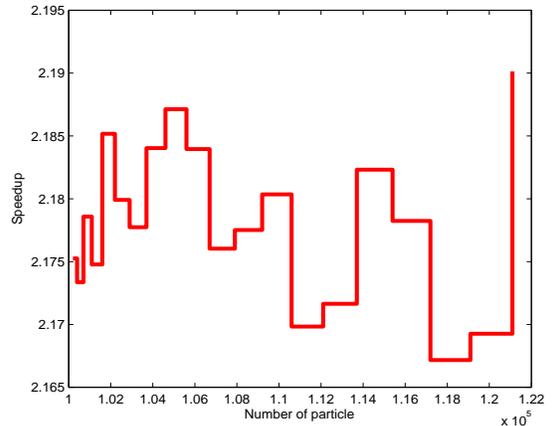}
}
\caption[]{Parallel MCO and PSO solving Zakharov function.}
\label{fZak}
\end{figure}

\begin{figure} \centering
\subfigure[Scalability for parallel MCO]{
\includegraphics[height=0.35\textwidth, width=0.45\textwidth]{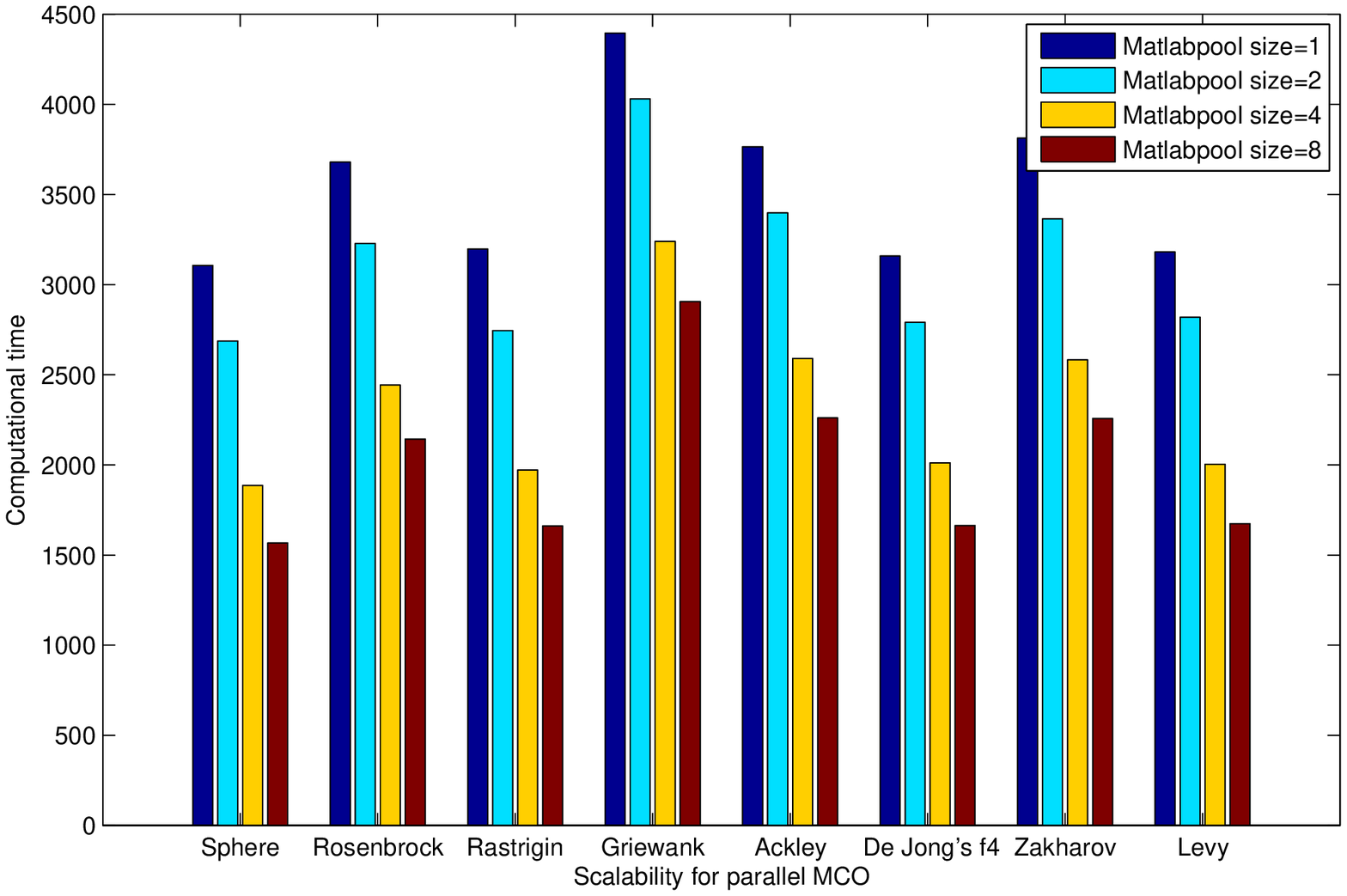}
 }
 \subfigure[Scalability for parallel PSO]{
\includegraphics[height=0.35\textwidth, width=0.45\textwidth]{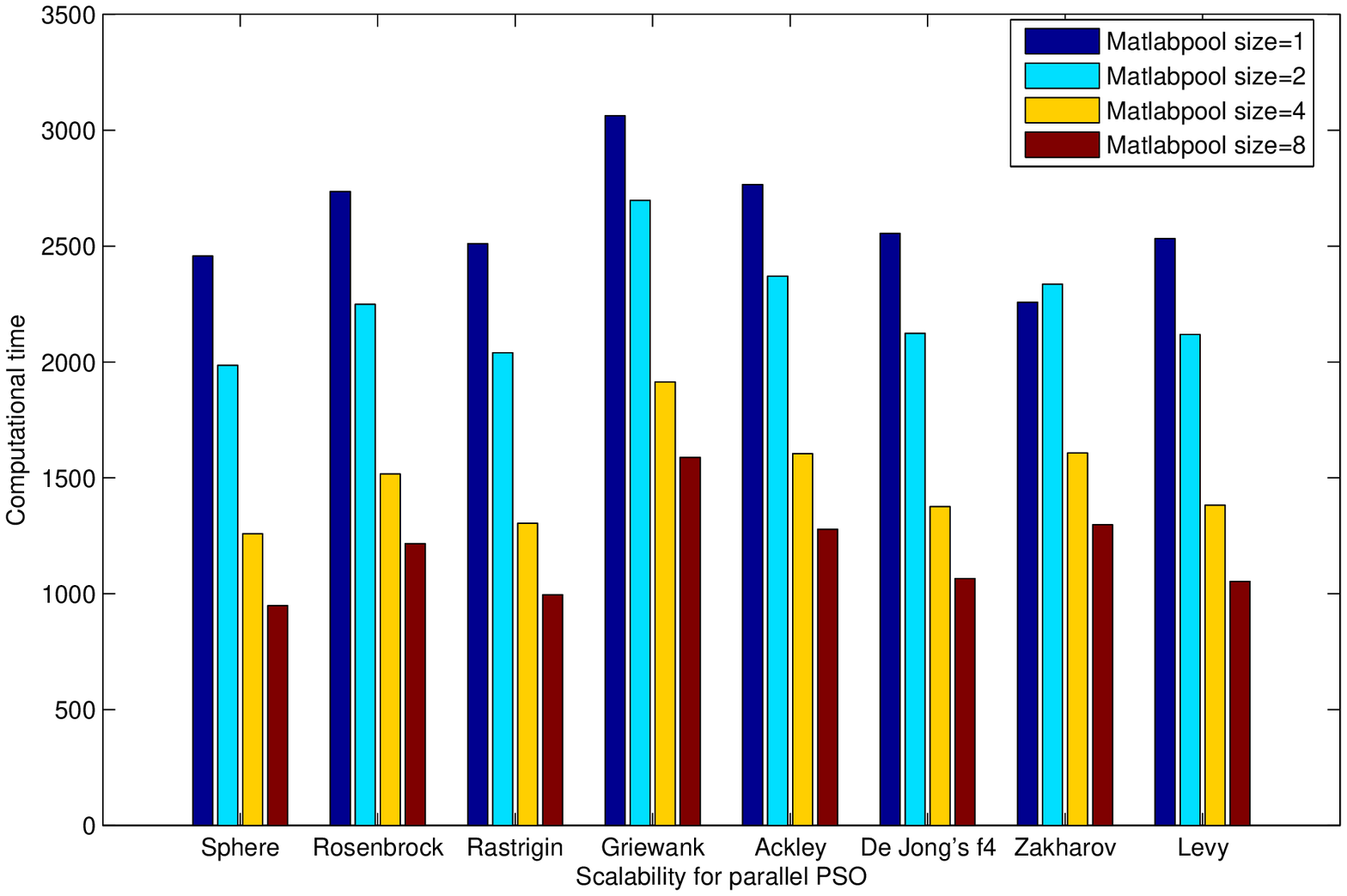}
}
\caption[]{Scalability chart.}
\label{fScalability}
\end{figure}

We first evaluate the computational time of the parallel MCO for different test functions. Specifically, eight 2.8 GHz cores equipped supercomputers in the High Performance Computing Center at Texas Tech University were used to run the parallel MCO and PSO algorithms for all the eight benchmark functions in which the search areas and dimensions of objective functions are listed in Subsection \ref{tf} with $n=30$. The $\mathtt{matlabpool}$ size is 8. We choose the communication graph $\mathcal{G}_{k}$ for MCO to be a complete graph. The simulation results are shown in Fig. \ref{fsphere}--\ref{fZak} and the scalability chart is also provided for different $\mathtt{matlabpool}$ sizes in Fig. \ref{fScalability}. The saving time $t_{saved}$ is calculated as
$t_{saved}=(t_{seri}-t_{para})/t_{seri}\times 100 \%$, and $speedup=\frac{t_{seri}}{t_{para}}$,
where $t_{seri}$ and $t_{para}$ are the computational times for the serial algorithm and parallel algorithm when solving the same optimization problem, respectively. From the simulation results, the parallel MCO algorithm can shorten the computational time by about $40\%$ to $70\%$ compared with the serial MCO. In the meantime, the parallel PSO algorithm also showed a similar improvement in computation time. Although the MCO algorithm conducts more computation than PSO, the accuracy of MCO is superb compared with PSO due to the information exchange between neighbors.

\subsection{Evaluation of Numerical Accuracy for the Parallel MCO}

To evaluate numerical accuracy for the parallel MCO, the statistical results of the best values obtained from the standard PSO, serial MCO and parallel MCO algorithms are compared numerically. Similarly, the search areas and dimensions of objective functions are listed in Subsection \ref{tf} with $n=30$. The maximum of the objective values, the minimum of the
objective values, the average of objective value, and the median
objective values are compared in Table \ref{table2}. Based on these results, it follows that the serial MCO and parallel MCO algorithms are more accurate for obtaining the best values than the PSO algorithm.

\begin{table*}[t]\tiny
\caption{Numerical Comparison Between serial and parallel PSO, and serial and parallel MCO for the Eight Test Functions} \label{table2}
\begin{center}
\begin{tabular}{|c |c| c|c | c| c|}
\hline
Function   &Min  &Max &Median  &Average  \\
\hline
& Ser. PSO  Par. PSO  Ser. MCO  Par. MCO &  Ser. PSO Par. PSO  Ser. MCO  Par. MCO&  Ser. PSO Par. PSO Ser. MCO Par. MCO& Ser. PSO  Par. PSO Ser. MCO  Par. MCO \\
\hline Sphere &2.3E-3 3.3E-3 4.4499E-7 4.1855E-7 &1.02E-2 9.8E-3 1.9988E-6 2.0584E-6  &6.5E-3 6.3E-3 8.7833E-7 1.1932E-6  &6.1E-3 6.4E-3 9.7107E-7 1.1520E-6\\
\hline Rosenbrock &4.0605E1 4.554E1 2.885E1 1.5600E-4 & 1.0455E2 9.741E1 8.059E1 9.155E1& 6.533E1 6.668E1  5.472E1 5.621E1 &6.973E1 6.668E1 5.413E1 5.826E1\\
\hline Rastrigin & 2.722E1 2.181E1 1.762E1 2.108E1 &9.266E1 7.976E1 7.962E1 7.748E1 &2.967E1 2.999E1 2.802E1 2.658E1 &4.288E1 3.143E1 3.87E1 3.364E1\\
\hline Griewank & 2.3E-3  3.2E-3 4.4499E-7  3.3862E-7&1.02E-2 1.03E-2 1.9988E-6 1.8827E-6&6.5E-3 5.6E-3 8.7833E-7 9.3364E-7 &6.0E-3 6.1E-3  9.7107E-7 9.5507E-7 \\
\hline Ackley &1.047E-1  9.454E-1 7.5209E-4 8.0127E-4 &2.6657 2.3269 2.0133 2.4078 &1.6489 1.8530 0.6711 1.9E-3&1.5231  1.7364 8.435E-1  7.522E-1\\
\hline De Jong's f4 &2.5018 1.2998 1.5389E-5 3.7159E-5&3.275E1 8.907E1 7.7E-3 4.9E-2&1.097E1  1.501E1  1.9E-3 1.097E1&1.230E1 1.9995E1 1.5E-3 1.9E-3 \\
\hline Zakharov &1.9932  1.3981 5.88E-2  1.34E-2 & 1.755E1 9.3668 3.283E-1  4.343E-1& 4.7803 5.5631 1.355E-1  1.02E-1 &4.9881  6.3180 1.576E-1  1.337E-1\\
\hline Levy &5.5330 5.1027 3.8140 4.9910 &3.0949E1 4.8200E1 3.7046E1  2.8139E1& 1.1543E1 1.2932E1 1.3546E1 1.6651E1 & 1.6841E1 1.3974E1 1.6066E1 1.6545E1\\
\hline
\end{tabular}
\end{center}
\end{table*}

\section{Conclusion}\label{co}
In this report, a parallel MCO algorithm is developed by introducing the MATLAB built-in parallel function $\mathtt{parfor}$ into the inner loop of the MCO algorithm. The numerical evaluation concludes that the parallel MCO algorithm can achieve similar accuracy compared with the serial MCO algorithm, but in a shorter computational time. A detailed convergence analysis of the MCO algorithm is presented. Future work will focus on the large-scale, real-time engineering applications of this parallel MCO algorithm, such as power grid network vulnerability problems.





\bibliographystyle{IEEEtran}
\bibliography{Reference}

\end{document}